\documentclass[oneside, a4paper]{amsart}

	\usepackage[a4paper, margin=1in]{geometry}
	\usepackage{amsmath, amsthm, amssymb, latexsym, mathtools, mathrsfs, eucal, mathscinet}
	\usepackage[all,cmtip]{xy}
	\usepackage{dsfont, soul, stackrel, tikz-cd, graphicx, etoolbox, multirow}
	\usepackage{stmaryrd}
	\DeclareGraphicsExtensions{.pdf,.png,.jpg}
	\usepackage{fancyhdr}
	\usepackage[shortlabels]{enumitem}
	\usepackage[pdfmenubar=true, pdfborder={0 0 0 [3 3]}, colorlinks, citecolor=magenta, urlcolor=blue]{hyperref}
    \usepackage{calrsfs}
    \usepackage{comment}	
\DeclareMathAlphabet{\pazocal}{OMS}{zplm}{m}{n}
    \usepackage{mathdots}
	\numberwithin{equation}{section}
	\setenumerate{listparindent=\parindent}
 \usepackage{booktabs}
 \usepackage{tabularx}
 \usepackage{makecell} 

	\newtheoremstyle{Mytheorem}%
	{1em}{1em}%
	{\slshape}{}%
	{\bfseries}{.}%
	{ }{}

	\newtheoremstyle{Mydefinition}%
	{1em}{1em}%
	{}{}%
	{\bfseries}{.}%
	{ }{}

	\theoremstyle{Mydefinition}
	\newtheorem{statement}{Statement}[section]
	\newtheorem{definition}[statement]{Definition}
	
	\newtheorem{remark}[statement]{Remark}
	
	\newtheorem{example}[statement]{Example}

	\newtheorem*{comment*}{Comment}

	\theoremstyle{Mytheorem}
	\newtheorem{theorem}[statement]{Theorem}
	\newtheorem{corollary}[statement]{Corollary}
	
	\newtheorem{proposition}[statement]{Proposition}
	\newtheorem{lemma}[statement]{Lemma}

	\newcommand{\nc}{\newcommand}
	\newcommand{\be}{\begin{eqnarray*}}
	\newcommand{\ee}{\end{eqnarray*}}
	\newcommand{\bea}{\begin{eqnarray}}
	\newcommand{\eea}{\end{eqnarray}}
	\newcommand{\bs}{\begin{split}}
	\newcommand{\es}{\end{split}}
	\newcommand{\bal}{\begin{align}}
	\newcommand{\eal}{\end{align}}

	\newcommand{\inj}{\mathrm{inj}}

	\newcommand{\red}{\mathrm{red}}

	\nc{\bei}{9}
	\nc{\eei}{\end{itemize}}
	\nc{\bee}{\begin{enumerate}}
	\nc{\eee}{\end{enumerate}}
	\nc{\bet}{\begin{thm}}
	\nc{\eet}{\end{thm}}
	\nc{\bed}{\begin{defn}}
	\nc{\eed}{\end{defn}}
	\nc{\bel}{\begin{lem}}
	\nc{\eel}{\end{lem}}
	\nc{\bep}{\begin{prop}}
	\nc{\eep}{\end{prop}}
	\nc{\bec}{\begin{corollary}}
	\nc{\eec}{\end{corollary}}
	\nc{\ber}{\begin{rem}}
	\nc{\eer}{\end{rem}}
	\nc{\beex}{\begin{example}}
	\nc{\eeex}{\end{example}}
	\nc{\bpm}{\begin{pmatrix}}
	\nc{\epm}{\end{pmatrix}}
	\nc{\bspm}{\left(\begin{smallmatrix}}
	\nc{\espm}{\end{smallmatrix}\right)}

	\newcommand{\cA}{\mathcal{A}}

	\newcommand{\pcC}{\pazocal{C}}

	\newcommand{\cF}{\mathcal{F}}
	\newcommand{\cG}{\mathcal{G}}
	\newcommand{\cH}{\pazocal{H}}

	\newcommand{\cK}{\mathcal{K}}
	\newcommand{\cL}{\mathcal{L}}
	\newcommand{\cM}{\mathcal{M}}
	
	\newcommand{\cO}{\mathcal{O}}

	\newcommand{\cU}{\mathcal{U}}
	\newcommand{\pcU}{\pazocal{U}}
	\newcommand{\cV}{\mathcal{V}}
	\newcommand{\pcV}{\pazocal{V}}
	
	\newcommand{\cX}{\pazocal{X}}

	\newcommand{\bA}{\mathbb{A}}
	
	\newcommand{\bC}{\mathbb{C}}

	\newcommand{\bP}{\mathbb{P}}
	
	\newcommand{\bR}{\mathbb{R}}

	\newcommand{\bZ}{\mathbb{Z}}

	\nc{\frf}{\mathfrak{f}}

	\nc{\frs}{\mathfrak{s}}  
	\nc{\frt}{\mathfrak{t}} 
	\nc{\fru}{\mathfrak{u}}
	\nc{\lsl}{\mathfrak{sl}}
	\nc{\lgl}{\mathfrak{gl}}
	\nc{\upsi}{\underline{\psi}}
	\nc{\uchi}{\underline{\chi}}
	 
	\DeclareMathOperator{\PV}{PV}
	
	\DeclareMathOperator{\Hom}{Hom}

	\DeclareMathOperator{\Lie}{Lie}

	\DeclareMathOperator{\id}{id}

	\DeclareMathOperator{\Ad}{Ad}

	\DeclareMathOperator{\Cl}{Cl}
	\DeclareMathOperator{\sgn}{sgn}

	\DeclareMathOperator{\dist}{dist}

	\newcommand{\lra}{\longrightarrow}    
	
	\nc{\surjto}{\twoheadrightarrow}
	\nc{\ts}{\times}
	\nc{\ds}{\displaystyle}
	\nc{\nd}{\noindent}  
	\newcommand{\ud}{\underline}
	\nc{\ov}{\overline}
	\nc{\maplra}[1]{\buildrel #1 \over \lra}
	\nc{\mapto}[1]{\buildrel #1 \over \to}
	\nc{\setb}[1]{\{  #1\}}
	\nc{\cHom}{\mathcal{H}om}

	\def\a{\alpha}
	\def\b{\beta}
	
	\def\d{\delta} 
	\def\e{\varepsilon} 
	 
	\def\g{\gamma} 
	\def\k{\kappa}

	\def\n{\nu}
	\def\r{\rho}
	
	 \def\O{\Omega}
	 
	\def\s{\sigma} 
	\def\th{\theta}

	\def\x{\xi}

	\def\C{\mathbb{C}}
	
	\def\Z{\mathbb{Z}}

	\def\Ker{\hbox{Ker}\;}

	\def\wt{\operatorname{wt}}
	\def\ch{\operatorname{ch}}

\title[$L^2$-Hodge theory of Hybrid LG/CY complete intersections]{Bounded Geometries on Hybrid Landau-Ginzburg models of Calabi-Yau complete intersections and $L^2$-Hodge Theory}

\author{Jeehoon Park and Jaewon Yoo}

\address{Jeehoon Park: QSMS, Seoul National University, 1 Gwanak-ro, Gwanak-gu, Seoul 08826, South Korea}
\email{jpark.math@gmail.com}

\address{Jaewon Yoo: Department of Mathematics, POSTECH, San 31, Hyoja-Dong, Nam-Gu, Pohang, Gyeongbuk, 790-784, South Korea}
\email{yooj1215@postech.ac.kr}


\begin{document}

\begin{abstract}
Given a Calabi-Yau smooth projective complete intersection variety $V$ over $\C$, a hybrid Landau-Ginzburg (LG) model may be associated using the Cayley trick. 
This hybrid LG model comprises a non-compact Calabi-Yau manifold $X_{CY}$, and a holomorphic function $W$, defined on $X_{CY}$, such that the critical locus of $W$ is isomorphic to $V$.
We construct a complete K\"ahler metric $\mathfrak{g}$  and a bounded Calabi-Yau volume form ${\Omega}$ on $X_{CY}$ such that $(X_{CY},\mathfrak{g}, {\Omega})$ is a bounded Calabi-Yau geometry (in fact, $(X_{CY},\mathfrak{g})$ is an asymptotically conical manifold) and the function $W$ is strongly elliptic; this enables us to apply the $L^2$-Hodge theory of Li-Wen \cite{LW} to $(X_{CY},\mathfrak{g}, {\Omega})$ and $W$, which leads to 
a Frobenius manifold structure on the twisted de Rham cohomology associated to $(X_{CY},W)$.
Furthermore, we prove that this twisted de Rham cohomology is isomorphic to the de Rham cohomology $H(V;\C)$, which results in a new $L^2$-Hodge theoretic construction of a Frobenius manifold structure on $H(V;\C)$.
{This paper provides the first explicit geometric verification of Li-Wen's theory for genuine non-isolated, compact critical loci using hybrid Landau-Ginzburg models.}
\end{abstract}

\maketitle

\tableofcontents

\section{Introduction}
The notion of a Frobenius manifold refers to a family of Frobenius algebras,
and its axiomatic treatment was first given by B. A. Dubrovin in his study of 2 dimensional topological field theory \cite{Dubrovin}. 
Even before Dubrovin's work \cite{Dubrovin}, K. Saito provided examples of Frobenius manifolds by analyzing period mappings of primitive forms for germs of holomorphic functions with isolated singularities \cite{Saito}. Sergey Barannikov and Maxim Kontsevich constructed Frobenius manifold structures on the cohomology of compact Calabi-Yau manifolds \cite{BK}. For a comprehensive understanding of Frobenius manifolds, we refer to the excellent books \cite{Her02} and \cite{Sab}.  

The Frobenius manifold structure plays an important role in singularity theory. Notably, the universal unfolding of a holomorphic function with isolated singularities gives rise to a Frobenius manifold structure. In general, the construction of Frobenius manifold structures has been limited to cases when the function $f$ has an isolated singularity. 
However, the works \cite{Sa} and \cite{OV} explored the case of non-isolated critical loci and \cite{OV} uses characteristic $p$ methods. More recently, Si Li and Hao Wen \cite{LW} introduced $L^2$-Hodge theoretic techniques to construct a Frobenius manifold for holomorphic functions with compact critical loci, drawing parallels to the work of Barannikov and Kontsevich \cite{BK}. Although these approaches share common ground in some cases, they also address different scenarios.

In Li-Wen's theory \cite{LW}, two important technical conditions must be satisfied for the construction of Frobenius manifolds.
Namely, when $X$ is a non-compact Calabi-Yau manifold with a non-vanishing holomorphic volume form ${\Omega}$ and $f:X\to \bC$ is a holomorphic function with compact critical locus, the following two technical conditions are required:
\begin{enumerate}
\item $X$ is a bounded Calabi-Yau geometry (Definition \ref{Definition:bcyg}),
\item $f:X\to \bC$ is a strongly elliptic holomorphic function (Definition \ref{Definition:strell}). 
\end{enumerate}
 
 Under these conditions, they constructed a commutative differential graded algebra $(PV^\bullet(X),\overline{\partial}_f)$ where $PV^\bullet(X_{CY})$ denotes the space of holomorphic polyvector fields valued in anti-holomorphic differential forms. They introduced a twisted differential $\overline{\partial}_f=\overline{\partial}+\{f,-\}$, and demonstrated that the cohomology $H(PV^\bullet(X_{CY}),\overline{\partial}_f)$ carries a Frobenius manifold structure (Theorem \ref{LWmain}). They gave three classes of examples of $(X,{\Omega},f)$ satisfying the two conditions. 

\begin{example}[Section 2.4 of \cite{LW}]
\ 
\label{LWexample}
\begin{enumerate}
\item Let $(X,{\Omega})$ be $(\bC^n,dz_1\wedge \cdots\wedge dz_n)$, $g$ be the standard Euclidean metric, and $f$ be a non-degenerate quasi-homogeneous polynomial. Then $(X,{\Omega},f)$ satisfies the conditions.
\item Let $(X,{\Omega})$ be $((\bC^*)^n, \frac{dz_1}{z_1}\wedge\cdots\wedge \frac{dz_n}{z_n})$, $g$ be the standard complete metric, and $f$ be a convenient nondegenerate Laurent polynomial. Then $(X,{\Omega},f)$ satisfies the conditions.
\item Let $\pi:X\to \bC^n/G$ be a crepant resolution of the quotient of $\bC^n$ by a finite group $G\subset SU(n)$, ${\Omega}$ be the pullback of the trivial Calabi-Yau form on $\bC^n$. There exists an asymptotically locally Euclidean K\"ahler metric $g$ on $X$, constructed by D. Joyce \cite{Joyce}. Let $f:\bC^n\to \bC$ be a $G$-invariant polynomial with an isolated singularity at the origin. Then $(X,{\Omega},\pi^*f)$ satisfies these conditions.
\end{enumerate}
\end{example}

Examples (1) and (2) are relatively well understood. Example (3) is a genuinely new case suggested by Li and Wen. Although $\pi^* f$ lacks an isolated singularity, McKay's correspondence \cite{Velez} suggests that $(X,\pi^* f)$ is equivalent to the orbifold Landau-Ginzburg B-model $(\bC^n,f,G)$. Thus, case (3) effectively reduces to the scenario where $f$ has an isolated singularity. 

This paper aims to provide concrete examples of \textit{genuine} non-isolated singularities that satisfy two technical conditions.
Our idea is to use the hybrid LG (Landau-Ginzburg) model of a smooth projective CY (Calabi-Yau) complete intersection variety, which was studied by Alessandro Chiodo and Jan Nagel in \cite{CN}. When $V(\ud{W}):=V(W_1,\ldots,W_r)$ is a Calabi-Yau smooth complete intersection in the weighted projective space $\bP(\ud{w}):=\bP(w_1,\ldots,w_n)$ with $d_i:=\deg(W_i)$, the Cayley trick leads them to define a holomorphic function $W=\sum_{j=1}^r p_r W_r$ on the total space of the locally free sheaf $\cO_{\bP(\ud{w})}(-\ud{d}):=\cO_{\bP(\ud{w})}(-d_1)\oplus\cdots \oplus \cO_{\bP(\ud{w})}(-d_r)$, and they established the cohomological LG/CY correspondence on the Chen-Ruan cohomologies $H_{CR}(V(\ud{W});\bC)$ and $H_{CR}(\cO_{\bP(\ud{w})}(-\ud{d}),W^{-1}(t_0);\bC)$, under the Calabi-Yau condition $d_1+\cdots+d_r=w_1+\cdots+w_n$.

We observe that $X_{CY}:=\mathrm{Tot}(\cO_{\bP(\ud{w})}(-\ud{d}))$ is a non-compact Calabi-Yau manifold and
$W:X_{CY}\to \C$ has a compact critical locus $V(\ud{W})$, which is embedded inside $X_{CY}$ as the zero section. A natural question is whether there exists a K\"ahler metric $\mathfrak{g}$ of bounded Calabi-Yau geometry on $X_{CY}$ that enables us to utilize the $L^2$-Hodge theory developed by Li and Wen \cite{LW}. Now we present our main theorem.

\begin{theorem}
\label{Main1}
Suppose $d_1,\ldots,d_r$ are positive integers such that $d_1+\cdots+d_r=n$ and $\max_j d_j \leq 2\min_j d_j$. Let $W_1,\ldots,W_r$ be homogeneous polynomials of degrees $d_1,\ldots,d_r$ such that $V(\ud{W})$ is a Calabi-Yau smooth complete intersection in $\bP^{n-1}$. Then there is a complete K\"ahler metric $\mathfrak{g}$ and a nonvanishing bounded holomorphic volume form ${\Omega}$ so that $(X_{CY},\mathfrak{g},{\Omega})$ is a bounded Calabi-Yau geometry, and $W=p_1 W_1 +\cdots + p_r W_r$ is strongly elliptic on $(X_{CY},\mathfrak{g})$. Consequently, there is a Frobenius manifold structure on $H(PV(X_{CY}),\overline{\partial}_W)$.
\end{theorem}

To prove Theorem~\ref{Main1}, we exploit the toric structure of \(X_{\mathrm{CY}}\) and write it in its standard quotient form. This quotient description allows us to construct K\"ahler metrics on \(X_{\mathrm{CY}}\) by K\"ahler reduction from an explicit \(U(1)\)-invariant K\"ahler metric on the ambient space. The key point is that the naive Euclidean choice on the ambient space, corresponding to \(s=0\) in our family, does not make the natural holomorphic volume form and its dual polyvector field simultaneously bounded after reduction. Our solution is therefore to replace the Euclidean metric by a one-parameter family of conic ambient metrics. The parameter \(s\) controls the asymptotic size of these tensors, and the balanced value \(s=\frac{1}{n+r-1}\) is exactly the one for which the reduced geometry becomes a bounded Calabi--Yau geometry.

In the compact case, Yau's solution of the Calabi conjecture provides a Ricci-flat K\"ahler metric in every K\"ahler class when \(c_1=0\) \cite{Yau} \cite{Yau2}. The noncompact situation is much subtler. At present, there is no general existence theorem in the noncompact setting comparable to the compact case, and complete Ricci-flat K\"ahler metrics are known only under additional geometric assumptions. Classical examples arise in the work of Tian--Yau on quasi-projective manifolds \cite{TY90} \cite{TY91}, Joyce on ALE crepant resolutions \cite{Joyce}, and later work on asymptotically conical metrics by van Coevering \cite{Coevering} and Conlon--Hein \cite{ConHein2} \cite{ConHein}. For this reason, even when one is ultimately interested in Ricci-flat geometry, it is natural to begin with other robust classes of complete K\"ahler metrics that can be constructed and analyzed directly.

One particularly important class is given by asymptotically conical K\"ahler metrics. Roughly speaking, a complete K\"ahler metric is asymptotically conical if, outside a compact set, it approaches a K\"ahler cone with quantitative decay of all derivatives. In the Ricci-flat setting, one further asks that the complex structure and holomorphic volume form also converge to those of a Calabi--Yau cone; this is the standard notion of an asymptotically conical Calabi--Yau metric. In our setting, the metric constructed on \(X_{\mathrm{CY}}\) is asymptotically conical in this K\"ahler sense, and this provides the bounded-geometric framework needed for the \(L^2\)-Hodge-theoretic arguments that follow. The bounded Calabi--Yau condition is then recovered by identifying the unique parameter \(s\) for which the natural holomorphic volume form and its dual are simultaneously bounded.

The condition $\max_j d_j \leq 2\min_j d_j$ is necessary to prove that $W$ is strongly elliptic; there is a counterexample such that $W$ is not strongly elliptic when $\max_j d_j > 2\min_j d_j$ in the current setting of the K\"ahler manifold $(\cX,\mathfrak{g}_s,J)$. The condition of strong ellipticity is subject to the setting of the K\"ahler metric, so it is plausible that $W$ becomes strongly elliptic when the K\"ahler metric is defined differently. See Ramarks \ref{remark3.40} and \ref{remark3.41} for details.

We provide a detailed computation of the cohomology $H(PV(X_{CY}),\overline{\partial}_W)$ in Theorem \ref{Main1}, and show that this cohomology is isomorphic to that of the smooth complete intersection $V(\ud{W})$.  In Section \ref{section4}, we prove the following theorem.
\begin{theorem}
\label{Main2}
For each $j=-(n-r-1),\ldots,n-r-1$, there is an isomorphism between vector spaces,
\begin{align}
\label{main:iso}
\Phi_j: H^j (PV(X_{CY}),\overline{\partial}_W)\stackrel{\cong}{\to} H^{j+n-r-1}(V(\ud{W});\bC).
\end{align}
Consequently, we have an induced isomorphism
$$
\Phi: H(PV(X_{CY}),\overline{\partial}_W) \stackrel{\cong}{\to} H(V(\ud{W});\bC).
$$
\end{theorem}
Our calculation relies on the computation of the cohomology of the complex $((\Omega^\bullet_{\bC[\ud{x},\ud{p}]/\bC})_{0},dW)$, which is a Koszul-type sequence; this was done by Alan Adolphson and Steven Sperber in \cite{AS}. 
It is highly nontrivial to relate the computations of cohomology in the smooth (analytic) setting to those in the algebraic setting
Although the cohomology group $H^\bullet (PV(X_{CY}),\overline{\partial}_W)$ is defined in the smooth setting, we will show that it is possible to compute it using algebraic methods once we prove that the cohomology group $H^p(H^q(X_{CY},\Omega_{X_{CY}}^\bullet),dW)$ is finite dimensional in the algebraic setting. Theorems \ref{Main1} and \ref{Main2} imply the following corollary.
\begin{corollary}
There is a $L^2$-Hodge theoretic construction of a Frobenius manifold structure on the total cohomology group $H(V(\ud{W});\bC)$ of the smooth Calabi-Yau complete intersection $V(\ud{W})$.
\end{corollary}

There is another known construction of a Frobenius manifold structure on $H(V(\ud{W});\bC)$ due to Barannikov-Kontsevich \cite{BK}.
In the case where $r=1$, we proved the above corollary in \cite{Paper1} based on the method
of the crepant resolution in \cite[Section 2.4.2]{LW}. Moreover, we provided a precise comparison of the two Frobenius algebra structures in the case $r=1$ in \cite{Paper1}; if $r > 1$, then 
a direct comparison result remains elusive due to technical obstructions (see Remark \ref{compare} for more details). 

\vspace{0.5em}

\textbf{Structure of the paper:} 
We briefly explain the structure of the paper.
In Section \ref{section2.1}, we explain the hybrid Landau-Ginzburg model $(X_{CY}, W)$ associated with a given smooth projective complete intersection variety $V(\ud W)$. In Section \ref{section2.2}, we briefly review Li-Wen's $L^2$-Hodge theory \cite{LW} for Landau-Ginzburg models.

In Section 
\ref{section3.0}, we give the definition of a Riemannian cone, a K\"ahler cone, and an asymptotically conical manifold. We prove that an asymptotically conical manifold has bounded geometry. In Sections \ref{section3.1}, \ref{section3.2}, and \ref{section3.2}, we equip the hybrid LG model $(\cX,W):=(X_{CY}, W)$ with a complete K\"ahler metric $\mathfrak{g}$ and a non-vanishing bounded holomorphic volume form ${\Omega}$ so that $(\cX,\mathfrak{g},{\Omega})$ is a bounded Calabi-Yau geometry, and $W:=p_1 W_1 +\cdots + p_r W_r$ is strongly elliptic on $(\cX,\mathfrak{g})$, i.e. we prove Theorem \ref{Main1}. More specifically, in Section \ref{section3.1}, we prove the existence of $\mathfrak{g}$ on $\cX$ such that $(\cX,\mathfrak{g})$ is a bounded geometry. In Section \ref{section3.2}, we show the existence of $\O$ such that $(\cX,\mathfrak{g},{\Omega})$ is a bounded Calabi-Yau geometry. In Section \ref{section3.3}, we prove that $W:\cX\to \C$ is strongly elliptic.

Sections \ref{section4.1}, \ref{section4.2}, and \ref{section4.3} are devoted to the proof of Theorem \ref{Main2}: we compute the cohomology of the commutative differential graded algebra $(\PV(\cX), \bar\partial_W)$ and show that it is isomorphic to the cohomology of the Calabi-Yau projective smooth complete intersection $V(\ud W)$.
To this end, we introduce two spectral sequences $E^{p,q}$ and $\tilde E^{p,q}$, namely, analytic and algebraic spectral sequences. In Section \ref{section4.1}, we compute the $E_2$-page of the algebraic spectral sequence $\tilde E_2^{p,q}$. In Section \ref{section4.2}, we prove that the $E_2$-pages of two spectral sequences are the same. In Section \ref{section4.3}, we compute the higher terms in the spectral sequences to conclude Theorem \ref{Main1}.

Finally, in Section \ref{appendix2}, we give technical details skipped in the proof of Proposition~\ref{ACisbg}, and in Section \ref{section5.2}, we explain a simultaneous weighted blow-up of $\cX=X_{CY}$ and $X_{LG}$, which generalizes the blow-up map $X_{CY}\to X_{LG}$ in the case $r=1$.

\vspace{0.5em}

\textbf{Acknowledgement}:

Jeehoon Park was supported by the National Research Foundation of Korea (NRF-2021R1A2C1006696) and the National Research Foundation of Korea (NRF) grant funded by the Korea government (MSIT) (No.2020R1A5A1016126). 

Jaewon Yoo thanks Prof. Gunhee Cho for the valuable discussions on the bounded geometry of K\"ahler manifolds.
The authors thank Prof. Vestislav Apostolov for answering a question on ALE scalar-flat K\"ahler metrics on non-compact weighted projective spaces.

\section{Hybrid Landau-Ginzburg models and Li-Wen's $L^2$-Hodge theory}\label{section2}

\subsection{Hybrid Landau-Ginzburg models}\label{section2.1}
We briefly recall the hybrid Landau-Ginzburg model studied by Chiodo and Nagel \cite{CN}. 

We let $d_1,\ldots,d_r$ and $w_1,\ldots,w_n$ be positive integers such that $(w_1,\ldots,w_n)$ are coprime and
\begin{equation}\label{CY}
    d_1+\cdots+d_r=w_1+\cdots+w_n.
\end{equation}

Let $V(\ud{W}):=V(W_1,\ldots,W_r)$ be a Calabi-Yau smooth complete intersection in the weighted projective space $\bP(\ud{w}):=\bP(w_1,\ldots,w_n)$ defined by the zero locus of $W_1, \ldots, W_r$, where $W_i$ is a weighted homogeneous polynomial of degree $d_i$ for each $i=1, \ldots, r$.
This condition \eqref{CY} also ensures that the algebraic varieties $X_{CY}$ and $X_{LG}$ below become Calabi-Yau. They are defined as the total space of locally free sheaf on the weighted projective spaces $\bP(d_1,\ldots,d_r)$ and $\bP(w_1,\ldots,w_n)$. We will write $\bP(\ud{d})=\bP(d_1,\ldots,d_r)$ and $\bP(\ud{w})=\bP(w_1,\ldots,w_n)$ for short and define
\begin{align*}
X_{CY} := \mathrm{Tot} (\cO_{\bP(\ud{w})}(-d_1)\oplus \cdots \oplus \cO_{\bP(\ud{w})}(-d_r)), \\
X_{LG} := \mathrm{Tot} (\cO_{\bP(\ud{d})}(-w_1)\oplus \cdots \oplus \cO_{\bP(\ud{d})}(-w_n)),
\end{align*}
which are both $(n+r-1)$-dimensional over $\bC$. 

Another characterization of $X_{CY}$ and $X_{LG}$ is to use the $\bC^*$-action on $\bC^n\times \bC^r$ which we call the \emph{charge} action: for $\lambda \in \bC^*$ and $(\ud{x},\ud{p})=(x_1,\ldots,x_n,p_1,\ldots,p_r)\in \bC^n\times \bC^r$,
\begin{align}
\lambda \cdot (x_1,\ldots,x_n,p_1,\ldots,p_r) =  (\lambda^{w_1} x_1,\ldots,\lambda^{w_n} x_n,\lambda^{-d_1} p_1,\ldots,\lambda^{-d_r} p_r),
\label{chargeaction}
\end{align}
$X_{CY}$ and $X_{LG}$ are isomorphic to the group quotients,
\begin{align*}
X_{CY} \cong \frac{(\bC^n\setminus \{0\})\times \bC^r}{(x_1,\ldots,x_n,p_1,\ldots,p_r) \sim (\lambda ^{w_1} x_1,\ldots,\lambda^{w_n} x_n,\lambda^{-d_1} p_1,\ldots,\lambda^{-d_r} p_r)}, \\
X_{LG} \cong \frac{\bC^n\times(\bC^r\setminus \{0\})}{(x_1,\ldots,x_n,p_1,\ldots,p_r) \sim (\lambda ^{w_1} x_1,\ldots,\lambda^{w_n} x_n,\lambda^{-d_1} p_1,\ldots,\lambda^{-d_r} p_r)}.
\end{align*}
\begin{comment}
In \cite{CN}, Chiodo and Nagel extended \cite{FJR} to consider
$$
\cH^{p,q}(W_1,\ldots,W_r) := H^{p+r,q+r}_{CR}(\cO_{\bP(\ud{d})}(-\ud{w});W^{-1}(t_0);\bC).
$$
where each $W_j$ is a quasi-homogeneous polynomial of degree $d_j$ which defines a smooth Calabi-Yau complete intersection $V(W_1,\ldots,W_r)\subset \bP(w_1,\ldots,w_n)$, and $W=p_1W_1 + \cdots + p_r W_r$. Then they proved the cohomological LG/CY correspondence \cite[Theorem 1.3]{CN} under the Calabi-Yau condition $d_1+\cdots+d_n=w_1+\cdots+w_r$, 
\begin{align}
H^{p-r,q-r}_{CR}(V(\ud{W});\bC)\cong H_{CR}^{p,q}\left(\bigoplus_{j=1}^r\cO_{\bP(w_1,\ldots,w_n)}(-d_j); W^{-1}(t_0);\bC\right),
\label{CNmain}
\end{align}
which is an orbifold version of the Thom isomorphism. In smooth case, the Chen-Ruan cohomology is nothing but the singular cohomology, so \eqref{CNmain} compares the cohomology of $V(\ud{W})$ and the relative cohomology $H^{p,q}\left(\bigoplus_{j=1}^r\cO_{\bP(w_1,\ldots,w_n)}(-d_j); W^{-1}(t_0);\bC\right)$ of the total space of the locally free sheaf $\bigoplus_{j=1}^r\cO_{\bP(w_1,\ldots,w_n)}(-d_j)$ relative to the fiber $W^{-1}(t_0)$. See \cite{CN} for details.
\end{comment}
There are projection maps $\pi_{CY}:X_{CY} \to \bP(\ud{w})$ and $\pi_{LG}:X_{LG} \to \bP(\ud{d})$.
Note that
$$
W:=\sum_{j=1}^r p_r W_r \in \C[\ud x, \ud p]
$$
defines a holomorphic function on $X_{CY}$ and $X_{LG}$. Note that $X_{CY}=\mathrm{Tot}(\cO_{\bP(\ud{w})}(-\ud{d}))$ is a non-compact Calabi-Yau manifold and
$W:X_{CY}\to \C$ has a compact critical locus $V(\ud{W})$, which is embedded as the zero section of $X_{CY}$.
\begin{definition}
    The pair $(X_{CY},W:X_{CY}\to \C)$ is called the \emph{hybrid Landau-Ginzburg model} associated with the Calabi-Yau complete intersection variety $V(\ud W)$ in $\bP(\ud{w})$.
\end{definition}

From now on, we assume that $w_1=\cdots=w_n=1$ so that $X_{CY}$ becomes a smooth complex manifold except for Section \ref{section5.2}. With this smooth condition, we can use the $L^2$-Hodge theory developed by Li-Wen \cite{LW}. Also in the later sections (except for Section \ref{section5.2}), we will denote $X_{CY}$ by $\cX$ if there is no need to distinguish $X_{CY}$ from $X_{LG}$.

When $r=1$, $X_{CY}$ is the crepant resolution of the locally Euclidean space $\bC^n / \mu_n$, where $\mu_n$ is the cyclic group of order $n$ acting diagonally on $\bC^n$. One can also observe that $X_{LG}\cong \bC^n /\mu_n$. Thus, in the $r=1$ case, $X_{CY}\to X_{LG}\simeq \bC^n / \mu_n$ is the blow-up along $0$.
However, this is not true in general. Instead, there is a simultaneous weighted blow-up $\widetilde{X}$ along the zero section $\bP(\ud{w})\hookrightarrow X_{CY}$ and $\bP(\ud{d})\hookrightarrow X_{LG}$. We provide a construction of this $\widetilde{X}$ in Section \ref{section5.2}.

\subsection{$L^2$-Hodge theory of Landau-Ginzburg models}
\label{section2.2}

We briefly review \cite{LW}. When $(X,g)$ is a Riemannian metric that preserves the almost complex structure $J$ on $TX$ (i.e. $g(X,Y)=g(JX,JY)$), the metric $g$ can be extended to the Hermitian metric on the tensor product $T^{\otimes p}_{\bC}X \otimes (T^*_\bC)^{\otimes q} X$. Using the eigenspace decomposition with respect to $J$, we have
\begin{align*}
T_{\bC}X = T^{(1,0)}X\oplus T^{(0,1)}X, && T_{\bC}^*X = T^*_{(1,0)}X\oplus T^*_{(0,1)}X.
\end{align*}
We define $T^{\otimes(p,q)}X$, $T^{(p,q)}X$, $T^*_{\otimes(p,q)}X$, and $T^*_{(p,q)}X$ as follows:
\begin{align*}
T^{\otimes(p,q)}X &= (T^{(1,0)}X)^{\otimes p}\otimes (T^{(0,1)}X)^{\otimes q}, &
T^{(p,q)}X &= (T^{(1,0)}X)^{\wedge p}\otimes (T^{(0,1)}X)^{\wedge q}, \\ 
T^*_{\otimes(p,q)}X &= (T^*_{(1,0)}X)^{\otimes p}\otimes (T^*_{(0,1)}X)^{\otimes q}, &
T^*_{(p,q)}X &= (T^{(1,0)}X)^{\wedge p}\otimes (T^{(0,1)}X)^{\wedge q}.
\end{align*}
When $\alpha, \beta \in T^{\otimes(p,q)}X \otimes T^*_{\otimes(s,t)} X$, the Hilbert-Schmidt inner product $g(\alpha,{\beta})$ is defined as follows:
\begin{align}
g(\alpha,{\beta}) := \sum_{\ud{i},\ud{j},\ud{k},\ud{l}} \alpha^{i_1\cdots i_p\overline{k_{p+1}\cdots k_{p+q}}}_{j_1\cdots j_s\overline{l_{s+1}\cdots l_{s+t}}}\ \overline{\beta^{k_1\cdots k_p\overline{i_{p+1}\cdots i_{p+q}}}_{l_1\cdots l_s\overline{j_{s+1}\cdots j_{s+t}}}}g_{i_1\overline{k_1}}\cdots g_{i_{p+q}\overline{k_{p+q}}}g^{j_1\overline{l_1}}\cdots g^{j_{s+t}\overline{l_{s+t}}}.
\label{GeneralHermitianMetric}
\end{align}
We let $|\alpha|_g := \sqrt{g(\alpha,{\alpha})}$. The operator norm $\| \cdot \|_{g}$ of tensors is defined as follows: for $\alpha\in (T_{\bC}X)^{\otimes p}\otimes (T_{\bC}^*X)^{\otimes q}$
\[
\|\a\|_{g} \ :=\ 
\sup_{v_1, \cdots , v_q\neq 0}\frac{|\a(v_1, \dots , v_q)|_{g}}{|v_1|_{g}\cdots |v_q|_{g}}.
\]
The Hilbert-Schmidt norm and the operator norm are equivalent: for each fixed $(r,s)$ there exists a constant $C=C(r,s)\ge 1$
such that pointwise on $X$,
\[
\|\a\|_{g}\ \le\ |\a|_{g}\ \le\ C\,\|\a\|_{g}.
\]
\color{black}

\begin{definition}[Bounded Geometry]
\label{Definition:bg}
Let $(X,g)$ be a complete Riemannian manifold. $(X,g)$ is said to have bounded geometry if 
\begin{enumerate}
\item 
the Riemannian curvature tensor and its derivatives are all uniformly bounded above:
$$
\sup_{p\in X} |\nabla^k R|_g < \infty,
$$
\item the injectivity radius is bounded from below by a positive number:
$$
\inf_{p\in X} \{\mathrm{inj.rad.}(p)\}>0.
$$
\end{enumerate}
\end{definition}
When condition (2) is dropped, we say that $(M,g)$ is of \emph{quasi-bounded geometry}.

These conditions appear in many literatures regarding non-compact Riemannian manifolds. They are useful since Sobolev spaces are guaranteed to be well behaved on a Riemannian manifold of bounded geometry. In particular, Sobolev's embedding theorem and Density theorem hold. See \cite{JEi} and \cite{JEl} for basic results on manifolds of bounded geometry. 

We say that a complex manifold $X$ is Calabi-Yau if its canonical line bundle is trivial
$$
K_X \cong \cO_X.
$$

\begin{definition}[{\cite[Definition 2.24]{LW}}]
Let $E$ be a Hermitian vector bundle on $X$ with the Hermitian metric $g$ and the corresponding Hermitian connection $\nabla$. The metric $g$ and its connection $\nabla$ extend uniquely onto the tensor products of $E$, $E^*$, $TX$, and $T^*X$. We define $C_{b}^\infty(E)$ as a collection of smooth sections $s$ of $E\to X$ such that the norm of $\nabla^k s$ is bounded on $X$ for each $k\geq 0$: for each $k\geq 0$, there is a positive real number $C_k$ such that
$$
|\nabla^k s|_g\leq C_k.
$$
\label{Definition:bounded}
\end{definition}

\begin{definition}[Bounded Calabi-Yau Geometry, {\cite[Definition 3.6]{LW}}]
\label{Definition:bcyg}
Let $(X,g)$ be a complete K\"ahler manifold of bounded geometry. We further assume that $X$ is Calabi-Yau and ${\Omega}\in K_X$ is a non-vanishing holomorphic $(n,0)$-form. Let $\Theta$ be the holomorphic section of $\wedge^n TX$ such that $\Theta\vdash{\Omega}=1$. A triple $(X,g,{\Omega})$ is called a \emph{bounded Calabi-Yau geometry} if ${\Omega}\in C_b^{\infty}\left(\bigwedge^n T^*_{\bC}X\right)$ and $\Theta\in C_b^{\infty}(\bigwedge^n T_{\bC}X)$.
\end{definition}

An additional condition on $f$ is needed to construct a Frobenius manifold structure, which is called \emph{strongly elliptic}.

\begin{definition}[Strongly Elliptic, {\cite[Theorem 2.6]{LW}}]
\label{Definition:strell}
Let $(X,g)$ be a complex manifold with K\"ahler metric $g$. A holomorphic function $f$ on $X$ is said to be \emph{strongly elliptic} if for all $\epsilon>0$, $k\geq 2$,
\begin{align}
\epsilon|\nabla f(z)|_g^k - |\nabla^k f(z)|_g \to +\infty
\label{Definition:strelleq}
\end{align}
as $z\to \infty$.
\end{definition}

The $k$-th covariant derivative $\nabla^k f$ is considered as a holomorphic symmetric $k$-form on $X$. The notation $z\to \infty$ means $d(z,z_0)\to \infty$, where $z_0\in X$ is arbitrary fixed point, and $d(z,z_0)$ is the distance between $z$ and $z_0$.

Let $\cA^{i,j}(X)$ be the space of smooth $(i,j)$-forms on $X$ and
\begin{align}
\cA^n(X) := \bigoplus_{i+j=n} \cA^{i,j}(X),\quad \cA(X) := \bigoplus_n \cA^n(X).
\label{spaceofsmoothforms}
\end{align}
The de Rham differential $d$ decomposes into $\partial + \overline{\partial}$,
$$
\partial: \cA^{i,j}(X) \to \cA^{i+1,j}(X),\, \overline{\partial}: \cA^{i,j}(X) \to \cA^{i,j+1}(X).
$$
Let $\cA_c^{i,j}(X)$ be a subspace of $\cA^{i,j}(X)$ consisting of compactly supported smooth $(i,j)$-forms. The twisted differential $\overline{\partial}_f$ is defined as $\overline{\partial}_f   :=\overline{\partial} + df\wedge$.

An element in $\cA^{i,j}(X)$ is a section of the bundle $T^*_{(i,j)}X$ on $X$, and by the injective map $\iota:T^*_{(i,j)}X\to T^*_{\otimes(i,j)}X$ that maps $v_1\wedge\cdots\wedge v_k \mapsto \frac{1}{k!}\sum_{\sigma\in S_k} \sgn(\sigma) v_{\sigma(1)}\otimes \cdots \otimes v_{\sigma(k)}$ that identifies the wedge product with the subspace inside the tensor product, an  element $\alpha\in \cA^{i,j}(X)$ induces a section $\iota \circ \alpha$ of the bundle $T^*_{\otimes(i,j)}X$ on $X$. We can define a pointwise inner product $(-,-)_\cA$ on $\cA^{i,j}(X)$ using \eqref{GeneralHermitianMetric}:
$$
(\alpha,\beta)_{\cA} = g(\iota\circ \alpha,\iota\circ \beta).
$$
We further define $L^2$ inner product:
$$
\langle \phi,\psi \rangle_{\cA} = \int_X (\phi,\psi)_\cA dv_g,
$$
where $\phi,\psi \in \cA^{i,j}(X)$, $v_g$ is the volume on $X$ induced by $g$. The $L^2$-norm is denoted by $\|-\|_\cA$. The $L^2$-norm is not finite for all elements in $\cA^{i,j}(X)$. We define $L^2_{\cA}(X)$ as the completion with respect to $\|-\|_\cA$ of the subspace of the forms in $\cA(X)$ that are bounded with respect to the same norm. Let $\overline{\partial}_f^*$ be defined as the adjoint of $\overline{\partial}_f$ with respect to the inner product $\langle -,-\rangle_{\cA}$ on $L^2_{\cA}(X)$.
\begin{definition}[{\cite[Definition 2.19]{LW}}]
The $f$-twisted Sobolev space $\cA_{f,k}(X)$ of smooth forms are subspaces of $L^2_{\cA}(X)$ defined as
$$
\cA_{f,k}(X):= \{\phi: (\overline{\partial}_f+\overline{\partial}_f^*)^i \phi \in L^2_\cA(X)\textrm{ for all $i=0,\ldots,k$}\}.
$$
The $\cA_{f,k}$-norm is defined as
$$
\|\phi\|_{\cA_{f,k}}:= \sum_{0\leq i\leq k}\|(\overline{\partial}_f+\overline{\partial}_f^*)^i \phi \|_{\cA}.
$$
$\cA_{f,\infty}(X)$ is defined as the intersection of all $\cA_{f,k}(X)$,
$$
\cA_{f,\infty}(X) := \bigcap_{k\geq 0} \cA_{f,k}(X).
$$
\end{definition}

We define similarly for $PV(X)$, $PV_{f,\infty}(X)$, and  $PV_c(X)$.
\begin{align*}
PV^{i,j}(X)&:= \cA^{0,j}(X,T^{(i,0)}X), \\
PV^k(X) &:= \bigoplus_{j-i=k} PV^{i,j}(X), \\
PV(X) &:= \bigoplus_{k} PV^j(X).
\end{align*}
$PV_c^{i,j}(X)$ is a subspace of $PV^{i,j}(X)$ consisting of elements with compact support. $PV_c(X)$ is defined as a direct sum,
$$
PV_c(X) := \bigoplus_{i,j} PV_c(X).
$$
Similar to the $L^2$-inner product $\langle -,- \rangle_{\cA}$, we define the pointwise Hermitian product $(-,-)_{PV}$ on $PV(X)$, $L^2$-inner product $\langle -,-\rangle_{PV}$, and the completion $L^2_{PV}(X)$ of the subspaces of $PV(X)$ with bounded $L^2$-norms.

\begin{definition}[{\cite[Definition 3.3]{LW}}]
The Sobolev spaces $PV_{f,k}(X)$ of polyvector fields are
$PV_{f,k}(X) := \{\alpha: |\nabla f|^i \nabla^j \alpha\in L^2_{PV}(X)\textrm{ for all $i+j \leq k,$}\}$
and the $PV_{f,k}$-norm is
$$
\|\alpha\|_{PV_{f,k}}:= \sum_{i+j\leq k} \| |\nabla f|^i \nabla^j \alpha \|_{PV}.
$$
$PV_{f,\infty}$ is the intersection of all $PV_{f,k}(X)$,
$$
PV_{f,\infty}(X) := \bigcap_{k\geq 0} PV_{f,k}(X).
$$

\end{definition}

We summarize the main result of \cite{LW}. 

\begin{theorem}[{\cite[Theorem 1.1, Lemma 3.7]{LW}}]
Let $(X,g,{\Omega})$ be a bounded Calabi-Yau geometry, and $f:X\to \bC$ be a strongly elliptic holomorphic function on $X$ with compact critical set.
\begin{enumerate}
\item $(PV_{f,\infty}(X),\overline{\partial}_f,\partial)$ forms a dGBV algebra with a trace pairing and the Hodge-to-de Rham degeneration holds.
\item We have quasi-isomorphic embeddings of complexes
$$
(PV_c(X),\overline{\partial}_f) \subset (PV_{f,\infty}(X),\overline{\partial}_f) \subset (PV(X),\overline{\partial}_f).
$$
In particular, the Hodge-to-de Rham degeneration holds for $(PV(X),\overline{\partial}_f,\partial)$ as well.
\item The trace pairing defines a non-degenerate pairing on the cohomology of $(PV_{f,\infty}(X),\overline{\partial}_f,\partial)$. It also induces a sesquilinear pairing on $H(PV_{f,\infty}(X)[\![u]\!],\overline{\partial}_f+u\partial)$ that generalizes K. Saito's higher residue pairing.
\item There exists a Frobenius manifold structure on the cohomology $H(PV(X),\overline{\partial}_f)$.
\end{enumerate}
\label{LWmain}
\end{theorem}

\begin{comment}
\begin{theorem}[{\cite[Theorem 2.34]{LW}}]
The natural embeddings of chain complexes
$$
(\cA_{c}(X),\overline{\partial}_W) \stackrel{i_1}{\hookrightarrow} (\cA_{f,\infty},\overline{\partial}_W) \stackrel{i_2}{\hookrightarrow} (\cA(X),\overline{\partial}_W)
$$
are quasi-isomorphisms.
\label{LWmain1}
\end{theorem}

\begin{definition}[{\cite[Theorem 2.38]{LW}}]
The pairing $\cK:\cA_{f,\infty}\times \cA_{-f,\infty}\to \bC$ on the $f$-twisted spaces defined as
$$
\cK(\alpha,\beta) = \int_X \alpha\wedge \beta
$$
induces a pairing on cohomologies,
\begin{align*}
&\cK:H(\cA_c(X),\overline{\partial}_f)\times H(\cA_c(X),\overline{\partial}_{-f}) \to \bC, \\
&\cK: H(\cA_{f,\infty}(X),\overline{\partial}_f)\times H(\cA_{f,\infty}(X),\overline{\partial}_{-f})\to \bC\\
&\cK: H(\cA(X),\overline{\partial}_f)\times H(\cA(X),\overline{\partial}_{-f})\to \bC.
\end{align*}
\end{definition}

\begin{lemma}[{\cite[Lemma 3.7]{LW}}]
$$
PV_{f,\infty}(X)\vdash \Omega_X = \cA_{f,\infty}(X).
$$
\label{LWmain2}
\end{lemma}

\begin{theorem}[{\cite[Theorem 3.26]{LW}}]
There exists a Frobenius manifold structure on the cohomology $H(PV(X),\overline{\partial}_f)$.
\label{LWmain4}
\end{theorem}
\end{comment}

\section{Bounded Calabi-Yau geometry and strong ellipticity}
\label{section3}

From now on, for simplicity, we will denote $X_{CY}$ by $\cX$ and $(\bC^n \setminus \{0\}) \times \bC^r$ by $\pcU$. Also, for a bundle $E\to B$, we denote the space of smooth sections of $E\to B$ by $\Gamma(E)$. We let $z_1,\cdots,z_{n+r}$ be the coordinate functions on $\pcU$ such that $z_j=x_j$ for $j=1,\ldots,n$, $z_j=p_{j-n}$ for $j=n+1,\cdots,n+r$. 
We define an additive grading on the variables, which we call the \textit{charge}: 
\begin{align} \label{eq:def_of_ch}
\ch(z_j)=1, \quad j=1,\ldots,n, \quad \ch(z_j)=-d_{j-n}, \quad j=n+1,\ldots,n+r. 
\end{align}
For each $j=1,\ldots, n$, we define $U_j\subset \cX$ to be the affine open subset defined by $x_j\neq0$, and we let $(Z_1,\ldots,Z_{j-1},Z_{j+1},\ldots,Z_{n+r})$ be the coordinate functions on $U_j$ defined by 
\[
Z_k = x_j^{-\ch(z_k)}z_k, \quad k=1,\ldots,j-1,j+1,\ldots,n+r.
\]
Then $\{U_1,\ldots,U_n\}$ is an open cover of $\cX$. The affine open subsets $U_1,\ldots,U_n$ carry a holomorphic structure, and the gluing maps are holomorphic, so it defines a natural holomorphic structure $J_\cX$ on $\cX$. We denote the quotient map $\pcU\to \cX$ by $\pi_{\cX}$.


{
Let $\widetilde{J}$ be the usual complex structure on $\pcU$ given by $\widetilde{J}\partial_{z}=i\partial_z$, $\widetilde{J}\partial_{\overline{z}}=-i\partial_{\overline{z}}$. We define 
\begin{align}
\label{eq:def_of_uw}
u(z):=\sum_{j=1}^{n+r}|z|^2, \quad w(z):=\sum_{j=1}^{n+r}\ch(z_j)|z|^2. 
\end{align}

{
Before proceeding with the explicit construction of the K\"ahler metric on $\cX$, we outline the geometric intuition behind the K\"ahler reduction and the introduction of the deformation parameter $s$. A naive approach might equip the ambient space $\pcU$ with the standard Euclidean K\"ahler potential corresponding to $s=0$. However, under K\"ahler reduction, this Euclidean choice induces a metric that causes the natural holomorphic volume form $\Omega$ to diverge at infinity, while its dual polyvector field $\Theta$ decays too rapidly, thereby failing the bounded Calabi-Yau geometry condition required for the subsequent $L^2$-Hodge theory. To resolve this imbalance, we introduce the parameter $s > 0$ to define a one-parameter family of K\"ahler potentials $u^{1+s}$. This parameter acts to artificially dampen the radial conical growth of the reduced metric. As demonstrated through the subsequent localized coordinate estimates, this parameter is determined by the topological dimension of the space to counterbalance the wedge product dimensions of $\Omega$ and $\Theta$, ultimately anchoring the geometry at the balanced value of $s = \frac{1}{n+r-1}$; see Proposition \ref{expcompvolumeform}.
}

We fix a real number $s\geq 0$. Note that
\begin{align}
u^{1+s}(\ud{z})=\left(\sum_{j=1}^{n+r}|z_j|^2\right)^{1+s}, \quad s \geq 0,
\label{Kahler potential}
\end{align}
is a $C^\infty$ real-valued function and is plurisubharmonic in $\pcU$, i.e. its complex Hessian matrix $H_{u^{s+1}} = \left( \frac{\partial^2 u^{s+1}}{\partial z_j \partial \bar{z}_k} \right)$ is positive definite at every point in $\pcU$.
Thus, the metric $\widetilde{g}_s$ given by  
$$\sum_{j,k} \frac{\partial^2 u^{1+s}}{\partial z_j \partial \bar{z}_k}(dz_j \otimes d\bar{z}_k + d\bar{z}_k \otimes dz_j)$$
is a well-defined Riemannian metric.
Moreover, the associated $(1,1)$-form \[i\partial\overline{\partial}u^{1+s}=\widetilde\omega_s = i \sum_{j,k} \frac{\partial^2 u^{s+1}}{\partial z_j \partial \bar{z}_k} dz_j \wedge d\bar{z}_k\] is closed. Therefore, $(\pcU, \widetilde{g}_s, \widetilde{J})$ is a K\"ahler manifold. 
}

There is a $U(1)$-action defined as
\[
e^{i\theta}\cdot (z_1,\cdots,z_{n+r}) = (e^{i\ch(z_1)\theta}z_1,\dots,e^{i\ch(z_{n+r})\theta}z_{n+r}), \qquad e^{i\theta}\in U(1),
\]
and it preserves the K\"ahler form $\widetilde{\omega}_s$. We let $\theta$ be the coordinate function on $U(1)$ given by
\[
\theta \mapsto e^{i\theta}\in U(1),
\]
and the Lie algebra $\Lie U(1)$ of $U(1)\subset \bC^*$ is isomorphic to $\bR$. For an element $\x \in \Lie U(1)$, we define the \emph{fundamental vector field} $\widehat{\xi}$ as a vector field on $\pcU$ evaluated at $z\in \pcU$:
\[
\widehat{\xi}(z) = \left.\frac{d}{dt}\exp(t\xi)\cdot z\right|_{t=0},
\]
where $\exp:\Lie U(1)\to U(1)$ is the exponential map. If we put $\xi=i\partial_{\theta}$, we get
\begin{align}\label{eq:std.vf}
\widehat{\xi} = \sum_{j=1}^{n+r}i\ch(z_j) \left(z_j\partial_{z_j}-\overline{z}_j\partial_{\overline{z}_j}\right).
\end{align}

\begin{definition}
Let $(M,\omega)$ be a symplectic manifold, and $G$ be a Lie group acting on $M$, whose action preserves $\omega$. Then $\mu:M\to (\Lie G)^*$ is called a \emph{moment map} if
\begin{align}
d(\langle \mu,\xi\rangle) = \iota_{\widehat{\x}}\omega, \label{moment_equation}
\end{align}
where $\iota_{\widehat{\x}}\omega$ is the contraction of the vector field $\widehat{\xi}$ with $2$-form $\omega$. If a moment map is given with such a condition, $G$-action on $(M,\omega)$ is called \emph{Hamiltonian}, and $(M,\omega)$ is called a \emph{Hamiltonian $G$-space}. For each $\tau \in (\Lie G)^*$, $\mu^{-1}(\tau)\subset M$ is called a \emph{moment level}.
\end{definition}

{
When the K\"ahler potential is given by $f(u)$ where $f:\bR^{>0}\to \bR^{>0}$ is a strictly convex increasing function, $f(u)$ is again plurisubharmonic since the complex Hessian 
\[
\frac{\partial^2 (f \circ u)}{\partial z_j \partial \bar{z}_k} = f'(u) \delta_{j\overline{k}} + f''(u) \overline{z}_j z_k
\]
is positive definite.
This is positive since $f$ is strictly convex and increasing, i.e. $f'(u),\ f''(u)>0$. We find the moment map condition \eqref{moment_equation} conceptually using the K\"ahler potential \(f(u)\),
\[
\tilde\omega_s \;=\; i\partial\bar\partial f(u) \;=\; dd^cf(u),
\qquad
d^c := \frac{i}{2}\,(\bar\partial-\partial).
\]
Since \(u=\sum_{j=1}^{n+r}|z_j|^2\) is \(U(1)\)-invariant, so is \(f(u)\); hence
\(\cL_{\hat\xi}f(u)=0\) where $\widehat{\xi}$ is the fundamental vector field given in \eqref{eq:std.vf}.
Define
\[
\mu \;:=\; -\,\iota_{\hat\xi}(d^cf(u)).
\]
Then, by Cartan's formula and \(\cL_{\hat\xi}(d^cf(u))=d^c(\cL_{\hat\xi}f(u))=0\),
\[
d\mu \;=\; -\,d\iota_{\hat\xi}(d^cf(u))
\;=\; \iota_{\hat\xi}dd^cf(u)
\;=\; \iota_{\hat\xi}\tilde\omega_s,
\]
which is exactly \eqref{moment_equation}.

It remains to identify \(\mu\) explicitly. Since \(d^c\rho=f'(u)\,d^c u\), we compute
\[
d^c u \;=\; \frac{i}{2}\sum_{j=1}^{n+r}\bigl(z_j\,d\bar z_j-\bar z_j\,dz_j\bigr),
\quad
dz_j(\hat\xi)= i\,\mathrm{ch}(z_j)\,z_j,\ \ d\bar z_j(\hat\xi)= -i\,\mathrm{ch}(z_j)\,\bar z_j,
\]
and hence
\[
\iota_{\hat\xi}(d^c u)
=\frac{i}{2}\sum_{j=1}^{n+r}\Bigl(z_j\,d\bar z_j(\hat\xi)-\bar z_j\,dz_j(\hat\xi)\Bigr)
=\sum_{j=1}^{n+r}\mathrm{ch}(z_j)\,|z_j|^2
=: w(z).
\]
Therefore,
\[
\mu \;=\; -\,\iota_{\hat\xi}(d^c\rho)
\;=\; -\,f'(u)\,\iota_{\hat\xi}(d^c u)
\;=\; -\,f'(u)\,w(z).
\]
Now we let $f(u)=u^{1+s}$, and define a map $\mu_s:\pcU \to (\Lie U(1))^*\cong \bR$ as
\[
\mu_s(z) = -f'(u)w(z)=-(1+s)u^s(z)w(z), \quad z\in \pcU.
\]
Then $\mu_s:\pcU \to \bR$ is a moment map for the $U(1)$-action on $(\pcU,\widetilde{\omega}_s).$
}
\begin{comment}
{\begin{proposition}
The map $\mu_s:\pcU \to i\bR\partial_{\th}^*$ is a moment map for $U(1)$-action on $(\pcU,\widetilde{\omega}_s).$
\end{proposition}
\begin{proof}
We verify that $\mu_s$ satisfies \eqref{moment_equation}; it is enough to put $\xi=i\partial_{\theta}$ because the equation \eqref{moment_equation} is linear with respect to $\xi$:
\[
d(\langle \mu,\xi\rangle) = -d\left((1+s)\rho_1\rho_0^s\right) =-(1+s)\sum_{j=1}^{n+r} \left(s\rho_1\rho_0^{s-1}+\ch(z_j)\rho_0^s\right)\left(z_jd\overline{z}_j+\overline{z}_jd{z}_j\right),
\]
and
\begin{align*}
\iota_{\widehat{\xi}}\widetilde{\omega}_s&=\sum_{1\leq j,k\leq n+r}i(\delta_{jk}(1+s)\rho_0^s+s(1+s)\rho_0^{s-1}z_k\overline{z}_j)(i\ch(z_j)z_jd\overline{z}_k+i\ch(z_k)\overline{z}_k dz_j)\\
&=-\sum_{j=1}^{n+r}((1+s)\ch(z_j)\rho_0^s)(z_jd\overline{z}_j+\overline{z}_j dz_j)-\sum_{j=}^{n+r}s(1+s)\rho_0^{s-1}\ch(z_j)|z_j|^2\sum_{k=1}^{n+r}(z_kd\overline{z}_k +\overline{z}_kdz_k)\\
&=-\sum_{j=1}^{n+r}((1+s)\ch(z_j)\rho_0^s+s(1+s)\rho_0^{s-1}\rho_1)(z_jd\overline{z}_j+\overline{z}_j dz_j).
\end{align*}
\end{proof}}
\end{comment}

\begin{proposition}[K\"ahler reduction]\label{kahler reduction} Let $G$ be a compact Lie group, $(M,g,J)$ be a K\"ahler manifold, and $\mu:M\to (\Lie G)^*$ be a moment map with respect to $G$-action on $(M,\omega)$ where $\omega(-,-)=g(-,J-)$ is the K\"ahler $2$-form. We assume further that $\tau\in (\Lie G)^*$ is a regular value of $\mu$, $\tau$ is invariant under the coadjoint action of $G$ on $(\Lie G)^*$, and $G$ acts freely, holomorphically, and isometrically on $M$. Then there is a complex structure and K\"ahler metric on the quotient space $M_\tau:=\mu^{-1}(\tau)/G$ induced from $M$.\label{KahlerReduction}
\end{proposition}
This is a standard application of symplectic reduction, originating from the Marsden-Weinstein-Meyer theorem \cite{MWM}. K\"ahler reduction (or K\"ahler quotient) in the present setting goes back to \cite{HKLR}; see also \cite{Kirwan} for a detailed exposition. Here we provide a short summary of the proof.
\begin{proof}
Since $\tau$ is a regular value of $\mu$, the moment level $\mu^{-1}(\tau)$ is a smooth submanifold of $M$.
Moreover, $\tau$ is $\Ad^*$-invariant, hence $\mu^{-1}(\tau)$ is $G$-invariant.
As $G$ acts freely and properly on $\mu^{-1}(\tau)$ (freely by assumption, and properly since $G$ is compact),
the quotient $M_\tau$ is a smooth manifold and the projection
\[
\pi:\mu^{-1}(\tau)\to M_\tau
\]
is a principal $G$-bundle.

In particular, there exists a unique symplectic form $\omega_\tau$ on $M_\tau$ characterized by
\[
\iota^*\omega=\pi^*\omega_\tau,
\]
where $\iota:\mu^{-1}(\tau)\hookrightarrow M$ is the inclusion.

Finally, since the $G$-action is holomorphic and isometric, the vertical distribution
$\mathcal V:=\ker(d\pi)=\{\,\widehat{\xi} \mid \xi\in\Lie(G)\,\}$ on $\mu^{-1}(\tau)$ is $J$-invariant.
Let $\mathcal H:=\mathcal V^{\perp_g}\subset T\mu^{-1}(\tau)$ be the $g$-orthogonal complement.
If $v\in\mathcal H$, then for any $\xi\in\Lie(G)$ we have
\[
g(Jv,\widehat{\xi})=\omega(v,\widehat{\xi})=-\omega(\widehat{\xi},v)=-d\langle\mu,\xi\rangle(v)=0,
\]
where the last equality uses $v\in T\mu^{-1}(\tau)=\ker(d\mu)$.
Hence $J(\mathcal H)=\mathcal H$, and $J$ induces an almost complex structure $J_\tau$ on $M_\tau$.
Moreover, $\omega_\tau$ is of type $(1,1)$ with respect to $J_\tau$, and thus
\[
g_\tau(-,-):=\omega_\tau(-,J_\tau-)
\]
defines a K\"ahler metric on $M_\tau$.
\end{proof}

\begin{proposition}
There is a K\"ahler structure on $\mu^{-1}_s(-1)/U(1)$.
\end{proposition}
\begin{proof}
We use Proposition \ref{KahlerReduction}. First, we show that $-1\in \bR$ is a regular value of $\mu_s$. Using the Jacobian criterion on the singular locus of $\mu^{-1}_s(-1)$, for each $j=1,\dots,n+r$, we must have
\[
\frac{\partial \mu_s}{\partial z_j}=0 \Leftrightarrow (\ch(z_j)u+sw)z_j =0,
\]
and since $(x_1,\ldots,x_n)\neq 0$ in $\pcU$, we must have
\[
u+sw=0,
\]
which implies that
\[
\mu_s(z) = \frac{1+s}{s}u^{1+s}.
\]
However, $u>0$, so $\mu_s(z)>0$. However, $\mu_s(z)=-1$. Thus, the singular locus of $\mu_s^{-1}(-1)$ is empty, so $-1$ is a regular value of $\mu_s$.

The group $U(1)$ is abelian, so its coadjoint action is trivial, so $-1$ is automatically invariant under the coadjoint action of $U(1)$. Finally, $U(1)$ obviously acts freely, holomorphically, and isometrically on $M$. Thus, by Proposition \ref{KahlerReduction}, we get the conclusion.
\end{proof}


Under the same assumption as Proposition \ref{kahler reduction}, we have the following fact regarding the complex structure $J_\tau$ on $M_\tau$.

\begin{proposition}\label{lem:descend-holomorphic}
Let $U\subset M$ be a $G$-invariant open set and $f:U\to\C$ a $G$-invariant holomorphic function.
Then $f|_{U\cap \mu^{-1}(\tau)}$ descends to a unique function $\bar f$ on $\pi(U\cap \mu^{-1}(\tau))\subset M_\tau$ such that
$\pi^*\bar f=f|_{U\cap \mu^{-1}(\tau)}$, and $\bar f$ is $J_\tau$-holomorphic.
\end{proposition}
\begin{proof}
Since $f$ is $G$-invariant, $f|_{U\cap \mu^{-1}(\tau)}$ is constant along $G$-orbits in $U\cap\mu^{-1}(\tau)$; hence it
descends to a unique function $\bar f$ on $\pi(U\cap \mu^{-1}(\tau))$ satisfying $\pi^*\bar f=f|_{U\cap \mu^{-1}(\tau)}$.

To prove holomorphicity, fix $x\in U\cap \mu^{-1}(\tau)$. Let $V_x:=\ker(d\pi_x)\subset T_xZ$ be the vertical
subspace; it is spanned by the fundamental vector fields $\hat\xi_x$ ($\xi\in\mathfrak g$).
Since $f$ is $G$-invariant, for every $\xi\in\mathfrak g$ we have
\[
df_x(\hat\xi_x)=0,
\]
hence $df_x|_{V_x}=0$. As $f$ is holomorphic on $U$, we also have
\[
df_x\circ J = i\,df_x,
\]
so for $v\in V_x$,
\[
df_x(Jv)= i\,df_x(v)=0,
\]
and therefore $df_x|_{J V_x}=0$.

Let $H_x\subset T_xZ$ be the horizontal subspace in K\"ahler reduction (i.e. the $g$-orthogonal
complement of $V_x$ in $T_xZ$). Then $T_xZ=V_x\oplus H_x$ and $J(H_x)=H_x$, and $d\pi_x|_{H_x}:H_x\to
T_{\pi(x)}M_\tau$ is an isomorphism. The reduced complex structure $J_\tau$ is characterized by
\[
d\pi_x\circ J = J_\tau\circ d\pi_x \quad\text{on }H_x.
\]

Define a complex-valued covector $\alpha$ on $T_{\pi(x)}M_\tau$ by
\[
\alpha(d\pi_x(h)) := df_x(h)\qquad (h\in H_x).
\]
This is well-defined because $d\pi_x|_{H_x}$ is an isomorphism. Moreover,
\[
\alpha(J_\tau\, d\pi_x(h)) = \alpha(d\pi_x(Jh)) = df_x(Jh) = i\,df_x(h)= i\,\alpha(d\pi_x(h)),
\]
so $\alpha$ is complex linear with respect to $J_\tau$. Since $d\bar f_{\pi(x)}=\alpha$, this shows that
$\bar f$ satisfies the Cauchy--Riemann equation with respect to $J_\tau$ at $\pi(x)$. As $x$ was arbitrary,
$\bar f$ is $J_\tau$-holomorphic on $\pi(U\cap \mu^{-1}(\tau))$.
\end{proof}

\begin{proposition}\label{biholomorphism}
We have a smooth biholomorphism $\mu^{-1}_s(-1)/U(1)\cong \cX$ which makes the diagram commute:
\[\begin{tikzcd}
	{\mu^{-1}_s(-1)} & \pcU \\
	{\mu^{-1}_s(-1)/U(1)} & {\cX=\pcU/\C^*}.
	\arrow[hook, from=1-1, to=1-2]
	\arrow[two heads, from=1-1, to=2-1]
	\arrow[two heads, from=1-2, to=2-2]
	\arrow["\cong"', from=2-1, to=2-2]
\end{tikzcd}\]
\end{proposition}

\begin{proof}
We first define a map $\phi:\mu^{-1}_s(-1)/U(1)\to \cX$. For $z\in \mu^{-1}_s(-1)$, we denote the $U(1)$-orbit of $z$ by $[z]_{U(1)}$. Similarly, for $z\in \pcU$, we denote the $\bC^*$-orbit of $z$ by $[z]_{\bC^*}$. Then we define $\phi$ as \[\phi([z]_{U(1)}):=[z]_{\bC^*}.\]

We describe the inverse of $\phi$ first. For any given $z=(z_1,\dots,z_{n+r})\in \pcU$, we consider the subset of $\bC^*$-orbit given by a parameter $t\in \bR$
\[
e^{t}\cdot z = e^t\cdot(z_1,\dots,z_{n+r})=(e^{\ch(z_1)t}z_1,\dots,e^{\ch(z_{n+r})t}z_{n+r})\in [z]_{\bC^*}.
\]
We compute the moment map $\mu_s$ at $e^t\cdot z$:
\[
\mu_s(e^{t}\cdot z) = -(1+s)w(e^{t}\cdot z)u^s(e^{t}\cdot z),
\]
and for any $t\in \bR$,
\[
\frac{dw(e^{t}\cdot z)}{dt} = \sum_{j=1}^{n+r}\ch(z_j)^2e^{2\ch(z_j)}|z_j|^2>0,
\]
so
\[
\frac{d\mu_s(e^{t}\cdot z)}{dt} = -(1+s)\frac{dw(e^{t}\cdot z)}{dt}u^s(e^{t}\cdot z) - s(1+s)w^2(e^{t}\cdot z)u^{s-1}(e^{t}\cdot z)<0,
\]
so $\mu_s(e^{t}\cdot z)$ is strictly decreasing. We consider the behavior of $\mu_s(e^{t}\cdot z)$ as either $t\to -\infty$ or $t\to +\infty$:
\begin{align*}
\lim_{t\to \infty} u(e^{t}\cdot z)=+\infty, \qquad \lim_{t\to \infty} w(e^{t}\cdot z)=+\infty,
\end{align*}
and
\[
 \lim_{t\to -\infty} w(e^{t}\cdot z)=\begin{cases}
 0 & \textrm{if \((p_1,\dots,p_r)=0\),} \\
 -\infty & \textrm{if \((p_1,\dots,p_r)\neq0\),}
 \end{cases} \qquad \lim_{t\to -\infty} u(e^{t}\cdot z)=\begin{cases}
 0 & \textrm{if \((p_1,\dots,p_r)=0\),} \\
 +\infty & \textrm{if \((p_1,\dots,p_r)\neq0\).}
 \end{cases} 
\]
Thus,
\[
\lim_{t\to \infty} \mu_s(e^{t}\cdot z)=-\infty, \qquad\lim_{t\to -\infty} \mu_s(e^{t}\cdot z)=\begin{cases}
 0 & \textrm{if \((p_1,\dots,p_r)=0\),} \\
 +\infty & \textrm{if \((p_1,\dots,p_r)\neq0\)}.
  \end{cases} 
\]
Hence, by the intermediate value theorem, for each $z\in \pcU$, there is $\lambda(z)>0$ such that
\begin{align} \label{lambdadef}
\mu_s(e^{\lambda(z)}\cdot z)=-1.
\end{align}
We define \[\psi([z]_{\bC^*}) :=[e^{\lambda(z)}\cdot z]_{U(1)}.\]
Then, we have $\phi^{-1}=\psi$ since $\phi\circ\psi([z]_{\bC^*})=[e^{\lambda(z)}\cdot z]_{\bC^*} = [z]_{\bC^*}$. 

The map $\phi$ is smooth by its definition, and $\psi$ is smooth since it is defined implicitly by the equation \eqref{lambdadef}, and the smoothness is guaranteed by the analytic implicit function theorem: choose $j\in \{1,\dots,n\}$ arbitrarily, and we consider the affine open set $U_j \subset \cX$ given by $x_j\neq 0$. Let $G:\bR\times U_j\to \bR$ as
\begin{align*}
G(t,Z_{\bullet/j}) := \mu_s(e^t\cdot(Z_{1/j},\ldots,Z_{j-1/j},1,Z_{j+1/j},\ldots,Z_{n+r/j}))+1, \\ Z_{\bullet/j}:=(Z_{1/j},\ldots,Z_{{j-1}/j},1,Z_{{j+1}/j},\ldots,Z_{{n+r}/j})\in U_j.
\end{align*}
Then $\partial G/ \partial t>0$, so there is a smooth function $\lambda_j(Z):U_j\to \bR$ such that
$$
G(\lambda_j(Z_{\bullet/j}),Z_{\bullet/j})=0.
$$
Then $\psi|_{U_j}:U_j \to \mu^{-1}_s(-1)/U(1)$ is characterized by
$$
\psi|_{U_j}([Z_{\bullet/j}]_{\bC^*}) = [e^{\lambda_j(Z_{\bullet/j})}\cdot Z_{\bullet/j}]_{U(1)}.
$$
Since $\lambda_j(Z_{\bullet/j})$ is smooth with respect to $Z_{\bullet/j}$, $\psi|_{U_j}$ is smooth. The index $j$ is chosen arbitrarily, and $\{U_j\}_{j=1,\,\dots\,,n}$ covers $\cX$, so $\psi$ is smooth.

Finally, we must check that $\phi$ and $\psi$ are holomorphic. We only need to show that $\phi$ is holomorphic, since the inverse of a holomorphic diffeomorphism is automatically holomorphic. Choose $j\in \{1,\dots,n\}$ arbitrarily, and we let the subset $\pcU_j\subset \pcU$ defined as
\[
\pcU_j = \{z=(z_1,\dots,z_{n+r})\in \pcU: x_j \neq 0\},
\]
and for each $k=1,\ldots,j-1,j+1,\ldots,n+r$, define $\widehat{Z}_{k/j}:\pcU_j\to \bC$ as
\[
\widehat{Z}_{k/j}(z) = x_j^{-\ch(z_k)}z_k, \quad z=(z_1,\ldots,z_{n+r})\in \pcU_j.
\]
Then this is a $\bC^*$-invariant holomorphic function, so it satisfies the condition in Proposition \ref{lem:descend-holomorphic}. Hence, by the conclusion of Proposition \ref{lem:descend-holomorphic}, we have a holomorphic map $\widehat{Z}_{{k/j},\red}:\pcU_j\cap \mu^{-1}(-1)\to \bC$ defined by
\[
\widehat{Z}_{k/j,\red}([z_1,\ldots,z_{n+r}]_{U(1)}) = x_j^{-\ch(z_k)}z_k,
\]
so that $\phi^*(Z_{k/j}) = \widehat{Z}_{{k/j},\red}$. Finally, we consider $\varphi_j:\pcU_j\cap \mu^{-1}(-1) \to U_j$ defined by
\[
\varphi_j([z_1,\ldots,z_{n+r}]_{U(1)}) := [\widehat{Z}_{1/j,\red},\ldots,\widehat{Z}_{j-1/j,\red},1,\widehat{Z}_{j+1/j,\red}, \ldots,\widehat{Z}_{n+r/j,\red}]_{\bC^*},
\]
and this is again holomorphic. Then this coincides with $\phi|_{\pcU_j\cap \mu^{-1}(-1)}$, so $\phi|_{\pcU_j\cap \mu^{-1}(-1)}$ is holomorphic. Since $j$ is chosen arbitrarily and $\{\pcU_j\cap \mu^{-1}(-1)\}_{1\leq j \leq n}$ covers $\mu^{-1}(-1)$, we prove that $\phi$ is holomorphic.
\end{proof}

\begin{definition}
We define a K\"ahler metric $\mathfrak{g}_s$ as the metric inherited from $\mu^{-1}_s(-1)/U(1)$ via the biholomorphism in Proposition \ref{biholomorphism}. We let $\nabla:=\nabla(s):\Gamma(T_{\bC}\cX)\to \Gamma(T_{\bC}\cX)\otimes \Gamma(T^*_{\bC}\cX)$ be the Levi-Civita connection of $\mathfrak{g}_s$.
\end{definition}

Using biholomorphism $\mu^{-1}_s(-1)/U(1)\cong \cX$, a $U(1)$-invariant function $f:\mu_s^{-1}(-1)\to\bR$ induces a function $\cX\to \bR$ with commutative diagram
\[\begin{tikzcd}
	{\mu^{-1}_s(-1)} \\
	\cX & \bR,
	\arrow[two heads, from=1-1, to=2-1]
	\arrow["f", from=1-1, to=2-2]
	\arrow[from=2-1, to=2-2]
\end{tikzcd}\]
which we will abuse the notation $f$. For example, $u|_{\mu_s^{-1}(-1)}:\mu_s^{-1}(-1)\to\bR_{>0}$ are $U(1)$-invariant, so we can define $u|_{\mu_s^{-1}(-1)}:\cX\to \bR_{>0}$.

{
For later purposes, we provide the detailed formula for the computation of $\mathfrak{g}_s$ on the local affine coordinate chart $U_j\subset \cX$. To compute $\mathfrak{g}_s (\partial_{Z_{k/j}},\partial_{\overline{Z}_{l/j}})$, we proceed in four steps:
(1) we take the arbitrary lift of $\partial_{Z_{k/j}}$ and $\partial_{\overline{Z}_{l/j}}$ on $\pcU$, 
(2) take the projection onto the horizontal distribution $\mathcal H:= (\ker d\pi_{\cX})^{\perp_{g_s}}$, 
(3) compute the metric $\widetilde{g}_s$ of two vector fields, 
(4) and take the restriction on the appropriate moment level $\mu^{-1}(-1)$.

\begin{enumerate} 
\item The simplest lift of $\partial_{Z_{k/j}}$ we can think of is $x_j^{\ch(z_k)}\partial_{z_k}$, which is defined on $\pcU$ and vanishes only if $x_j=0$. 
\item The vertical distribution is generated by two vector fields
\[
\mathcal{V}=\left\langle\sum_{j=1}^{n+r}\ch(z_j) z_j\partial_{z_j},  \sum_{j=1}^{n+r}\ch(z_j)\overline{z}_j\partial_{\overline{z}_j}\right\rangle.
\]
We define a vector field $\mathbf{n}$ on $\pcU$ 
\begin{align}\label{eq:def_of_n}
\mathbf{n}:=\sum_{j=1}^{n+r}\ch(z_j) z_j\partial_{z_j},
\end{align} and the horizontal projection $\cH:T_{\bC}\pcU\to \mathcal{H}$ yields
\begin{align}\label{eq:def_of_H}
\cH \partial_{z_k} = \partial_{z_k} - \frac{\widetilde{g}_s(\partial_{z_k},\overline{\mathbf{n}})}{\widetilde{g}_s(\mathbf{n},\overline{\mathbf{n}})}\mathbf{n}.
\end{align}
\item Thus, if we compute $\widetilde{g}_s$ for these two projected vector fields, we get
\begin{align}\label{eq:def_of_G}
G_{k\overline{l};j}(z):=\widetilde{g}_s(x_j^{\ch(z_k)}\cH \partial_{z_k},\overline{x}_j^{\ch(z_l)}\cH\partial_{\overline{z}_l}),
\end{align}
for each $k,l=1\ldots,n+r$, which is a function on the subset $\pcU_j:=\{x_j\neq 0\}$ of $\pcU$.
\item Finally, we restrict the value of this function on the moment level $\mu_s^{-1}(1)$ to obtain the result $\mathfrak{g}_s(\partial_{Z_{k/j}},\partial_{\overline{Z}_{l/j}})$. More precisely, for $[Z_{\bullet/j}]_{\bC^*}\in U_j$, we have $\lambda_j(Z_{\bullet/j})\in \bR$ such that
\[
\mu_s(e^{\lambda_j(Z_{\bullet/j})}\cdot Z_{\bullet/j})=-1.
\]
Thus, we evaluate $G_{k\overline{l};j}$ at $e^{\lambda_j(Z_{\bullet/j})}\cdot Z_{\bullet/j}$ to obtain $\mathfrak{g}_s(\partial_{Z_{k/j}},\partial_{\overline{Z}_{l/j}})|_{[Z_{\bullet/j}]_{\bC^*}}$:
\begin{align}\label{eq:compute_g1}
G_{k\overline{l};j}(e^{\lambda_j(Z_{\bullet/j})}\cdot Z_{\bullet/j}) = \mathfrak{g}_s(\partial_{Z_{k/j}},\partial_{\overline{Z}_{l/j}})|_{[Z_{\bullet/j}]_{\bC^*}}.
\end{align}
\end{enumerate}
}


The proof of Theorem \ref{Main1} is divided into three parts (Sections \ref{section3.1}, \ref{section3.2}, and \ref{section3.3}). In Section \ref{section3.1}, we prove the following theorem.

\begin{theorem}\label{mainthm3.7}
The triple $(\cX,\mathfrak{g}_s,J_{\cX})$ is a K\"ahler manifold of bounded geometry.
\label{bg}
\end{theorem}
The proof of Theorem~\ref{mainthm3.7} is completed by proving in Section~\ref{section3.1} that $(\cX,\mathfrak{g}_s)$ is
asymptotically conical with respect to the K\"ahler cone
\[
\pcC:=\mu_s^{-1}(0)/U(1),
\]
and then invoking Proposition~\ref{ACisbg}, which shows that every asymptotically conical
manifold has bounded geometry.

More precisely, the asymptotically conical statement is established in three steps.
First, Propositions~\ref{cone1} and~\ref{cone2} identify the ambient space $\pcU$ and the reduction
$\pcC$ as natural K\"ahler cones, thereby providing the model geometry at infinity.
Second, Proposition~\ref{prop:diffeo} constructs a diffeomorphism
\[
\Phi: \pcC \longrightarrow \cX \setminus \{p_1=\cdots=p_r=0\}
\]
by moving points from the moment level $\mu_s^{-1}(0)$ to $\mu_s^{-1}(-1)$ along
the $\mathbb C^*$-orbit. Finally, Theorem~\ref{keylemma} proves the decay estimate
\[
|(\nabla^C)^k(\Phi^*{\mathfrak{g}}_s-\mathfrak{g}_{s,\pcC})|=O(\rho_{s,\pcC}^{-2-k}),
\qquad k\ge 0.
\]

The main difficulty lies in this last step. In local coordinates, negative powers of
the affine coordinates appear naturally, so one must first shrink the standard affine
cover to charts $V_j$ on which $|x_j| \asymp \sqrt{u}$; this is the content of
Lemma~\ref{prop:shrinking} and is the key geometric input. Once this is in place, the Christoffel
symbols and the metric error terms, together with all of their iterated derivatives,
admit rational expressions with controlled denominators, and the required decay
estimate follows from a systematic chartwise bookkeeping argument.

\begin{remark}
Theorem \ref{bg} holds without the Calabi-Yau condition $d_1+\cdots+d_r=n$.
\end{remark}

In Section \ref{section3.2}, we prove that there exists a bounded Calabi-Yau form in the sense of Definitions \ref{Definition:bounded} and \ref{Definition:bcyg}. 

\begin{definition} 
\label{volumeform}
We define a holomorphic volume form ${\Omega}$ on $\cX$ whose local description on the affine open set $U_j:=\{x_j\neq 0\}$;
$$U_j\cong\bA_{\bC[Z_{1/j},\ldots,Z_{j-1/j},Z_{j+1/j},\ldots,Z_{n+r/j}]}$$ 
where $Z_k = \frac{z_k}{x_j^{\ch(z_k)}}$, is
$$
{\Omega}|_{U_j} = (-1)^j dZ_{1/j} \wedge \cdots \wedge dZ_{j-1/j}\wedge dZ_{j+1/j}\wedge\cdots \wedge dZ_{n+r/j}.
$$
Similarly, we define a holomorphic polyvector field $\Theta\in \Gamma(T^{(n+r-1,0)}\cX)$ of type $(n+r-1,0)$ by
$$
\Theta|_{U_j} = (-1)^j \partial_{Z_{1/j}} \wedge \cdots \wedge \partial_{Z_{j-1/j}}\wedge \partial_{Z_{j+1/j}}\wedge \cdots \wedge \partial_{Z_{n+r/j}}. 
$$
\end{definition}
The assumption $d_1+\cdots + d_r =n$ ensures that $\Theta$ and ${\Omega}$ are well defined on $\cX$. The boundedness of the norm of $\Omega$ and $\Theta$ is controlled by the constant $s$ which is introduced to define the K\"ahler potential of $\pcU$, and in particular, when $s=\frac{1}{n+r-1}$, both $\Omega$ and $\Theta$ are bounded.

\begin{theorem}
The 4-tuple $(\cX,\mathfrak{g}_{s},J,{\Omega})$ is a bounded Calabi-Yau geometry if and only if $s=\frac{1}{n+r-1}$. If $d_{\max}\leq 2d_{\min}$, then $W:\cX\to \bC$ is strongly elliptic.
\label{bcyg}
\end{theorem}

By Theorem \ref{bg}, we know that $(\cX,\mathfrak{g}_s,J)$ is a complete K\"ahler manifold of bounded geometry. According to Definition \ref{Definition:bcyg}, it remains to show that $|\nabla^k {\Omega}|_{\mathfrak{g}_{s}}$ and $|\nabla^k \Theta|_{\mathfrak{g}_{s}}$ are bounded for all $k\geq 0$. In Proposition \ref{4.16}, we prove that $|\nabla^k \Omega|_{\mathfrak{g}_{s}}/|\Omega|_{\mathfrak{g}_{s}}$ and $|\nabla^k \Theta|_{\mathfrak{g}_{s}}/|\Theta|_{\mathfrak{g}_{s}}$ are bounded. In Proposition 
\ref{expcompvolumeform}, we prove that $\Omega$ is bounded if $\frac{1}{n+r-1}\leq s$, and $\Theta$ is bounded if $0\leq s\leq \frac{1}{n+r-1}$. In the proof, we compute $|\nabla^k\Omega|_{\mathfrak{g}_{s}}$ and $|\nabla^k\Theta|_{\mathfrak{g}_{s}}$ on the open subset $V_f$ defined in Proposition \ref{prop:shrinking}, and use the asymptotic computations established in Section~\ref{section3.1} including asymptotics of the Christoffel symbols in Proposition \ref{prop:formula_of_Gamma} and that of $|\nabla^k dZ_j|$ in Lemma \ref{cor:cov_asymp}. 

In Section \ref{section3.3}, we prove that $W=p_1W_1+\cdots+p_r W_r$ is strongly elliptic in $(X,\mathfrak{g}_s$) if $d_{\max} \leq 2d_{\min}$ and $s<1$. Note that the condition is satisfied when $s=\frac{1}{n+r-1}$ so that $(\cX,\mathfrak{g}_s,\Omega)$ is a bounded Calabi-Yau geometry. In the proof, the fact that $V(W_1,\ldots,W_r)$ is a smooth complete intersection is crucially used; it implies that either $|W_1|^2+\cdots + |W_r|^2 >0
$ is true or the Jacobian matrix
\[
\mathbf{J}:=\begin{pmatrix}
\frac{\partial W_1}{\partial x_1} & \cdots & \frac{\partial W_r}{\partial x_1} \\
\vdots & \ddots & \vdots \\
\frac{\partial W_1}{\partial x_n} &\cdots & \frac{\partial W_r}{\partial x_n} 
\end{pmatrix}
\]
is injective. Using this fact, we find the lower bound of $|\nabla W|^2$ with respect to $u|_{\mu^{-1}_s(-1)}$:
\[
|\nabla W|_{\mathfrak{g}_s}^k \geq C^2{(u|_{\mu^{-1}_s(-1)})}^{d_{\min}-s}
\]
for some constant $C$. Finally, we prove that
\[
|\nabla^k W|_{\mathfrak{g}_{s}}= O\left((u|_{\mu^{-1}_s(-1)})^{\frac{1}{2}({d_{\max}+1}-k{(1+s)})}\right)
\]
using the asymptotic computations established in Section~\ref{section3.1}.

\subsection{Asymptotically conical manifolds have bounded geometry}
\label{section3.0}
In this subsection, we first recall the definition of an asymptotically conical manifold and then prove that asymptotically conical manifolds have bounded geometry. Although this implication appears to be well known to experts, we were unable to locate a convenient reference. For completeness, we therefore include a proof that relies only on a result \cite{Ehr74} of Ehrlich.

\begin{definition}
A Riemannian manifold $(C,g)$ is called a \emph{Riemannian cone} if there exists a complete vector field $E\in \Gamma(TC)$ that is nowhere vanishing with $\nabla^g E = \id_{TM}$. The vector field $E$ is called an \emph{Euler vector field}. A nonvanishing smooth function 
\[\rho_C:C\to \bR^{>0}\]
defined by the length of the vector field $E$
\[
\rho_C(x) := \sqrt{g(E,E)|_{x}}
\]
is called the \emph{radial coordinate} of $C$. The submanifold $L_C$ of $C$ defined by
\[
L_C := \r_C^{-1}(1), \qquad g_{L_C} := g\big|_{L_C}.
\]
is called the \emph{link} of $C$. Note that $g$ can be written as
\be
d\rho_C \otimes d\rho_C + \rho_C^2 g_{L_C}.
\ee
\end{definition}

\begin{definition}
\label{ACmetric}
A geodesically complete Riemannian manifold $(M,g)$ is called an \emph{asymptotically conical Riemannian manifold}, or an \emph{AC manifold} for short, if there exist a Riemannian cone $(C,g_C)$ with a compact, connected link $L_C$, a diffeomorphism 
\be
\Phi:C\setminus K'\to M\setminus K
\ee
for a compact set $K \subset M$ and a compact set $K' \subset C$, and a real number $\alpha >0$
such that for all $k=0,1,\ldots$,
$$
|(\nabla^{g_C})^{k}(\Phi^* g - g_{C})|_{g_{C}} = O(\rho_C^{-\a-k}).
$$
The real number $\alpha$ is called the \emph{decay rate}. We say that an \emph{AC manifold} $(M,g)$ is associated with a Riemannian cone $(C,g_C)$.
\end{definition}

Recall from Definition \ref{Definition:bg} that a Riemannian manifold $(M,g)$ is said to have bounded geometry if the covariant derivatives $\nabla^k R$ of the Riemannian curvature tensor $R$ are bounded, and the injectivity radius is positive. 
{ We use Ehrlich's lower semicontinuity theorem \cite{Ehr74} for the injectivity radius in the proof of the following proposition. Because Ehrlich's continuity of the injectivity radius requires a compact manifold without boundary, we perform a doubling trick along the boundary of the cylinder to extract a uniform lower bound.
}

\begin{proposition}[asymptotically conical \(\Rightarrow\) bounded geometry]
An asymptotically conical manifold $(M,g)$ has bounded geometry.
\label{ACisbg}
\end{proposition}

\begin{proof}
Every Riemannian cone is isometric to $(\bR_{>0}\times L,d\rho^2+\rho^2g_{L})$ where $(L,g_L)$ is a (compact and connected) link of the Riemannian cone, and $\rho$ is the standard coordinate of $\bR_{>0}$. Since $(M,g)$ is an asymptotically manifold, we may assume that
$(C,g_C):=(\bR_{>0}\times L,d\rho^2+\rho^2g_{L})$ is the cone associated with $(M,g)$. From now on, we write $\nabla^C:= \nabla^{g_C}$ for simplicity. The proof is split into three parts: (1) we give an asymptotic estimate of covariant derivatives of the Riemannian curvature tensor and the injectivity radius on the cone $(C,g_C)$, (2) we compare the asymptotic estimate of covariant derivatives of the Riemannian curvature tensors of $(C,g_C)$ and that of $(M,g)$, (3) and finally we give the pointwise lower bound on the injectivity radius of $(M,g)$.

(1) This part of the proof is a standard exercise in Riemannian geometry; hence, we only state the conclusion. We use the following notation for the Riemannian curvature tensor of $g_C$:
\[
R^{C}(X,Y) := \nabla_X^{C}\nabla_Y^{C}-\nabla_Y^{C}\nabla_X^{C}-\nabla_{[X,Y]}^{g_C}, \quad
R^{C}(X,Y,Z,W) := g_C(R^{C}(X,Y)Z,W).
\] 
Then we have
\[
|(\nabla^{C})^k R^{C}|_{g_C} = O(\rho^{-2-k}).
\]
The injectivity radius on $C$ satisfies the following pointwise equality: for any $\rho>0$ and $x\in L$,
\[
\inj((\rho,x),g_C) =\rho\,\inj((1,x),g_C). 
\]

(2) We let $\nabla$ be the Levi-Civita connection of $\Phi^*g$ on $C$. We let $R$ be the Riemannian curvature tensor of $\nabla$. We define $h:=\Phi^*g - g_C$ and \[A:=\nabla -\nabla^C, \quad A(X,Y) = \nabla_X Y - \nabla^C _X Y.\] After increasing $a$ if necessary, we may assume $\|h\|_{g_{C}}\le \frac12$ on $(a,\infty)\times L$. We have the following two equalities.
For all vector fields $X,Y,Z$ on $C$, we have
\begin{equation*}
2\Phi^*g(A(X,Y),Z)
=(\nabla^{C}h)(X,Y,Z)+(\nabla^{C}h)(Y,X,Z)-(\nabla^{C}h)(Z,X,Y).
\end{equation*}
Using this, we can show that
\begin{align}\label{eq:A-decay}
|(\nabla^C)^k A|_{g_C}=O(\rho^{-1-\alpha-k}).
\end{align}
Also,
For all vector fields $X,Y,Z$ on $C$, we have
\begin{equation}\label{eq:curv-diff}
R(X,Y)Z = R^{C}(X,Y)Z
+(\nabla^{C}A)(X,Y,Z)-(\nabla^{C}A)(Y,X,Z)
+A(X,A(Y,Z))-A(Y,A(X,Z)),
\end{equation}
or schematically,
\begin{equation*}
R= R^C + \nabla^C A + A*A.
\end{equation*}
Taking $\nabla^{C}$-derivatives of \eqref{eq:curv-diff} and using Leibniz's rule yields, for each $k\ge 0$,
\begin{equation}\label{eq:Rm-diff-est}
\big|(\nabla^{C})^{k}(R-R^C)\big|_{g_{0}}
\ \le\ D_{k}\Big(\big|(\nabla^{C})^{k+1}A\big|_{g_{C}}
+\sum_{a+b=k}\big|(\nabla^{C})^{a}A\big|_{g_{C}}\ \big|(\nabla^{C})^{b}A\big|_{g_{C}}\Big),
\end{equation}
for some constant $D_k$.
Combining \eqref{eq:Rm-diff-est} with \eqref{eq:A-decay} gives
\[
\big|(\nabla^{C})^{k}(R-R^C)\big|_{g_{C}}=O(\r^{-2-\a-k}).
\]
By the first step, we have $|(\nabla^{C})^{k}R^C|_{g_{C}}=O(\r^{-2-k})$, and therefore
\begin{equation}\label{eq:Rmbound-nabla0}
\big|(\nabla^{C})^{k}R\big|_{g_{C}}=O(\r^{-2-k}).
\end{equation}
Since $\nabla=\nabla^{C}+A$ and $A$ decay as in \eqref{eq:A-decay}, each $\nabla^{k}R$
is a finite sum of contractions of $(\nabla^{C})^{j}R^C$ and $(\nabla^{C})^{j'}A$.
Using \eqref{eq:A-decay} and \eqref{eq:Rmbound-nabla0},
we conclude that for every $k\ge 0$,
\[
|\nabla^{k}R|_{\Phi^*g}=O(\r^{-2-k}).
\]

(3) We first give the rough sketch of the proof while skipping some technical parts. This part of the proof is based on the lower semicontinuity of the injectivity radius proved in {\cite[Chapter~6, Theorem~1]{Ehr74}}:
for a closed manifold $N$, the map
\[
\operatorname{inj}\colon C^{2}(\mathrm{Met}(N)) \times N \to \mathbb{R}
\]
is lower semicontinuous, where $\mathrm{Met}(N)$ is a space of Riemannian metrics on $X$, and $C^2(\mathrm{Met}(N))$ is a topological space on $\mathrm{Met}(N)$ with $C^2$-topology. For each $t>1$, we consider the submanifold $M_t:=\Phi([t/2,2t]\times L) \subset M$ and the metric $g|_{M_t}$. Using the dilation map $\delta_t:C\to C$, $\delta_t(\rho,x):= (t\rho, x)$, we can define the metric $g_t:= t^{-2}\delta_t^*\Phi^* g|_{M_t}$ on \[P:=[1/2,2]\times L.\] 
Then as $t\to \infty$, \[g_t\to g_C|_{P}\] in $C^\infty$-topology. We may apply {\cite[Chapter~6, Theorem~1]{Ehr74}} to show that 
\[
\inj((\rho,x),{g}_C|_{P}) \leq \liminf_{t\to \infty}\inj((\rho,x),{g}_t)
\]
for $(\rho,x)\in P$. If we let $m=\inj(P,{g}_C|_{P})$, we get
\[
\inj((\rho,x),{g}_t) \geq \frac{m}{2}
\]
when $t>t_0$ for $t_0$ large enough. If we put $g_t:= t^{-2}\delta_t^*\Phi^* g|_{M_t}$, we get 
\[
\inj(\delta_t(\rho,x),\Phi^* g)=t\;\inj((\rho,x),t^{-2}\delta_t^*\Phi^* g)\geq \frac{m}{2}t
\]
when $t>t_0$. Hence, we have
\[
\inj(M,g) = \min \left(\inj(\Phi([0,t_0]\times L),g|_{\Phi([0,t_0]\times L)}) , \frac{m}{2}t_0\right).
\]
The value of $\inj(\Phi([0,t_0]\times L),g|_{\Phi([0,t_0]\times L)}$ is positive since $[0,t_0]\times L$ is compact.

The technical part we skipped here is that {\cite[Chapter~6, Theorem~1]{Ehr74}} is not directly applicable to $([1/2,2]\times L,g_t)$ since $[1/2,2]\times L$ is not a closed manifold. In order to apply {\cite[Chapter~6, Theorem~1]{Ehr74}}, we construct a closed manifold $\widehat{P}:=P\cup_{\partial P} \overline{P}$ using the ``doubling" trick; here $\overline{P}$ is a copy of $P$ with flipped orientation. For more detail of this part, see Appendix \ref{appendix2}.

\end{proof}

The following fact which compares the asymptotic order of $|(\nabla^C)^k\Phi^*T|_{g_C}$ on $(C,g_C)$ and that of $|\nabla^k T|_g$ on the asymptotically conical manifold $(M,g)$, where $T$ is a tensor field on $(M,g)$, will be useful in Sections \ref{section3.2} and \ref{section3.3}.
\begin{proposition}\label{prop:tensor_compariton}
Let $(M,g)$ be the asymptotically conical manifold with the diffeomorphism $\Phi:C\setminus K \to M\setminus K'$, where $(C,g)$ is a Riemannian cone, and $K\subset C$, $K'\subset M'$ are compact subsets. We let $T$ be a tensor field on $M$. Then, for some $\n\in \bR$,
\begin{align}\label{assumption1}
|(\nabla^C)^k \;\Phi^*T|_{g_C} = O(\rho_C^{-\n-k})
\end{align}
for each $k=0,1,\ldots$, if and only if
\begin{align}\label{assumption2}
|(\nabla)^k \;\Phi^*T|_{\Phi^*g} = O(\rho_C^{-\n-k})
\end{align}
for each $k=0,1,\ldots$; 
the connection $\nabla:=\Phi^*\nabla^{g}$ is a pullback connection of the Levi-Civita connection $\nabla^{g}$ of the Riemannian metric $g$ on $M$ onto the Riemannian cone $C$.
\end{proposition}
\begin{proof}
First,
\[
|(\nabla^C)^k \;\Phi^*T|_{g_C} ^2 - |(\nabla^C)^k \;\Phi^*T|_{\Phi^*g} ^2 = h\left((\nabla^C)^k\Phi^*T,\overline{(\nabla^C)^k\Phi^*T}\right)
\]
and by the assumption that $(M,g)$ is asymptotically conical,
\[
\left|h\left((\nabla^C)^k\Phi^*T,\overline{(\nabla^C)^k\Phi^*T}\right)\right| \leq \|h\|_{g_C}|(\nabla^C)^k \;\Phi^*T|_{g_C} ^2 =O(\rho_C^{-2\n-2k-\a})
\]
where $\a$ is the decay rate of $(M,g)$. Hence,
\[
|(\nabla^C)^k \;\Phi^*T|_{g_C}  = O(\rho_C^{-\n-k})
\textrm{\quad if and only if\quad}
|(\nabla^C)^k \;\Phi^*T|_{\Phi^*g}= O(\rho_C^{-\n-k}).
\]
Recall that in the proof of Proposition \ref{ACisbg}, we defined $A:=\nabla-\nabla^C$ and obtained that
\[
|(\nabla^C)^k A|_{g_C}=O(\rho^{-1-\alpha-k})
\]
in \eqref{eq:A-decay}. Now schematically, we have
\begin{align} \label{schematic_rel}
\nabla^k \;\Phi^*T = \sum_{\substack{a+b+c=k\\a,b,c\geq 0}} (\nabla^C)^a \;\Phi^*T * (\nabla^C)^b(\underbrace{A*A*\cdots * A}_{c \textrm{ times}}).
\end{align}
Thus, if we assume that \eqref{assumption1} is true for each $k=0,1,\ldots,$ from \eqref{schematic_rel}
\begin{align*}
|\nabla^k \;\Phi^*T|_{g_C} &\leq \sum_{\substack{a+b+c=k\\ b_1+\cdots+b_c=b\\a,b,b_1,\ldots,b_c,c\geq 0}} |(\nabla^C)^a \;\Phi^*T |_{g_C} |(\nabla^C)^{b_1}A|_{g_C} \cdots |(\nabla^C)^{b_c}A|_{g_C}\\&=\sum_{\substack{a+b+c=k\\a,b,c\geq 0}}  O(\rho_C^{-\n-k-c\a}) =  O(\rho_C^{-\n-k}).
\end{align*}
To prove the converse, we use induction: we assume that \eqref{assumption2} is true for each $k=0,1,\ldots,$ and \eqref{assumption1} is true for $k=0,1,\ldots,l-1$ for some $l>0$. It is enough to prove that \eqref{assumption1} is true for $k=l$. Then from \eqref{schematic_rel},
\begin{align} 
\nabla^l \;\Phi^*T - (\nabla^C)^l \;\Phi^*T = \sum_{\substack{a+b+c=l\\a>0,c\geq 1,b\geq 0}} (\nabla^C)^a \;\Phi^*T * (\nabla^C)^b(\underbrace{A*A*\cdots * A}_{c \textrm{ times}}).
\end{align}
Thus,
\[
|\nabla^l \;\Phi^*T - (\nabla^C)^l \;\Phi^*T|_{g_C} = O(\rho_C^{-\n-k-\a}).
\]
Since $\a>0$, by the triangular inequality, we prove that  \eqref{assumption1} is true when $k=l$.
\end{proof}

When a compatible complex structure is attached to a Riemannian cone, we say that $(M,g,J)$ is a \emph{K\"ahler cone}, and we have an easier criterion for a K\"ahler cone that does not involve computing the Levi-Civita connection.

\begin{definition}[K\"ahler cone]\label{def:kahler-cone} A K\"ahler manifold $(C,g,J)$ is called a \emph{K\"ahler cone} if there exists an Euler vector field $E$ a complete vector field on $C$ such that
\begin{align} \label{eq:KahlerCone}
\cL_E g = 2g, \qquad \cL_E J = 0, \qquad dg(E,-)=0.
\end{align}
\end{definition}

\begin{proposition} \label{prop:KahlerCone1}
When $(C,g,J)$ is a K\"ahler cone with an Euler vector field $E$, then $(C,g)$ is a Riemannian cone with the Euler vector field $E$.
\end{proposition}
\begin{proof}
We write $\theta:= g(E,-)$. Then by the definition of the exterior derivative,
\begin{align}\label{eq:SymmetryOfnablaE}
d\theta(X,Y) = Xg(E,Y)-Yg(E,X)-g(E,[X,Y])= g(\nabla_X E,Y) -g(\nabla_Y E, X)=0.
\end{align}
On the other hand,
\begin{align*}
\cL_E g(X,Y) &= Eg(X,Y)-g([E,X],Y)-g(X,[E,Y])\\
&=g(\nabla_E X,Y)+g(X,\nabla_E Y)-g([E,X],Y)-g(X,[E,Y])
\\&=g(\nabla_X E,Y)+g(\nabla_X E,Y).
\end{align*}
By \eqref{eq:KahlerCone} and \eqref{eq:SymmetryOfnablaE},
\[
2g(\nabla_X E,Y)=2g(X,Y).
\]
This is true for all vectors $Y$, so we must have
$
\nabla_X E = X.
$
\end{proof}

\subsection{The K\"ahler manifold $X_{CY}$ is asymptotically conical}

\label{section3.1}

We shall prove that $(\cX,\mathfrak{g}_s)$ is an AC manifold. We have two K\"ahler cones that naturally arise in our setting.

{

\begin{proposition}\label{cone1}
The K\"ahler manifold $(\pcU,\widetilde{g}_s,\widetilde{J})$ is a K\"ahler cone with an Euler vector field $E_s$ being
\[
E_{s} = \frac{1}{1+s}\sum_{j=1}^{n+r}(z_j \partial_{z_j} + \overline{z}_j\partial_{\overline{z}_j}).
\]
\end{proposition}
\begin{proof}
We need to prove four things: we prove that $E_s$ is a complete vector field, and we prove the three equalities given in Definition \ref{def:kahler-cone}.

Completeness of $E_s$:
The flow of $E_s$ is $\varphi_t(z)=e^{t/(1+s)}z$. Since $U=\bC^n\setminus\{0\}\times\bC^r$ is invariant
under scalar multiplication, $\varphi_t$ is defined for all $t\in\bR$, therefore $E_s$ is complete.

$\mathcal L_{E_s}\widetilde J=0$:
$\varphi_t$ is holomorphic for each $s$, so $d\varphi_t \circ \widetilde J=\widetilde J\circ d\varphi_t$.
Now view $\widetilde J$ as a $(1,1)$-tensor. We get
\[
(\varphi_t^*\widetilde J)_p
=(d\varphi_t|_p)^{-1}\circ \widetilde J_{\varphi_t(p)}\circ d\varphi_t|_p
=\widetilde J_p,
\]
hence, $\varphi_t^*\widetilde J=\widetilde J$ for all $t$.
By definition of Lie derivative via flow $\varphi_t$,
\[
\mathcal L_{E_s}\widetilde J
=\left.\frac{d}{dt}\right|_{t=0}\varphi_t^*\widetilde J
=\left.\frac{d}{dt}\right|_{t=0}\widetilde J
=0,
\]
so
$\mathcal L_{E_s}\widetilde J=0$.

$\mathcal L_{E_s}\widetilde g_s=2\widetilde g_s$:
Since $E_s(u)=\frac{2}{1+s}u$, we have
\[
E_s(u^{1+s})=(1+s)u^s E_s(u)=2u^{1+s}.
\]
Because $\varphi_t$ is holomorphic, $\mathcal L_{E_s}$ commutes with $\partial,\bar\partial$, hence
\[
\mathcal L_{E_s}\widetilde\omega_s
=i\partial\bar\partial\bigl(E_s(u^{1+s})\bigr)
=i\partial\bar\partial(2u^{1+s})
=2\widetilde\omega_s.
\]
Together with $\mathcal L_{E_s}\widetilde J=0$, this gives $\mathcal L_{E_s}\widetilde g_s=2\widetilde g_s$.

$d\,\widetilde g_s(E_s,-)=0$:
Write $\theta:=\widetilde g_s(E_s,-)$. In complex coordinates,
\begin{align} \label{metricmatrix}
\widetilde g_{j\bar k}
=\frac{\partial^2 u^{1+s}}{\partial z_j\,\partial \bar z_k}
=(1+s)u^s\delta_{jk}+s(1+s)u^{s-1}\bar z_j z_k.
\end{align}
For any $m$,
\[
\theta(\partial_{z_m})
=\widetilde g_s(E_s,\partial_{z_m})
=\frac{1}{1+s}\sum_{j=1}^{n+r}\bar z_j\,\widetilde g_{m\bar j}
=(1+s)u^s\bar z_m,
\]
and similarly $\theta(\partial_{\bar z_m})=(1+s)u^s z_m$. Therefore
\[
\theta=(1+s)u^s\sum_{m=1}^{n+r}\bigl(\bar z_m\,dz_m+z_m\,d\bar z_m\bigr)
=(1+s)u^s\,du
=d(u^{1+s}),
\]
so $d\theta=0$.

Thus, $(U,\widetilde g_s,\widetilde J)$ satisfies \eqref{eq:KahlerCone}, hence is a K\"ahler cone with Euler field $E_s$.
\end{proof}

Similarly to $\cX$, $\mu^{-1}_s(0)/U(1)$ is given a K\"ahler structure via K\"ahler reduction in Proposition \ref{KahlerReduction}. We let $\pazocal{C}:= \mu^{-1}_s(0)/U(1)$, and let $\mathfrak{g}_{s,\pazocal{C}}$, $J_{\pazocal{C}}$ be the Riemannian metric and complex structure. We denote the quotient map $\mu_s^{-1}(0)\to \pcC$ by $\pi_{\pcC}$.

\begin{proposition}\label{cone2}
The Euler vector field $E_{s}$ induces a vector field $E_{s,\pazocal{C}}$ on $\pazocal{C}$, and the K\"ahler manifold $(\pazocal{C},\mathfrak{g}_{s,\pazocal{C}},J_{\pazocal{C}})$ is a K\"ahler cone with the Euler vector field $E_{s,\pazocal{C}}$.
\end{proposition}
\begin{proof}
First, we prove that $E_s$ is tangent to $\mu_s^{-1}(0)$ and descends.
Recall $\mu_s(z)=-(1+s)u^s w(z)$ with $u=\sum|z_j|^2$ and $w(z)=\sum \mathrm{ch}(z_j)|z_j|^2$.
Let $E_0:=\sum(z_j\partial_{z_j}+\bar z_j\partial_{\bar z_j})$, so $E_s=\frac{1}{1+s}E_0$.
Since $u,w$ are homogeneous of degree $2$,
\[
E_0(u)=2u,\qquad E_0(w)=2w,\qquad E_0(u^s)=2s\,u^s.
\]
Hence
\[
E_0(\mu_s)=-(1+s)\bigl(E_0(u^s)w+u^sE_0(w)\bigr)
=-(1+s)\bigl(2s\,u^s w+2u^s w\bigr)
=2(1+s)\mu_s,
\]
so $E_s(\mu_s)=2\mu_s$. In particular, on $\mu_s^{-1}(0)$ we have $E_s(\mu_s)=0$, i.e. $E_s$ is tangent.

Moreover, the $U(1)$-action is linear and commutes with scalar multiplication, hence preserves $E_s$.
Therefore, $E_s$ is $\pi$-projectable and defines a unique vector field $E_{s,\pcC}$ on $\pcC$ by
$d\pi(E_s)=E_{s,\pcC}\circ\pi$.

Second, we prove that $E_{s,\pcC}$ is complete.
The flow of $E_s$ is $\varphi_t(z)=e^{t/(1+s)}z$, it preserves $\mu_s^{-1}(0)$ (since $E_s(\mu_s)=2\mu_s$),
and it commutes with $U(1)$. Hence, it descends to a flow on $\pcC$ defined for all $t$, so $E_{s,\pcC}$ is complete.

Third, we prove the equalities $\mathcal L_{E_{s,\pcC}}g_{s,\pcC}=2g_{s,\pcC}$ and $\mathcal L_{E_{s,\pcC}}J_\pcC=0$.
By K\"ahler reduction, $g_{s,\pcC},J_C$ are induced from $(\widetilde g_s,\widetilde J)$ via the Riemannian
submersion $\pi_\pcC$. Since $E_s$ and $E_{s,\pcC}$ are $\pi_\pcC$-related, for any tensor $T$ on $\pcC$,
$\pi_\pcC^*(\mathcal L_{E_{s,\pcC}}T)=\mathcal L_{E_s}(\pi_\pcC^*T)$.
Using Proposition~\ref{prop:KahlerCone1}, $\mathcal L_{E_s}\widetilde g_s=2\widetilde g_s$ and $\mathcal L_{E_s}\widetilde J=0$,
and these identities restrict to the horizontal distribution defining $\pi_\pcC^*g_{s,\pcC}$ and $\pi_\pcC^*J_\pcC$.
Therefore
\[
\pi_\pcC^*(\mathcal L_{E_{s,\pcC}}g_{s,\pcC})=2\pi^*g_{s,\pcC},\qquad
\pi_\pcC^*(\mathcal L_{E_{s,\pcC}}J_\pcC)=0.
\]
Since $\pi_\pcC$ is a submersion, $\pi_\pcC^*$ is injective on forms/tensors, hence
$\mathcal L_{E_{s,\pcC}}g_{s,\pcC}=2g_{s,\pcC}$ and $\mathcal L_{E_{s,\pcC}}J_\pcC=0$.

Finally, we only have to prove $d\,g_{s,\pcC}(E_{s,\pcC},-)=0$.
Let $\theta_\pcC:=g_{s,\pcC}(E_{s,\pcC},-)$. On $\mu_s^{-1}(0)$,
\[
\pi_\pcC^*\theta_\pcC=\widetilde g_s(E_s,-).
\]
Indeed, $E_s$ is orthogonal to the vertical direction: since $u$ is $U(1)$-invariant, $du(\hat\xi)=0$ for any fundamental vector field $\hat\xi$, and by Proposition~\ref{prop:KahlerCone1} we have
$\widetilde g_s(E_s,-)=d(u^{1+s})$, hence $\widetilde g_s(E_s,\hat\xi)=d(u^{1+s})(\hat\xi)=0$.
Thus, $E_s$ is horizontal, so the pullback identity holds.
But $\widetilde g_s(E_s,-)=d(u^{1+s})$ is closed, hence $d(\pi_\pcC^*\theta_\pcC)=0$, i.e. $\pi_\pcC^*(d\theta_\pcC)=0$.
Injectivity of $\pi_\pcC^*$ gives $d\theta_\pcC=0$.

Therefore, $(\pcC,g_{s,\pcC},J_\pcC)$ satisfies \eqref{eq:KahlerCone} with Euler field $E_{s,\pcC}$, so it is a K\"ahler cone.
\end{proof}

The radial coordinate $\rho_{s,\pcC}$ on $\pcC$ is given by
\begin{equation}\label{eq:radial_coordinate}
\rho_{s,\pcC} = \sqrt{\widetilde{g}_s(E_s,E_s)}=\sqrt{\frac{1}{(1+s)^2}\sum_{j,k}2\widetilde{g}_s(z_j\partial_{z_j},\overline z_k \partial_{\overline{z}_k})}=\sqrt{2}(u|_{\mu_s^{-1}(0)})^{\frac{1+s}{2}}.
\end{equation}

We claim that the Riemannian manifold $(\cX,\mathfrak{g}_s)$ is an asymptotically conical manifold associated to the Riemannian cone $(\pazocal{C},\mathfrak{g}_{s,\pazocal{C}})$. We first need to show that $(\cX,\mathfrak{g}_s)$ is geodesically complete.

\begin{proposition}
The Riemannian manifold $(\cX,\mathfrak{g}_s)$ is geodesically complete.
\end{proposition}
\begin{proof}
Since $(\cX,\mathfrak{g}_s)$ is given a quotient metric from $(\mu^{-1}_s(0),{g}_s|_{\mu^{-1}_s(0)})$, we need to show that the Riemannian manifold $(\mu^{-1}_s(0),{g}_s|_{\mu^{-1}_s(0)})$ is geodesically complete. We use the Hopf-Rinow theorem \cite{HR}: closed bounded subsets are compact if and only if the Riemannian manifold is geodesically complete.

Let $C$ be a closed and bounded subset of $(\mu^{-1}_s(0),{g}_s|_{\mu^{-1}_s(0)})$. We let $P_0,P_1\in C$ and let $\gamma:I\to \mu^{-1}_s(0)$ be the shortest geodesic such that $\gamma(0)=P_0$, $\gamma(1)=P_1$. The length of this curve is the distance between $P_0$ and $P_1$. In $\mu^{-1}_s(0)$, we have the inequality $\rho_0\geq1$. The metric $\widetilde{g}_s$ can be represented by the matrix in \eqref{metricmatrix}. Hence, for any tangent vector $X\in T_{P}\mu^{-1}_s(0)$, we have
$$
\widetilde{g}_s(X,X) \geq g_{\mathrm{Euc}}(X,X)
$$
where $g_{\mathrm{Euc}}$ is the Euclidean metric on $\bC^{n+r}$. Hence, the length of the path $\gamma$ with respect to the Euclidean metric is shorter than that with respect to the metric $\widehat{g}_s$. Moreover, the Euclidean distance $d_{\mathrm{Euc}}(P_0,P_1)$ between $P_0$ and $P_1$ is shorter than the length of the path $\gamma$ with respect to either of the metrics. This is true for any pair $(P_0,P_1)$ of points in $C$, so $C$ is bounded in the Euclidean metric as well. Thus, $C$ is closed and bounded in the Euclidean metric and therefore compact.
\end{proof}

\begin{proposition}\label{prop:diffeo}
There is a diffeomorphism
\[
\Phi:\pazocal{C} \to \cX\setminus \{p_1=\ldots=p_r=0\}
\]
defined by
\[
\Phi([z]_{U(1)}) := [z]_{\bC^*}
\]
where $z\in \mu^{-1}_s(0)$, and $[z]_{U(1)}$ is the $U(1)$-orbit of $z$.
\end{proposition}
\begin{proof}
Recall that from Proposition \ref{biholomorphism}, we have a biholomorphism $\mu_s^{-1}(-1)/U(1)\cong \cX$. We construct a $U(1)$-equivariant map \[\widetilde{\Phi}:\mu^{-1}_s(0)\to \mu^{-1}_s(-1)\setminus \{p_1=\cdots=p_r=0\},\]
and it naturally induces $\Phi:\pcC\to \cX\setminus \{p_1=\cdots=p_r=0\}$ satisfying the commutative diagram
\[\begin{tikzcd}
\mu_s^{-1}(0) \arrow[r, "\widetilde{\Phi}"] \arrow[d, "\pi_{\pcC}"'] &
\mu_s^{-1}(-1)\setminus\{p_1=\cdots=p_r=0\} \arrow[d, "\pi"] \\
\pcC=\mu_s^{-1}(0)/U(1) \arrow[r, "\Phi"] &
X=\mu_s^{-1}(-1)/U(1).
\end{tikzcd}
\] 
Using the analytic inverse function theorem, we have a smooth function $\lambda:\mu^{-1}_s(0) \to \bR$ such that
\[
\mu_s(e^{\lambda(z)}\cdot z) = -1.
\]
Since $\mu_s$ is $U(1)$-invariant, $\lambda$ is also $U(1)$-invariant. We define $\widetilde{\Phi}:\mu^{-1}_s(0)\to \mu^{-1}_s(-1)$ as 
\[
\Phi(z) = e^{\lambda(z)}\cdot z.
\]
Since $\bC^*$ is an abelian Lie group, $\widetilde{\Phi}$ is $U(1)$-equivariant, so the proof follows.
\end{proof}

We introduce symbols $\overline{1},\ldots,\overline{n+r}$ to denote $\overline{Z}_{k}$ and $\overline z_{k}$ by $Z_{\overline{k}}$ and ${z_{\overline k}}$. For simplicity, we write
\[
h := \Phi^* \mathfrak{g}_s - \mathfrak{g}_{s,\pcC}.
\]
Now, our theorem is as follows.


\begin{theorem} \label{keylemma}
The K\"ahler manifold $(\cX,\mathfrak{g}_s,J)$ is asymptotically conical along the diffeomorphism $\Phi:\pcC \to \cX\setminus\{p_1=\cdots=p_r=0\}$ with decay rate $2$, i.e.
for each $k\geq 0$,
\begin{align} \label{eq:goal_equation}
|(\nabla^{\pcC})^k h|_{\mathfrak{g}_{s,\pcC}} = O(\rho_{s,\pcC}^{-2-k}).
\end{align}
\end{theorem}
\begin{proof}
The definition of the diffeomorphism $\Phi$ is given in Proposition \ref{prop:diffeo}. In Lemma \ref{prop:shrinking}, we define an open cover $\{V_{j}\}_{1\leq j \leq n}$ of $\pcC$ as
\[
V_j := \left\{[z]_{U(1)}\in \pcC: \mu_s(z)=0,\ |x_j|>\frac{\sqrt{u(z)}}{\sqrt{8n}}\right\}.
\]

We prove that \eqref{eq:goal_equation} is true in each $V_j$; once the open cover $\{V_j\}$ is in place, the remainder of the proof is reduced to a chartwise computation. 
In local coordinates, one checks that Christoffel symbols $\Gamma^{m}_{kl;j}$ and error terms $h_{k\overline{l};j}$,
together with all iterated derivatives of arbitrary order, can be written as rational functions (in the relevant local variables,
with explicitly controlled denominators). 
This ``rational'' structure is precisely what makes the asymptotic analysis feasible: once the local formulas are established,
the desired decay estimates follow by straightforward bookkeeping. We choose a holomorphic coordinate chart $Z_{1/j},\ldots,Z_{n+r/j}$ defined as
\[
Z_{k/j} = x_1^{-\ch(z_k)}z_k.
\]
Then locally, the symmetric bilinear two form $h$ is written as
\[
h = \sum_{k,l} h_{k\overline{l};j} dZ_{k/j}\odot d\overline{Z}_{l/j}, 
\]
where $v\odot w := v\otimes w + w\otimes v$. Then by Leibniz's rule, the covariant derivatives $(\nabla^{\pcC})^m h$ is written as a linear combination of
\begin{align} \label{eq:goal2}
\partial_{Z_{j_1}}\cdots \partial_{Z_{j_{m'}}} h_{k\overline{l};j}\; (\nabla^{\pcC})^{q_1} dZ_{j_1/j}\otimes \cdots \otimes (\nabla^{\pcC})^{q_{m'}}dZ_{j_{m'}/j} \otimes ((\nabla^{\pcC})^{q_{m'+1}}dZ_{k/j}\odot (\nabla^{\pcC})^{q_{m'+2}}d\overline{Z}_{l/j})
\end{align}
where $j_1,\ldots,j_{m'}\in\{2,\overline{2},\ldots,n+r,\overline{n+r}\}$, $m'+q_1+\cdots + q_{m'+2}=m$. Hence, by the triangle inequality, the norm of the covariant derivatives $(\nabla^{\pcC})^m h$ is controlled by the norm of \eqref{eq:goal2}. In Lemma \ref{cor:cov_asymp}, we will prove that the following is true for each $l=1,2,\ldots$, $j=1,\ldots,n$, and $k=1,\ldots,j-1,j+1,\ldots,n+r$:
\[
|(\nabla^{\pcC})^l dZ_{k/j}|_{\mathfrak{g}_{s,\pcC}} = O(\rho_{s,\pcC}^{-l} |dZ_{j}|_{\mathfrak{g}_{s,\pcC}}).
\]
Hence, we have
\begin{align} 
|\partial_{Z_{j_1/j}}\cdots \partial_{Z_{j_{m'}/j}} h_{k\overline{l};j}\; (\nabla^{\pcC})^{q_1} dZ_{j_1/j}\otimes \cdots \otimes (\nabla^{\pcC})^{q_{m'}}dZ_{j_{m'}/j} \otimes ((\nabla^{\pcC})^{q_{m'+1}}dZ_{k/j}\odot (\nabla^{\pcC})^{q_{m'+2}}d\overline{Z}_{l/j})|_{\mathfrak{g}_{s,\pcC}}\nonumber\\
=O(\rho_{s,\pcC}^{-q_1-\cdots-q_{m'+2}})|\partial_{Z_{j_1}}\cdots\partial_{Z_{j_{m'}}} h_{k\overline{l};j}dZ_{j_1}\otimes \cdots \otimes dZ_{j_{m'}}\otimes dZ_j \otimes d\overline{Z}_k|_{\mathfrak{g}_{s,\pcC}}.\label{eq:goal3}
\end{align}
In Lemma \ref{cor:h_diff_asymp}, we will prove that for $m'\geq 0$, $j_1,\ldots,j_{m'}\in \{2,\overline{2},\ldots,n+r,\overline{n+r}\}$, $j=1,\ldots,n$, and $k,l\in \{1,2,\ldots,j-1,j+r,\ldots,n+r\}$,
\[
|\partial_{Z_{j_1}}\cdots\partial_{Z_{j_m'}} h_{k\overline{l};j}dZ_{j_1/j}\otimes \cdots \otimes dZ_{j_{m'}/j}\otimes dZ_{k/j} \otimes d\overline{Z}_{l/j}|_{\mathfrak{g}_{s,\pcC}} = O(\rho_{s,\pcC}^{-2-m'}).
\]
Hence, \eqref{eq:goal3} becomes
\begin{align*}
O(\rho_{s,\pcC}^{-q_1-\cdots-q_{m'+2}})|\partial_{Z_{j_1}}\cdots\partial_{Z_{j_{m'}}} h_{k\overline{l};j}dZ_{j_1}\otimes \cdots \otimes dZ_{j_{m'}}\otimes dZ_j \otimes d\overline{Z}_k|_{\mathfrak{g}_{s,\pcC}}\\
=O(\rho_{s,\pcC}^{-q_1-\cdots-q_{m'+2}})O(\rho_{s,\pcC}^{-2-m'})=O(\rho_{s,\pcC}^{-2-m}).
\end{align*}
Thus, it remains to prove Lemmas \ref{prop:shrinking}, \ref{cor:cov_asymp}, and \ref{cor:h_diff_asymp}.
\end{proof}

Using \eqref{eq:radial_coordinate}, we will compute the asymptotic estimate with respect to $u|_{\mu^{-1}_s(0)}$ instead of $\rho_{s,\pcC}$, since we have a useful fact that
\[
|x_j| = O(\sqrt{u})
\]
on $\pcU$. For example, \eqref{eq:goal_equation} is equivalent to
\[
|(\nabla^{\pcC})^k h|_{\mathfrak{g}_{s,\pcC}} = O\left((u|_{\mu^{-1}_s(0)})^{-(2+k)\frac{1+s}{2}}\right).
\]

Our strategy is to set an open cover with explicit coordinate functions and prove that Theorem \ref{keylemma} is true for each open subset. For each $j=1,\ldots,n$, we define $U_{j,0}\subset \pcC$ as
\[
U_{j,0} = \{[z_1,\ldots,z_{n+r}]_{U(1)}\in \pcC: x_j \neq 0\}.
\]
Then $\{U_{j,0}\}_{1\leq j \le n}$ is an open cover of $\pcC$, and it has a holomorphic coordinate system $\bC^{n-1}\times (\bC^r\setminus \{0\})\cong U_{j,0}$ given by
\[
\bC^{n-1}\times (\bC^r\setminus \{0\})\cong U_{j,0},\quad (Z_{1/j},\ldots,Z_{j-1/j},Z_{j+1/j},\ldots,Z_{n+r/j})\mapsto (z_1,\ldots,z_{n+r}),\quad Z_{k/j} = x_j^{-\ch(z_k)}z_k.
\]
However, this open cover does not behave well; in our computation of $h_{k\overline{l};j}$ and its partial derivatives on $U_{j,0}$, negative powers of $x_j$ appears, but $|x_j|$ has no lower bounds on $U_{j,0}$. Thus, we want smaller open subsets $V_j\subset U_{j,0}$ so that they cover $\pcC$ and $|x_j|\geq c\sqrt{u}$ in $V_j$ for some constant $c>0$.

On $\pcC$, when $f,g$ are two smooth real valued functions defined on a noncompact open subset of $\pcC$, we write 
\[
f\asymp g
\]
if there is a constant $c_1$ and $c_2$ such that when the radial coordinate $\rho_{s,\pcC}(z)$ is large enough,
\[
c_1f(z)< g(z)< c_2f(z).
\]

\begin{lemma}[Shrinking of the open cover $\{U_{j,0}\}_{1\leq j \leq n}$]
\label{prop:shrinking} \label{cor3.25}
We define $V_{j}\subset U_{j,0}$ as
\[
V_j := \left\{[z]_{U(1)}\in \pcC: \mu_s(z)=0,\ |x_j|>\frac{\sqrt{u(z)}}{\sqrt{8n}}\right\}
\]
for each $j=1,\ldots,n$. Then, $\{V_j\}_{1\leq j \leq n}$ is an open cover of $\pcC$. In particular, $|x_j|\asymp \sqrt{u|_{\mu^{-1}_s(0)}}$ on $V_j$.
\end{lemma}
\begin{proof}
Choose $(x_1,\ldots,x_n,p_1,\ldots,p_r)\in \pcC$ arbitrarily. Since $(x_1,\ldots,x_n)\neq 0$, $\max(|x_1|,\ldots,|x_n|) \neq 0$. Suppose that for some $j\in \{1,\ldots,n\}$, $|x_j|= \max(|x_1|,\ldots,|x_n|)$. Then obviously,
\[
|x_j|> \frac{1}{2}\max(|x_1|,\ldots,|x_n|),
\]
and
\[
\max(|x_1|,\ldots,|x_n|) \geq \frac{\sqrt{|x_1|^2+\cdots +|x_n|^2}}{\sqrt{n}}.
\]
Since $\mu_s(z)=0$, we have
\[
|x_1|^2+\cdots +|x_n|^2 = d_1 |p_1|^2 +\cdots +d_r|p_r|^2,
\]
so
\[
\sqrt{|x_1|^2+\cdots +|x_n|^2} = \frac{1}{\sqrt{2}} \sqrt{|x_1|^2+\cdots +|x_n|^2+d_1|p_1|^2+\cdots +d_r|p_r|^2}\geq \frac{\sqrt{u(z)}}{\sqrt{2}}.
\]
Hence,
\[
|x_j|> \frac{\sqrt{u(z)}}{\sqrt{8n}},
\]
so $[z]_{U(1)}\in V_j$.
\end{proof}

\begin{remark}
In Lemma \ref{prop:shrinking}, the constant $\sqrt{8n}$ is arbitrarily chosen; for any $0<\e<1$, we define $V_j^\e$ as
\[
V_j^\e = \left\{[z]_{U(1)}\in \pcC:\mu_s(z)=0, |x_j|> \frac{\e}{\sqrt{2n}}\sqrt{u(z)}\right\},
\]
and $\{V_j^\e\}_{1\leq j\leq n}$ is still an open cover of $\pcC$. In this paper, we chose $\e=1/2$.
\end{remark}

Using the analytic implicit function theorem, we have two smooth functions $\lambda_j,\eta_{j}:V_j\to \bR$ such that for $Z_{\bullet/j}=(Z_{1/j},Z_{2/j},\ldots,Z_{j-1/j},1,Z_{j+1/j},\dots,Z_{n+r/j})\in V_j$,
\[
\mu_s(e^{\lambda_j(Z_{\bullet/j})}\cdot Z_{\bullet/j}) = -1, \quad \mu_s(e^{\eta_{j}(Z_{\bullet/j})}\cdot Z_{\bullet/j}) = 0.
\]
That is, $e^{\lambda_j(Z_{\bullet/j})}\cdot Z_{\bullet/j}$ (or $e^{\eta_{j}(Z_{\bullet/j})}\cdot Z_{\bullet/j}$) is a lift of $Z_{\bullet/j}\in V_j\subset \pcC$ onto the moment level $\mu^{-1}_s(-1)$ (or $\mu^{-1}_s(0)$, respectively.) 

From now on, we fix $j=1$; the proof is similar for other $j$'s. Thus, for simplicity, we drop the index $j$ and write $Z_k := Z_{k/1}$, $Z=(1,Z_2,\ldots,Z_{n+r})\in V_1$. Also, we define $\widehat{z}_k:V_1\to \bC$ as 
\begin{align*} 
\widehat{z}_k(Z) := e^{\ch(z_k)\eta_1(Z)}Z_k
\end{align*}
so that $(\widehat{z}_1,\cdots,\widehat{z}_{n+r})\in \mu^{-1}_s(0)$ and $(\widehat{z}_1,\cdots,\widehat{z}_{n+r})\mapsto(1,Z_2,\ldots,Z_{n+r})\in V_1$ via the quotient map $\mu^{-1}_s(0)\to \pcC$. On $V_1$, we let
\begin{align*} 
\widehat{u}(Z) := u(\widehat{z}(Z)) = u (e^{\eta_1(Z)}\cdot Z).
\end{align*}
The following proposition tells us why Proposition~\ref{prop:shrinking} is useful.

\begin{proposition}\label{prop:dZ_asymp}
On $V_1$,
\[
|dZ_j|_{g_{s,\pazocal{C}}}\asymp \widehat u^{-\frac12({\ch(z_j)+s})}.
\]
\end{proposition}

\begin{proof}
Let $\pi_{\pcC}:\mu_s^{-1}(0)\to\pazocal{C}$ be the quotient
map. Since each $Z_j$ is $U(1)$-invariant, the $1$-form $dZ_j$ is basic on $\mu_s^{-1}(0)$; hence its norm with
respect to the quotient metric equals the norm computed upstairs:
\[
|dZ_j|_{g_{s,\pazocal{C}}}\circ \pi_{\pcC} \;=\; |d(x_1^{-\ch(z_j)}z_j)|_{\widetilde g_s|_{\mu_s^{-1}(0)}} .
\]
We use the explicit inverse matrix of $(\widetilde g_s)_{j \overline k}=\partial_{z_j}\partial_{\overline{z}_k}(u^{1+s})$ (see \eqref{metricmatrix}) on $\pazocal{U}$:
\[
(\widetilde g_s)^{j\bar k}
=\frac{1}{(1+s)u^s}\delta_{jk}-\frac{s}{(1+s)^2u^{s+1}}\bar z_j z_k .
\]
In particular there are constants $c_1,c_2>0$ (independent of $u$ for $u\gg 1$) such that
\[
c_1\,u^{-s}\le (\widetilde g_s)^{a\bar a}\le c_2\,u^{-s},\qquad a=1,\dots,n+r.
\]
Differentiate $Z_j=x_1^{-\ch(z_j)}z_j$:
\[
dZ_j=x_1^{-\ch(z_j)}dz_j-\ch(z_j)\,x_1^{-\ch(z_j)-1}z_j\,dx_1 .
\]
Thus on $V_1$, we have
\[
dZ_j=\widehat x_1^{-\ch(z_j)}dz_j-\ch(z_j)\,\widehat x_1^{-\ch(z_j)-1}\widehat z_j\,dx_1 .
\]
Therefore, using $(\widetilde g_s)^{j\bar j}\le c_2 u^{-s}$ and $|\widehat z_j|^2\le u$, we get
\[
|dZ_j|^2_{\widetilde g_s}
\le c_2 \widehat u^{-s}\Bigl(|\widehat x_1|^{-2\ch(z_j)}+|\ch(z_j)|^2|\widehat x_1|^{-2\ch(z_j)-2}|\widehat z_j|^2\Bigr)
\le c_2 \widehat u^{-s}|\widehat x_1|^{-2\ch(z_j)}\Bigl(1+|\ch(z_j)|^2\frac{\widehat u}{|\widehat x_1|^2}\Bigr).
\]
On $V_1$ we have $|\widehat{x}_1|>\sqrt{\widehat u}/\sqrt{8n}$ by definition (Proposition~\ref{prop:shrinking}), hence
$\frac{\widehat u}{|\widehat x_1|^2}\le 8n$, and also \[|\widehat x_1|^{-2\ch(z_j)}\le c_2 \widehat u^{-\ch(z_j)}.\] Thus
\[
|dZ_j|^2_{\widetilde g_s}\le c_2(1+8|\ch(z_j)|^2n)\,\widehat u^{-s-\ch(z_j)},
\qquad\text{i.e.}\qquad
|dZ_j|_{g_{s,\pazocal{C}}}=O\!\left(\widehat u^{-\frac{\ch(z_j)+s}{2}}\right).
\]

For the second estimate, let $\alpha=\sum_{a=1}^{n+r}\alpha_a\,dz_a$ be a $(1,0)$--form on $\pazocal{U}$.
Using the explicit inverse matrix
\[
(\widetilde g_s)^{a\bar b}
=\frac{1}{(1+s)u^s}\delta_{ab}
-\frac{s}{(1+s)^2u^{s+1}}\bar z_a z_b,
\]
we compute
\begin{align*}
|\alpha|_{\widetilde g_s}^2
&=\sum_{a,b}(\widetilde g_s)^{a\bar b}\alpha_a\overline{\alpha_b} \\
&=\frac{1}{(1+s)u^s}\sum_{a}|\alpha_a|^2
-\frac{s}{(1+s)^2u^{s+1}}\Bigl|\sum_a \alpha_a\bar z_a\Bigr|^2.
\end{align*}
By the Cauchy-Schwarz inequality,
\[
\Bigl|\sum_a \alpha_a\bar z_a\Bigr|^2
\le \Bigl(\sum_a|\alpha_a|^2\Bigr)\Bigl(\sum_a|z_a|^2\Bigr)
=u\sum_a|\alpha_a|^2,
\]
hence
\[
|\alpha|_{\widetilde g_s}^2
\ge \left(\frac{1}{(1+s)u^s}-\frac{s}{(1+s)^2u^{s+1}}\cdot u\right)\sum_a|\alpha_a|^2
=\frac{1}{(1+s)^2}\,u^{-s}\sum_a|\alpha_a|^2.
\]
In particular,
\begin{equation}\label{eq:alpha-lower}
|\alpha|_{\widetilde g_s}^2 \ge \frac{1}{(1+s)^2}\,u^{-s}\,|\alpha_j|^2
\qquad (j=1,\dots,n+r).
\end{equation}

Applying \eqref{eq:alpha-lower} to $\alpha=dZ_j$ (computed first on $\pazocal{U}$ as
$dZ_j=x_1^{-\ch(z_j)}dz_j-\ch(z_j)x_1^{-\ch(z_j)-1}z_j\,dx_1$) and then restricting to $\mu_s^{-1}(0)$,
we obtain
\[
|d Z_j|_{\widetilde g_s}^2 \ge \frac{1}{(1+s)^2}\,\widehat u^{-s}\,|\widehat x_1|^{-2\ch(z_j)}.
\]
Since $|\widehat x_1|\asymp \sqrt{\widehat u}$ by Proposition~\ref{prop:shrinking}, we have the desired inequality.
\end{proof}

We now consider the asymptotic estimate of the Christoffel symbols of the Levi-Civita connection $\nabla^{\pcC}$ on $V_1$. We define a function $w_j:\pcU\to \bR$ as 
\begin{align}\label{eq:def_of_w2}
w_j = \sum_{k=1}^{n+r}\ch(z_k)^j |z_k|^2
\end{align}
for each $j=0,1,\ldots.$ In particular, $u=w_0$ and $w=w_1$. On $V_1$, we let
\[
\widehat{w}_j(Z) := w_j (\widehat{z}(Z)) = w_j (e^{\eta_1(Z)}\cdot Z). 
\]
Then we have $\widehat{w}_1\equiv 0$. 

We let $\Gamma^l_{jk}$ be the Christoffel symbol of $\nabla^{\pcC}$ on $V_1$, i.e.
\[
\nabla^{\pcC}_{\partial_{Z_j}}\partial_{Z_k} = \sum_{l=2}^{n+r} \Gamma^l_{jk} \partial_{Z_l}.
\]
Since $Z$ is a holomorphic coordinate on $V_1$, for $j,k,l\in \{2,\ldots,n+r\}$, $\Gamma^{\overline{l}}_{jk} = \Gamma^{l}_{\overline jk}=0$ and $\overline{\Gamma^l_{jk}} = \Gamma^{\overline{l}}_{\overline{jk}}$. We let ${}^{\pcU}{\Gamma}^j_{kl}$ be the Christoffel symbol of the Levi-Civita connection $\nabla^{\pcU}$ of the metric $\widetilde{g}_s$. Then,
\begin{align} \label{eq:formula_chstf}
\!^{\pcU}{\Gamma}^j_{kl}(z) = s\frac{\overline{z}_k \delta_{jl}+\overline{z}_l \delta_{jk}}{u}- s\frac{\overline{z}_k \overline z_{l}z_j}{u^2}.
\end{align}
We let $\widehat\Gamma^j_{kl}$ as a function on $V_1$ defined as
\begin{align} \label{eq:formula_chstf2}
\widehat\Gamma^j_{kl}:=\!^{\pcU}{\Gamma}^j_{kl}(\widehat{z})=s\frac{\widehat{\overline{z}}_k \delta_{jl}+\widehat{\overline{z}}_l \delta_{jk}}{u}- s\frac{\widehat{\overline{z}}_k \widehat{\overline z}_{l}\widehat{z}_j}{u^2}.
\end{align}
For each $j=1,\dots,n+r$, we define $\pcV_j:\pcU\to \bC$ as
\[
\pcV_j(z) = \frac{\widetilde{g}_s(\partial_{z_j},\overline{\mathbf{n}})}{\widetilde{g}_s(\mathbf{n},\overline{\mathbf{n}})} = \overline{z}_k\frac{\ch(z_k)u+sw}{w_2 u+ sw^2},
\]
see \eqref{eq:def_of_n} for the definition of $\mathbf{n}$ and \eqref{eq:def_of_w2} for the definition of $w_2$. In particular, we have
\begin{align} \label{eq:formula_H}
\cH \partial_{z_j} = \partial_{z_j} - \pcV_j \mathbf{n}.
\end{align}
On $V_1$, we have
\[
\widehat{\pcV}_j:=\pcV_j(\widehat{z}) = \frac{\ch(z_k) \widehat{\overline{z}}_k}{\widehat w_2}
\]
where $\cH$ is the horizontal projection defined in \eqref{eq:def_of_H}.

{
To systematically bound the Christoffel symbols and metric error terms, the analysis must move away from the singular locus defined by the condition $x_j=0$. This is achieved by shrinking the initial open cover to localized charts $V_j$ where the coordinate lower bounds $|x_j| \asymp \sqrt{u}$ hold. Within these refined charts, the coordinates are lifted to the moment levels to define localized variables $\widehat{z}$ and $\widetilde{z}$, alongside their corresponding radial norms $\widehat{u}$ and $\widetilde{u}$. This localization isolates the dominating terms with respect to $\widehat{u}$, enabling the expression of the Levi-Civita connection and the K\"ahler metric components as rational functions with explicitly controlled denominators. We summarize in Table~\ref{tab:asymptotic-dict1} the local quantities used in Lemma~\ref{cor:cov_asymp} and their asymptotics with respect to $\widehat u$ on the chart $V_1$.  The purpose of this table is twofold: it fixes notation and provides uniform scale estimates for the basic building blocks.  
With these estimates in place, the proof reduces to systematic bookkeeping, since higher-order derivatives are obtained by iterating differentiation rules on rational expressions.
}

\begin{table}[h]
  \centering
  \caption{Asymptotics of Functions and Tensors on $V_1$.}
  \label{tab:asymptotic-dict1}

  \renewcommand{\arraystretch}{1.3}
  \setlength{\tabcolsep}{6pt}

  \begin{tabularx}{\linewidth}{@{}  l | l | X | r @{}}
    \toprule
      & \textbf{Definition} & \textbf{Asymptotics on $V_1$} & \textbf{Reference} \\
    \midrule    $\rho_{s,\pcC}$  & $\sqrt {\widetilde{g}_s(E_s,E_s)}$. & $\sqrt{2}\widehat{u}^{\frac{1+s}{2}}$. & \eqref{eq:radial_coordinate}  \\
    $\widehat{x}_1$  & $e^{\eta_1(Z)}$. & $|x_1|\asymp \widehat{u}^{1/2}.$ & Prop \ref{prop:shrinking}.  \\    $\widehat{z}_j$  & $e^{\ch(z_j)\eta_1(Z)}Z_j$. & $|\widehat{z}_j|=O( \widehat{u}^{1/2})$. & By def.  \\
    $\widehat{w}_{2j}$  & $\sum_k \ch(z_k)^{2j}|\widehat{z}_k|^2$. & $\widehat{w}_{2j}\asymp {\widehat{u}}$. & By def.\\
    $\widehat{\pcV}_j$  & $\ch(z_j)\widehat{\overline{z}}_j/\widehat{w}_2$. & $|\widehat{\pcV}_j|=O(\widehat{u}^{-\frac{1}{2}})$. & By def.\\
    $\widehat{\Gamma}^j_{jl}$  & ${}^\pcU\Gamma^j_{kl}(\widehat{z})$. & $|\widehat{\Gamma}^j_{jl}|=O(\widehat{u}^{-\frac{1}{2}})$. &  \eqref{eq:formula_chstf2}.\\
    $dZ_j$  & Differential forms on $V_1$. & 
    $|dZ_j|_{\mathfrak{g}_{s,\pcC}}=O(\widehat{u}^{\frac{1}{2}(-\ch(z_j)-s)})$. & Prop \ref{prop:dZ_asymp}.\\
    $\Gamma^j_{kl}$  & Christoffel symbols of $\nabla^\pcC$. & $|\Gamma^j_{kl}|=O(\widehat{u}^{\frac{1}{2}(\ch(z_k)+\ch(z_l)-\ch(z_j)-1)})$. & Prop \ref{prop:formula_of_Gamma}.\\
    $\partial_{Z_m}\Gamma^j_{kl}$  & Partial derivative of $\Gamma^j_{kl}$.& $|\partial_{Z_m}\Gamma^j_{kl}| =O(\widehat{u}^{\frac{1}{2}(\ch(z_m)+\ch(z_k)+\ch(z_l)-\ch(z_j)-2)})$. & Prop \ref{differentiation1}.\\
    \bottomrule
  \end{tabularx}
\end{table}

\begin{proposition}\label{prop:formula_of_Gamma} For each $j,k,l=2,\ldots,n+r$, we have
\[
\Gamma^j_{kl} = \widehat{x}_1^{\ch(z_k)+\ch(z_l)-\ch(z_j)}\left(\widehat\Gamma^j_{kl}-\delta_{jl}\ch(z_l)\widehat \pcV_k-\delta_{jk}\ch(z_k)\widehat \pcV_l+(\ch(z_j)^2-\ch(z_j))\widehat \pcV_k\widehat \pcV_l \widehat{z}_j\right).
\]
\end{proposition}
\begin{proof}
Since $\pi|_{\mu^{-1}_s(0)}:\mu^{-1}_s(0) \to \cX$ is a Riemannian submersion, the Levi-Civita connection $\nabla^{\pcC}$ satisfies
\[
\nabla^{\pcC}_{\partial_{Z_k}} \partial_{Z_l} = d\pi|_{\mu^{-1}_s(0)} \nabla_{x_1^{\ch(z_k)}\cH\partial_{z_k}}^{\pcU} x_1^{\ch(z_l)}\cH\partial_{z_l}.
\]
The proof follows from the direct computation using \eqref{eq:formula_chstf2} and \eqref{eq:formula_H}.
\end{proof}

We observe that $\Gamma^j_{kl}$ can be written in the form of a rational function
\[
\widehat{x}_1^{\ch(z_k)+\ch(z_l)-\ch(z_j)}\frac{P(\widehat{z})}{\widehat{u}^2 \widehat{w}^2_2}
\]
where $P$ is a polynomial with degree at most $7$; here the \emph{degree} is the ordinary degree of a polynomial, not the charge or weight. We write
\[
\bC[z,\overline{z}]_{\deg\leq k}
\]
be a set of polynomials in $\bC[z_1,\overline{z}_1,\ldots,z_{n+r},\overline{z}_{n+r}]$ whose degrees are at most $k$. From Proposition \ref{prop:formula_of_Gamma}, we have 
\[
\Gamma^j_{kl}(Z) = P(\widehat{z}) \quad \textrm{for $P\in x_1^{\ch(z_j)+\ch(z_k)-\ch(z_l)} u^{-2} w_2^{-2} \bC[z,\overline{z}]_{\deg\leq 7}$. }
\]
We show that such a rational form is preserved by partial differentiation along $Z_j$.

\begin{definition}\label{def:degree1}
We define $\widehat{R}$ as a $\bC$-algebra of smooth $\bC$-valued functions on $V_1$ such that
\[
\widehat{R} := \{f:V_1\to \bC: f(Z) = P(\widehat{z}(Z)) \textrm{ for $P\in \bC[z,\overline{z}][x_1^{-1},\overline{x}_1^{-1},u^{-1},w_2^{-1}]$}\}.
\]
There is a filtration 
\[\cdots\subset \widehat{R}_{-k} \subset \widehat{R}_{-k+1}\subset\cdots \subset \widehat{R}_{0} \subset \widehat{R}_{1}\subset \cdots \subset \widehat{R}
\]
where
\[
\widehat{R}_{k}:=\{f\in \widehat{R}:f(Z)=P(\widehat{z}(Z)) \textrm{ for $P\in x_1^{\a}\overline{x}_1^\b u^\g w_2^\d \bC[z,\overline{z}]_{\deg \leq l}$ where $\a+\b+2\g+2\d+l\leq k$}\}.
\]
Then 
\[
\bigcup_{k\in \bZ} \widehat{R}_k = \widehat{R}, \qquad \bigcap_{k\in \bZ} \widehat{R}_k =\{0\}.
\]
Then, we have an additive grading $\deg:\widehat{R}\to \bZ \cup \{-\infty\}$ as 
\[
\deg f := \min\{k:f\in \widehat{R}_k\} \textrm{\quad for $f\neq0$},
\]
and $\deg 0 := -\infty.$
\end{definition}

From the fact that $\widehat{w}_2\asymp \widehat{u}$, $|x_1|\asymp \widehat{u}^{1/2}$ by Proposition \ref{prop:shrinking}, and $z_j =O(\widehat{u}^{1/2})$, we have 
\[
|f| = O(\widehat{u}^{\deg f/2})
\]
for $f\in \widehat{R}$. Also, by Proposition \ref{prop:formula_of_Gamma}, we have $\Gamma^j_{kl}\in \widehat{R}_{\ch(z_k)+\ch(z_l)-\ch(z_j)-1}$.

\begin{proposition}\label{differentiation1}
For each $j\in \{2,\overline{2},\ldots,n+r,\overline{n+r}\}$, $\widehat{R}$ is closed under the operator $\partial_{Z_j}$, and $\partial_{Z_j}$ is an operator of degree $\ch(z_j)-1$, i.e.,
\[
\partial_{Z_j} : \widehat{R}_k \to \widehat{R}_{k+\ch(z_j)-1}
\]
is well-defined.
\end{proposition}
\begin{proof}
Since $\eta_1$ satisfies the equation
\[
\mu_s(e^{\eta_1(Z)}\cdot Z)=0,
\]
by differentiating both sides with $Z_j$, we get
\begin{align}
\frac{\partial \eta_1}{\partial Z_j} = -\widehat{x}_1^{\ch(z_j)} \frac{\ch(z_j)\widehat{\overline{z}}_j}{2\widehat{w}_2}
\end{align}
for each $j=2,\overline{2},\ldots,n+r,\overline{n+r}$. By definition of $\widehat{z}_k$,
\[
\partial_{Z_j}\widehat{z}_k = (\partial_{Z_j}\eta_1)\ch(z_k)\widehat{z}_k + \widehat{x}_1^{\ch(z_k)}\delta_{jk} = \widehat{x}_1^{\ch(z_j)} \left(-\frac{\ch(z_j)\ch(z_k)\widehat{\overline{z}}_j\widehat{{z}}_k}{2\widehat{w}_2}+\delta_{jk}\right).
\]

By the chain rule, for $P\in \bC[z,\overline{z}][z_1^{-1},\overline{z}_1^{-1},u^{-1},w_2^{-1}]$,
\[
{\partial}_{Z_j} P(\widehat{z}(Z)) = \widehat{x}_1^{\ch(z_j)}\sum_{k=1}^{n+r} \left( \left(-\frac{\ch(z_j)\ch(z_k){\widehat {\overline{z}}}_j{\widehat{z}}_k}{2\widehat{w}_2}+\delta_{jk}\right) (\partial_{z_k}P(\widehat z))+ \left(-\frac{\ch(z_j)\ch(z_k){\widehat{\overline{z}}}_j{\widehat{\overline{z}}}_k}{2\widehat{w}_2}\right) (\partial_{\overline z_k}P(\widehat{z}))\right).
\]
Then, ${\partial}_{Z_j}(P(\widehat{z}))\in \widehat{R}_{\ch(z_j)+k-1}$ if $P(\widehat{z})\in \widehat{R}_{k}$.
\end{proof}

\begin{lemma} \label{cor:cov_asymp}
For each $k\geq 0$,
\[
|(\nabla^{\pcC})^k dZ_j|_{\mathfrak{g}_{s,\pcC}} = O(u^{-k\frac{1+s}{2}} |dZ_j|_{\mathfrak{g}_{s,\pcC}}).
\]
\end{lemma}
\begin{proof}
We use induction on $k$. When $k=1$, this is given by Proposition \ref{prop:formula_of_Gamma} since 
\[\nabla^{\pcC} dZ_j = \sum_{k,l} \Gamma^j_{kl} dZ_k \otimes dZ_l,\]
and by triangle inequality and $|\widehat{x}_1|\asymp \widehat{u}^{1/2}$ (Proposition \ref{prop:shrinking}),
\[
|\nabla^{\pcC} dZ_j |_{\mathfrak{g}_{s,\pcC}} \leq \sum_{k,l} |\Gamma^j_{kl} dZ_k \otimes dZ_l|_{\mathfrak{g}_{s,\pcC}} = \sum_{k,l}O(u^{\frac{1}{2}(-1+\ch(z_k)+\ch(z_l)-\ch(z_j))}|dZ_k|_{\mathfrak{g}_{s,\pcC}}|dZ_l|_{\mathfrak{g}_{s,\pcC}}),
\]
and by Proposition \ref{prop:dZ_asymp},
\[
O(u^{\frac{1}{2}(-1+\ch(z_k)+\ch(z_l)-\ch(z_j))}|dZ_k|_{\mathfrak{g}_{s,\pcC}}|dZ_l|_{\mathfrak{g}_{s,\pcC}}) = O(u^{\frac{1}{2}(-1-2s-\ch(z_j))}) =O(u^{-\frac{1+s}{2}}|dZ_j|_{\mathfrak{g}_{s,\pcC}}).
\]
Now we assume that the statement is true for all $l=0,1,\ldots,k$. Then $(\nabla^{\pcC})^{k+1}dZ_j$ can be written as the linear combination of
\[
\partial_{Z_{j_1}}\cdots\partial_{Z_{j_m}} \Gamma^j_{kl}((\nabla^{\pcC})^{k_1}dZ_{j_1})\otimes \cdots \otimes ((\nabla^{\pcC})^{k_m}dZ_{j_m})\otimes ((\nabla^{\pcC})^{k_{m+1}}dZ_k) \otimes ((\nabla^{\pcC})^{k_{m+2}}dZ_l)
\]
where $m+k_1+\cdots +k_{m+2}=k$. Then, by the induction hypothesis,
\begin{align}
&|\partial_{Z_{j_1}}\cdots\partial_{Z_{j_m}} \Gamma^j_{kl}((\nabla^{\pcC})^{k_1}dZ_{j_1})\otimes \cdots \otimes ((\nabla^{\pcC})^{k_m}dZ_{j_m})\otimes ((\nabla^{\pcC})^{k_{m+1}}dZ_k) \otimes ((\nabla^{\pcC})^{k_{m+2}}dZ_l)|\nonumber\\
&=O(u^{-(k_1+\cdots+k_m)\frac{1+s}{2}})|\partial_{Z_{j_1}}\cdots\partial_{Z_{j_m}} \Gamma^j_{kl}dZ_{j_1}\otimes \cdots \otimes dZ_{j_m}\otimes dZ_k \otimes dZ_l|_{\mathfrak{g}_{s,\pazocal{C}}}. \label{eq:inductionstep1}
\end{align}
 
From Proposition \ref{prop:formula_of_Gamma}, $\Gamma^j_{kl}$ can be represented as $P(\widehat{z})$ where $P\in x_1^{\ch(z_k)+\ch(z_l)-\ch(z_j)}u^{-2}w_2^{-2}\bC[z,\overline{z}]$, and thus
\[
\Gamma^j_{kl}\in \widehat{R}_{\ch(z_k)+\ch(z_l)-\ch(z_j)-1}
.\]
Hence, by Proposition \ref{differentiation1}, we have 
\[\partial_{Z_{j_1}}\cdots \partial_{Z_{j_m}}\Gamma^j_{kl}
\in 
\Gamma^j_{kl}\in \widehat{R}_{\ch(z_k)+\ch(z_l)-\ch(z_j)-1-m+\sum_{q=1}^m \ch(z_{j_q}))},\]
and by $|\widehat{x}_1|\asymp \widehat{u}^{1/2}$ (Proposition \ref{prop:shrinking}), $\widehat{w}_2\asymp \widehat{u}$, and $|\widehat{z}_j|=O(\widehat{u}^{1/2})$
\begin{align} \label{eq:m=0}
|\partial_{Z_{j_1}}\cdots \partial_{Z_{j_m}}\Gamma^j_{kl}|= O(u^{\frac{1}{2}(\ch(z_k)+\ch(z_l)-\ch(z_j)-1-m+\sum_{q=1}^m \ch(z_{j_q}))}).
\end{align}
Using Proposition \ref{prop:dZ_asymp}, we obtain
\begin{align*}
&|\partial_{Z_{j_1}}\cdots\partial_{Z_{j_m}} \Gamma^j_{kl}dZ_{j_1}\otimes \cdots \otimes dZ_{j_m}\otimes dZ_k \otimes dZ_l|_{\mathfrak{g}_{s,\pazocal{C}}} \\
&=O(u^{\frac{1}{2}(\ch(z_k)+\ch(z_l)-\ch(z_j)-1-m+\sum_{q=1}^m \ch(z_{j_q}))})|dZ_{j_1}|_{\mathfrak{g}_{s,\pazocal{C}}} \cdots |dZ_{j_m}|_{\mathfrak{g}_{s,\pazocal{C}}}|dZ_{k}|_{\mathfrak{g}_{s,\pazocal{C}}}|dZ_{l}|_{\mathfrak{g}_{s,\pazocal{C}}}\\
&=O(u^{\frac{1}{2}(\ch(z_k)+\ch(z_l)-\ch(z_j)-1-m+\sum_{q=1}^m \ch(z_{j_q}))})O(u^{\frac{1}{2}(-(m+2)s-\ch(z_k)-\ch(z_l)-\sum_{q=1}^m \ch(z_{j_q})))})\\
&=O(u^{\frac{1}{2}(-\ch(z_j)-1-m-(m+2)s)})
 \\
 &= O(u^{-\frac{(m+1)(1+s)}{2}}|dZ_j|_{\mathfrak{g}_{s,\pazocal{C}}}).\end{align*}
By \eqref{eq:inductionstep1},
\begin{align*}
|\partial_{Z_{j_1}}\cdots\partial_{Z_{j_m}} \Gamma^j_{kl}((\nabla^{\pcC})^{k_1}dZ_{j_1})\otimes \cdots \otimes ((\nabla^{\pcC})^{k_m}dZ_{j_m})\otimes ((\nabla^{\pcC})^{k_{m+1}}dZ_k) \otimes ((\nabla^{\pcC})^{k_{m+2}}dZ_l)|_{\mathfrak{g}_{s,\pcC}}\\
=O(u^{-\frac{1+s}{2}(m+1+k_1+\cdots+k_{m+2})}|dZ_j|_{\mathfrak{g}_{s,\pcC}})=O(u^{-\frac{(k+1)(1+s)}{2}}|dZ_j|_{\mathfrak{g}_{s,\pcC}}).
\end{align*}
Thus, by triangle inequality, we conclude that
\[
|(\nabla^{\pcC})^{k+1}dZ_j|_{\mathfrak{g}_{s,\pcC}}=O(u^{-\frac{(k+1)(1+s)}{2}}|dZ_j|_{\mathfrak{g}_{s,\pcC}}).
\]
\end{proof}

We write
\[
h_{k\overline{l}} := h(\partial_{Z_{k}},\partial_{\overline{Z}_{l}}).
\]
Then, on $V_1$,
\[
h = \sum_{k,l} h_{k\overline{l}} (dZ_{k} \odot d\overline{Z}_{l}).
\]
We consider the asymptotic estimate of the derivatives of $h_{k\overline{l}}$.

We define $\Lambda:\bR\times V_1 \to \bR$ as a function defined as
\begin{align}\label{eq:def_of_Lambda}
\mu_s(e^{\Lambda(t,Z)}\cdot Z) = -t
\end{align}
for $t\in \bR$ and $Z\in V_1$; $\Lambda(t,Z)$ is smooth by the analytic inverse function theorem. Then, $\Lambda(1,Z)=\lambda_1(Z)$, $\Lambda(0,Z)=\eta_1(Z)$. 
We write $\widetilde{z}_j(t,Z) = e^{\ch(z_j)\Lambda(t,Z)} Z_j$ for $j=2,\ldots,n+r.$ Similar to $\widehat{u}$, $\widehat{w}_k$, $\widehat{\pcV}_j$, and $\widehat{\Gamma}^j_{kl}$, we define
\begin{align*}
\widetilde{u}(t,Z) &= u(\widetilde{z}(t,Z)), &
\widetilde{w}_k(t,Z) &= w_k(\widetilde{z}(t,Z)) \quad(k=0,1,\ldots), \\
\widetilde{\pcV}_j(t,Z) &= \pcV_j(\widetilde{z}(t,Z)) \quad(j=1,\ldots n+r), &
\widetilde{\Gamma}^j_{kl}(t,Z) &= \!^{\pcU}\Gamma^j_{kl}(\widetilde{z}(t,Z))\quad(j,k,l=2,\ldots n+r).
\end{align*}
Then by \eqref{eq:def_of_Lambda},
\begin{align}\label{tildew}
(1+s)\widetilde{u}^s\widetilde{w}=t.
\end{align}
We define $\varsigma:\pcU \to \bR_{>0}$ as
\[
\varsigma(z) := {w}_2{u} +s w^2, \quad \widetilde{\varsigma}(Z) := \varsigma(\widetilde{z}).
\]
Recall that in \eqref{eq:def_of_G}, for each $j,k=1,\ldots,n+r$, $G_{j\overline{k}}:=G_{j\overline{k};1}$ is a function on $\pcU_{1}$ defined as
\[
G_{j\overline{k}}(z) := \widetilde{g}_s (x_1^{\ch(z_j)}\cH \partial_{z_j},\overline{x}_1^{\ch(z_k)}\cH \partial_{\overline{z}_k}).
\]
We define $\widetilde{G}_{j\overline{k}}(t,Z) := G_{j\overline{k}}(\widetilde{z}(t,Z))$ as a function of $\bR\times V_1$. Then, $h_{j\over{k}}$ can be written as follows:
\begin{align} \label{FTA}
h_{j\overline{k}} = \widetilde{G}_{j\overline{k}}(e^{\Lambda(1,Z)}\cdot Z) -  \widetilde{G}_{j\overline{k}}(e^{\Lambda(0,Z)}\cdot Z) = \int^1_0  \frac{\partial}{\partial t}\widetilde{G}_{j\overline{k}}(e^{\Lambda(t,Z)}\cdot Z)dt.
\end{align}

We extend Table \ref{tab:asymptotic-dict1} to Table \ref{tab:asymptotic-dict2} by adding functions and tensors used in the proof of Lemma \ref{cor:h_diff_asymp}. Note that the variable $t$ is restricted to the interval $[0,1]$, and the asymptotics are uniform in $t\in [0,1]$;
when $f(t,Z):[0,1]\times V_1\to \bR_{\geq 0}$ is a smooth function, we say $f(t,Z)=O(g)$ \emph{uniformly in $t\in [0,1]$} for some $g:V_1 \to \bR_{\geq 0}$ if there is a constant $c>0$ independent of $t$ such that
\[
f(t,Z) \leq cg(Z)
\]
whenever $\widehat{u}(Z)$ is large enough. Also, we say $f(t,Z) \asymp g(Z)$ \emph{uniformly in $t\in [0,1]$} if there is a constant $c_1,c_2>0$ independent of $t$ such that
\[
c_1g(Z) \leq f(t,Z) \leq c_2g(Z)
\]
whenever $\widehat{u}(Z)$ is large enough. We are not claiming that the asymptotics are uniformly true for $t\in \bR$, which might not be true.

\begin{table}[h]
  \centering
  \caption{Asymptotics of Functions and Tensors on $V_1$ (Extended from Table \ref{tab:asymptotic-dict1}).}
  \label{tab:asymptotic-dict2}

  \renewcommand{\arraystretch}{1.3}
  \setlength{\tabcolsep}{6pt}

  \begin{tabularx}{\linewidth}{@{}  l | l | X | r @{}}
    \toprule
      & \textbf{Definition} & \textbf{Asymptotics on $V_1$ Uniform in $t\in [0,1]$} & \textbf{Reference} \\
    \midrule       
    $e^{\Lambda(t,Z)-\Lambda(0,Z)}$  &  &    $e^{\Lambda(t,Z)-\Lambda(0,Z)} = 1+O(\widehat{u}^{-(1+s)/2})$. & \eqref{lambda-eta} \\
    $\widetilde{x}_1$  & $e^{\Lambda(t,Z)}$. & $|\widetilde{x}_1|\asymp \widehat{u}^{1/2}$. & Prop \ref{cor:tilde-hat}. \\
    $\widetilde{z}_j$  & $e^{\ch(z_j)\Lambda(t,Z)}Z_j$. & $|\widetilde{z}_j|=O(\widehat{u}^{1/2}).$ & Prop \ref{cor:tilde-hat}.  \\    
    $\widetilde{w}_{2j}$  & $\sum_k \ch(z_k)^{2j}|\widetilde{z}_k|^2$. & $\widetilde{w}_{2j}\asymp \widehat{u}$. & Prop \ref{cor:tilde-hat}.  \\
    $\widetilde{w}$  & $\sum_k \ch(z_k)|\widetilde{z}_k|^2$. & $\widetilde{w}=O(\widehat{u}^{-s})$ & \eqref{tildew}\\   
    $\widetilde{\varsigma}$  & $\widetilde{w}_2\widetilde{u}+s\widetilde{w}^2$. & $\widetilde{\varsigma}\asymp \widehat{u}^2$. & By def.  \\
    $\widetilde{\pcV}_j$  & $\widetilde{\overline{z}}_j(\ch(z_j)\widetilde{u}+s\widetilde{w})/\widetilde{\varsigma}$. & $|\widetilde{\pcV}_j|=O(\widehat{u}^{-\frac{1}{2}})$. & By def.\\
    $\partial_t\widetilde{G}_{jk}$  & $\partial_tG_{jk}(\widetilde{z})$. & 
    $|\partial_t\widetilde{G}_{jk}|=O(\widehat{u}^{\frac{1}{2}(\ch(z_j)+\ch(z_k)-2)})$. & Prop \ref{prop:G_diff_eq}.\\
    ${h}_{j\overline{k}}$  & $h(\partial_{Z_j},\partial_{\overline{Z}_k})$. & $|{h}_{j\overline{k}}|=O(\widehat{u}^{\frac{1}{2}(\ch(z_j)+\ch(z_k)-2)})$. & \eqref{FTA}\\
    $\partial_{Z_l}\widetilde{G}_{jk}$  & Partial derivative of $\widetilde{G}_{jk}$. & $|\partial_{Z_l}\widetilde{G}_{jk}| =O(\widehat{u}^{\frac{1}{2}(\ch(z_j)+\ch(z_k)+\ch(z_l)-3)})$. & Props  \ref{differentiation2}.\\
    $\partial_{Z_l}{h}_{j\overline{k}}$  & Partial derivative of ${h}_{j\overline{k}}$.& $|\partial_{Z_l}{h}_{j\overline{k}}|=O(\widehat{u}^{\frac{1}{2}(\ch(z_j)+\ch(z_k)+\ch(z_l)-3)})$. & 
    \eqref{FTA2}
    \\
    \bottomrule
  \end{tabularx}
\end{table}

\begin{proposition}\label{prop:G_diff_eq}
For each $j,k=1,\ldots,n+r$, we have
\[
\frac{\partial \widetilde{G}_{j\overline{k}}}{\partial t} = \frac{\partial \Lambda}{\partial t}\sum_{l=1}^{n+r} (H_{j\overline{k}}^l \widetilde G_{l\overline{k}} + H^{\overline{l}}_{j\overline{k}}\widetilde G_{j\overline{l}}),
\]
where
\begin{align*}
H^l_{j\overline{k}} :=& \delta^l_j \left(\ch(z_j)+s\frac{\widetilde w}{\widetilde u}\right)+\widetilde{x}_1^{\ch(z_j)-\ch(z_l)}\left(-s\frac{\widetilde{\overline{z}}_j\widetilde w}{\widetilde u^2}-\widetilde{\pcV}_j \ch(z_l)^2 \widetilde z_l+s\widetilde{\pcV}_j\frac{\widetilde w^2}{\widetilde u^2}\widetilde z_l\right), \quad H^{\overline{l}}_{j\overline{k}} := \overline{H^l_{k\overline{j}}}.
\end{align*}
In particular, we have $P\in x_1^{\ch(z_j)}\overline{x}_1^{\ch(z_k)}u^{-2}\varsigma^{-3}\bC[z,\overline{z}]_{\deg\leq 14}$ such that
\[
\frac{\partial \widetilde{G}_{j\overline{k}}}{\partial t}(Z) = P(\widetilde{z}).
\]
\end{proposition}
\begin{proof}
By the chain rule, we have
\[
\frac{\partial \widetilde{G}_{j\overline{k}}}{\partial t} = \sum_{l=1}^{n+r}\frac{\partial\widetilde{z}_l}{\partial t}\frac{\partial G_{j\overline{k}}}{\partial z_l}(\widetilde{z}) + \frac{\partial\widetilde{\overline z}_l}{\partial t}\frac{\partial G_{j\overline{k}}}{\partial \overline z_l}(\widetilde{z}).
\]
We have
\[
\frac{\partial\widetilde{z}_l}{\partial t} = \ch(z_l)\frac{\partial \Lambda}{\partial t}\widetilde{z}_l, \quad \frac{\partial\widetilde{\overline z}_l}{\partial t} = \ch(z_l)\frac{\partial \Lambda}{\partial t}\widetilde{\overline z}_l. 
\]
For $\partial_{z_l}G_{j\overline{k}}$, consider first when $l=1$:
\begin{align*}
x_1\partial_{x_1}G_{j\overline{k}} &= x_1\partial_{x_1} \widetilde{g}_s(x_1^{\ch(z_j)}\cH \partial_{z_j},\overline{x}_1^{\ch(z_k)}\cH\partial_{\overline z_k}) \\
=& \ch(z_j) G_{j\overline{k}} + x_1\widetilde{g}_s(x_1^{\ch(z_j)}\nabla_{\partial_{x_1}}^{\pcU}\cH\partial_{z_j},\overline{x}_1^{\ch(z_k)}\cH\partial_{\overline{z}_{k}})
\end{align*}
For any $l=1,\ldots,n+r$,
\begin{align*}
&\ch(z_l)z_l \widetilde{g}_s(x_1^{\ch(z_j)}\nabla_{\partial_{z_l}}^{\pcU}\cH\partial_{z_j},\overline{x}_1^{\ch(z_k)}\cH\partial_{\overline z_k}) \\
&=\ch(z_l)x_1^{\ch(z_j)}z_l\widetilde{g}_s(\nabla^{\pcU}_{\partial_{z_l}}\partial_{z_j}-\pcV_j \nabla^{\pcU}_{\partial_{z_l}}\mathbf{n},\overline{x}_1^{\ch(z_k)}\cH\partial_{\overline z_k})\\
&=\ch(z_l)x_1^{\ch(z_j)}z_l\widetilde{g}_s\left(-\pcV_j \ch(z_l)\partial_{z_l}+\sum_{m=1}^{n+r}\!^\pcU\Gamma^m_{lj}\partial_{z_m}-\pcV_j\sum_{m,q}\ch(z_q)z_q\!^\pcU \Gamma^m_{lq}\partial_m,\overline{x}_1^{\ch(z_k)}\cH\partial_{\overline z_k}\right).
\end{align*}
Using the fact that $\mathbf{n} = \sum\ch(z_j)z_j \partial_{z_j}$ and $g_s(\mathbf{n},\cH(-))=0$,
\begin{align*}
&\sum_{l=1}^{n+r} \ch(z_l)z_l \widetilde{g}_s(x_1^{\ch(z_j)}\nabla_{\partial_{z_l}}^{\pcU}\cH\partial_{z_j},\overline{x}_1^{\ch(z_k)}\cH\partial_{\overline z_k})\\
&=s\frac{w}{u}G_{j\overline{k}}+x_1^{\ch(z_j)}\sum_{l=1}^{n+r}\left(-s\frac{\overline{z}_jw}{u^2}-\pcV_j \ch(z_l)^2 z_l+s\pcV_j\frac{w^2}{u^2}z_l\right)x_1^{-\ch(z_l)}G_{l\overline{k}}.
\end{align*}
Hence, we have
\begin{align*}
\frac{\partial \widetilde{G}_{j\overline k}}{\partial t} =& \frac{\partial \Lambda}{\partial t}\sum_{l=1}^{n+r}\left(\ch(z_l) \widetilde z_l\partial_{z_l}G_{j\overline{k}}+\ch(z_l)\widetilde{ \overline z}_l\partial_{\overline z_l}G_{j\overline{k}}\right)\\
=& \frac{\partial \Lambda}{\partial t}\left(\ch(z_j)+\ch(z_k) +2s\frac{\widetilde{w}}{\widetilde{u}}\right)\widetilde G_{j\overline{k}} \\&+  \frac{\partial \Lambda}{\partial t}\widetilde{x}_1^{\ch(z_j)}\sum_{l=1}^{n+r}\left(-s\frac{\widetilde{\overline{z}}_j\widetilde w}{\widetilde u^2}-\widetilde{\pcV}_j \ch(z_l)^2 \widetilde z_l+s\widetilde \pcV_j\frac{\widetilde w^2}{\widetilde u^2}\widetilde z_l\right)\widetilde x_1^{-\ch(z_l)}\widetilde G_{l\overline{k}}\\
&+ \frac{\partial \Lambda}{\partial t}\widetilde x_1^{\ch(z_k)}\sum_{l=1}^{n+r}\left(-s\frac{\widetilde{z}_k\widetilde w}{\widetilde u^2}-\widetilde{\overline {\pcV}}_k \ch(z_l)^2 \widetilde{\overline{z}}_l+s\widetilde{\pcV}_j\frac{\widetilde{w}^2}{\widetilde{u}^2}\widetilde{\overline{z}}_l\right)\widetilde x_1^{-\ch(z_l)}\widetilde G_{j\overline{l}}.
\end{align*}

In summary, we have
\[
Q^l_{j\overline{k}}\in x_1^{\ch(z_j)-\ch(z_l)}u^{-2}\varsigma^{-1}\bC[z,\overline{z}]_{\deg\leq 8}
\]
such that
\[
H^l_{j\overline{k}} = Q^l_{j\overline{k}}(\widetilde{z}).
\]
We let $Q^{\overline l}_{j\overline{k}} = \overline{Q^l_{k\overline{j}}}$.
Also, by differentiating the defining equation
\[
\mu_s(e^{\Lambda(t,Z)}\cdot Z)=-t
\]
of $\Lambda(t,Z)$ by $t$, we have
\begin{align}
\frac{\partial \Lambda}{\partial t} &=\frac{1}{(1+s)\widetilde{u}^{s-1}\widetilde{\varsigma}}, \label{eq:Lambda_formula1}
\end{align}
and
\[
\widetilde{G}_{l\overline{k}} = \widetilde {x}_1^{\ch(z_l)} \widetilde {\overline x}_1^{\ch(z_k)} \left(\delta_{lk}(1+s) \widetilde{u}^s + s(1+s) \widetilde u^{s-1}\widetilde{z}_l \widetilde{\overline z}_k-  (1+s)\widetilde u^{s-1}\widetilde \varsigma\widetilde \pcV_l\widetilde{\overline{\pcV}}_k\right),
\]
so for all $l,k=1,\ldots,n+r$, we have
\[
R_{l\overline{k}}\in x_1^{\ch(z_l)}\overline x_1^{\ch(z_l)}\varsigma^{-2}\bC[z,\overline{z}]_{\deg\leq 6},
\]
such that
\[
\frac{\partial \Lambda}{\partial t}\widetilde{G}_{l\overline{k}} = R_{l\overline{k}}(\widetilde{z}).
\]
We let
\[
P = \sum_{l=1}^{n+r}\left( Q^l_{j\overline{k}} R_{l\overline{k}} +  Q^{\overline l}_{j\overline{k}} R_{j\overline{l}}\right) \in  x_1^{\ch(z_j)}\overline x_1^{\ch(z_k)}u^{-2}\varsigma^{-3}\bC[z,\overline{z}]_{\deg\leq 14},
\]
and we have
\[
\frac{\partial \widetilde{G}_{j\overline{k}}}{\partial t}(Z) = P(\widetilde{z}).
\]
\end{proof}

We have the following definition and proposition analogous to Proposition \ref{differentiation1}, which describes the partial differentiation of $P(\widetilde{z})$, where $P$ is a rational function in $\bC[z,\overline{z}][z_1^{-1},\overline{z}_1^{-1},u^{-1},w_2^{-1}]$, along $Z_j$

\begin{definition}\label{def:degree2}
We define $\widetilde{R}$ as a $\bC$-algebra of smooth $\bC$-valued functions on $[0,1]\times V_1$ such that
\[
\widetilde{R} := \{f:[0,1]\times V_1\to \bC: f(Z) = P(\widetilde{z}(t,Z)) \textrm{ for $P\in \bC[z,\overline{z}][x_1^{-1},\overline{x}_1^{-1},u^{-1},\varsigma^{-1}]$}\}.
\]
There is a filtration 
\[\cdots\subset \widetilde{R}_{-k} \subset \widetilde{R}_{-k+1}\subset\cdots \subset \widetilde{R}_{0} \subset \widetilde{R}_{1}\subset \cdots \subset \widetilde{R}
\]
where
\[
\widetilde{R}_{k}:=\{f\in \widetilde{R}:f(Z)=P(\widetilde{z}(Z)) \textrm{ for $P\in x_1^{\a}\overline{x}_1^\b u^\g \varsigma^\d \bC[z,\overline{z}]_{\deg \leq l}$ where $\a+\b+2\g+4\d+l\leq k$}\}.
\]
Then 
\[
\bigcup_{k\in \bZ} \widetilde{R}_k = \widetilde{R}, \qquad \bigcap_{k\in \bZ} \widetilde{R}_k =\{0\}.
\]
Then, we have an additive grading $\deg:\widetilde{R}\to \bZ \cup \{-\infty\}$ as 
\[
\deg f := \min\{k:f\in \widetilde{R}_k\} \textrm{\quad for $f\neq0$},
\]
and $\deg 0 := -\infty.$
\end{definition}

\begin{proposition} \label{differentiation2}
For each $j\in \{2,\overline{2},\ldots,n+r,\overline{n+r}\}$, $\widetilde{R}$ is closed under the operator $\partial_{Z_j}$, and $\partial{Z_j}$ is an operator of degree $\ch(z_j)-1$, i.e.,
\[
{\partial}_{Z_j}: \widetilde{R}_k \to \widetilde{R}_{k+\ch(z_j)-1}
\]
is well-defined.
\end{proposition}
\begin{proof}
By differentiating the defining equation
\[
\mu_s(e^{\Lambda(t,Z)}\cdot Z) = -t
\]
of $\Lambda(t,Z)$ by $Z_j$, we get
\begin{align}
\frac{\partial \Lambda}{\partial Z_j}&=-\frac{1}{2}\widetilde{x}_1^{\ch(z_j)}\widetilde{\overline{z}}_j\frac{\ch(z_j)\widetilde{u}+s\widetilde{w}}{\widetilde \varsigma} \quad(j=2,\overline{2},\ldots,n+r,\overline{n+r}). \label{eq:Lambda_formula2}
\end{align}
By the definition of $\widetilde{z}_k$, for $j=2,\overline{2},\ldots,n+r,\overline{n+r},\;k=1,\overline{1},\ldots,n+r,\overline{n+r}$
\begin{align*}
\partial_{Z_j} \widetilde{z}_k =& \widetilde{x}_1^{\ch(z_j)}\left(-\frac{1}{2}\ch(z_k)\widetilde{\overline z}_j\widetilde{z}_k\frac{\ch(z_j)\widetilde{u}+s\widetilde{w}}{\widetilde \varsigma}+\delta_{jk}\right).
\end{align*}
By the chain rule, for $P\in \bC[z,\overline{z}][z_1^{-1},\overline{z}_1^{-1},u^{-1},w_2^{-1}]$, we define
\begin{align*}
{\partial}_{Z_j} (P(\widetilde{z})) = \widetilde{x}_1^{\ch(z_j)}\sum_{k=1}^{n+r}&  \left(-\frac{1}{2}\ch(z_k)\widetilde{\overline{z}}_j{\widetilde{z}}_k\frac{\ch(z_j)\widetilde{u}+s\widetilde{w}}{ \widetilde\varsigma}+\delta_{jk}\right) (\partial_{z_k}P)(\widetilde{z}) \\
&+\left(-\frac{1}{2}\ch(z_k)\widetilde{\overline{z}}_j\widetilde{\overline{z}}_k\frac{\ch(z_j)\widetilde{u}+s\widetilde{w}}{ \widetilde\varsigma}\right) (\partial_{\overline z_k}P)(\widetilde{z}).
\end{align*}
Thus, ${\partial}_{Z_j}(P(\widetilde{z}))\in \widetilde{R}_{k+\ch(z_j)-1}$ if $P(\widetilde{z})\in \widetilde{R}_{k}$.
\end{proof}

In the following proposition, we compare the asymptotic estimate of $|\widetilde{z}_j|$ and $|\widehat{z}_j|$. 

\begin{proposition}\label{prop:limit}
\label{cor:tilde-hat}
If $f\in \widetilde{R}$, $f=O(\widehat{u}^{\deg{f}/2})$ uniformly in $t\in [0,1]$.
\end{proposition}
\begin{proof}
We first prove that as $\widehat{u}\to \infty$, $\lambda_1(Z)-\eta_1(Z)\to 0^+$. 
By \eqref{eq:Lambda_formula1}, $\Lambda(t,Z)$ is strictly increasing with respect to $t$, so $\lambda_1(Z)-\eta_1(Z)>0$. Also,
\[
\frac{\partial u(e^{\Lambda(t,Z)}\cdot Z)}{\partial t} = \frac{\partial \Lambda}{\partial t}2w(e^{\Lambda(t,Z)}\cdot Z) = \frac{\partial \Lambda}{\partial t} \frac{2t}{(1+s)u^s(e^{\Lambda(t,Z)}\cdot Z)},
\]
so when $t> 0$, $\frac{\partial u(e^{\Lambda(t,Z)}\cdot Z)}{\partial t}>0$. Hence, $u(e^{\lambda_1(Z)}\cdot Z) > u(e^{\eta_1(Z)}\cdot Z)$.
By definition of $\lambda_1(Z)$,
\begin{align*}
1&= u^s(e^{\lambda_1(Z)}\cdot Z) w(e^{\lambda_1(Z)}\cdot Z)
\\ &\geq u^s(e^{\eta_1(Z)}\cdot Z) \sum_{j=1}^{n+r} e^{2\ch(z_j)(\lambda_1(Z)-\eta_1(Z))}\ch(z_j)e^{2\ch(z_j)\eta_1(Z)}|Z_j|^2\\
&\geq \widehat u^s(Z) \;\left(e^{2(\lambda_1(Z)-\eta_1(Z))}\sum_{j=1}^{n} |\widehat{x}_j|^2 -\sum_{k=1}^r d_k|\widehat{p}_k|^2\right).
\end{align*}
Since $w(\widehat{z})=0$,
\[
\sum_{j=1}^{n} |\widehat{x}_j|^2 =\sum_{k=1}^r d_k|\widehat{p}_k|^2,
\]
and
\begin{align*}
1&\geq \widehat u^s(Z) \left(e^{2(\lambda_1(Z)-\eta_1(Z))}-1\right) \sum_{j=1}^{n} |\widehat{x}_j|^2\\
&\geq \widehat u^s(Z) \left(e^{2(\lambda_1(Z)-\eta_1(Z))}-1\right) \frac{1}{2}\sum_{j=1}^{n} 2|\widehat{x}_j|^2 \\
&= \widehat u^s(Z) \left(e^{2(\lambda_1(Z)-\eta_1(Z))}-1\right) \frac{1}{2}\left(\sum_{j=1}^{n} |\widehat{x}_j|^2+\sum_{k=1}^r d_k|\widehat{p}_k|^2\right)\\
&\geq \frac{1}{2} \widehat u^{s+1}(Z) \left(e^{2(\lambda_1(Z)-\eta_1(Z))}-1\right).
\end{align*}
Thus,
\begin{align} \label{lambda-eta}
1\leq e^{2(\lambda_1(Z)-\eta_1(Z))} \leq 1+2\widehat{u}^{-1-s}.
\end{align}
Hence, as $\widehat{u}\to \infty$, $e^{2(\lambda_1(Z)-\eta_1(Z))}\to 1^+$ and thus $\lambda_1(Z)-\eta_1(Z)\to 0^+$.

By \eqref{eq:Lambda_formula1}, $\Lambda(t,Z)$ is strictly increasing with respect to $t$, so $\lambda_1(Z)-\eta_1(Z)>\Lambda(t,Z)-\eta_1(Z)>0$ for all $t\in [0,1]$. Hence,
\[
e^{-2|\ch(z_j)|\,(\lambda_1(Z)-\eta_1(Z))}|\widehat{z}_j(Z)| \leq |\widetilde{z}_j(t,Z)| \leq e^{2|\ch(z_j)|\,(\lambda_1(Z)-\eta_1(Z))}|\widehat{z}_j(Z)|.
\]
We may assume that $\widehat{u}$ is large enough so the $e^{2(\lambda_1(Z)-\eta_1(Z))}\leq 2$. Then,
\[
\frac{1}{2^{|\ch(z_j)|}}|\widehat{z}_j(Z)| \leq |\widetilde{z}_j(t,Z)| \leq 2^{|\ch(z_j)|}|\widehat{z}_j(Z)|.
\]
Hence,
\[
\frac{1}{2^{2\max_j|\ch(z_j)|}}\widehat{u}(Z) \leq \widetilde{u}(t,Z) \leq 2^{2\max_j|\ch(z_j)|}\widehat{u}(Z).
\]
For an even integer $k=2l$, we have
\[
\widetilde{u}(t,Z)\leq  \widetilde{w}_{2l}(t,Z) \leq\max_{1\leq m\leq n+r}  (\ch(z_m))^{2l} \widetilde{u}(t,Z),
\]
so we conclude that
\[
\frac{1}{2^{2\max_j|\ch(z_j)|}}\widehat{u}(Z) \leq \widetilde{w}_{2l}(t,Z)  \leq 2^{2\max_j|\ch(z_j)|}\max_{1\leq m\leq n+r} (\ch(z_m))^{2l}\widehat{u}(Z).
\]
From Proposition \ref{cor3.25}, we have $|\widetilde{x}_1|\asymp |\widehat{x}_1|\asymp \sqrt{\widehat u}$.

Finally, if $P\in \bC[z,\overline{z}]_{\deg \leq k}$, then since $\widetilde{u}\asymp \widehat{u}$ and $\widetilde\varsigma\asymp \widehat{u}^2$, so when $\a,\b,\g\in \bZ$,
\[
|\widetilde{x}_1^\a \widetilde{u}^\b\widetilde{\varsigma}^\g P(\widetilde{z})| = O(\widehat{u}^{\frac{1}{2}(\a+2\b+4\g+k)}).
\]
\end{proof}

\begin{lemma} \label{cor:h_diff_asymp}
For $m\geq 0$, $j_1,\ldots,j_m\in \{2,\overline{2},\ldots,n+r,\overline{n+r}\}$,  $j,k\in \{2,\ldots,n+r\}$,
\[
|\partial_{Z_{j_1}}\cdots\partial_{Z_{j_m}} h_{j\overline{k}}dZ_{j_1}\otimes \cdots \otimes dZ_{j_m}\otimes (dZ_j \odot d\overline{Z}_k)| = O(\widehat{u}^{\frac{-(2+m)(1+s)}{2}})
\]
\end{lemma}
\begin{proof}
By Propositions \ref{prop:G_diff_eq}, we have $P\in x_1^{\ch(z_j)}\overline{x}_1^{\ch(z_k)}u^{-2}\varsigma^{-3}\bC[z,\overline{z}]_{\deg\leq 14}$ such that
\[
\frac{\partial \widetilde{G}_{j\overline{k}}}{\partial t}(Z) = P(\widetilde{z}),
\]
and by Propositions \ref{differentiation2}
\[
\frac{\partial^{m+1} \widetilde{G}_{j\overline{k}}}{\partial t\, \partial Z_{j_1}\cdots \partial Z_{j_m}} \in \widetilde{R}_{\ch(z_j)+\ch(z_k)-2-m\sum_{l=1}^{m}\ch(z_{j_m})}.
\]
By Proposition \ref{prop:limit}, 
\[
|\widetilde{\partial}_{Z_{j_1}}\cdots\widetilde{\partial}_{Z_{j_m}} P(\widetilde{z})| = O(\widehat{u}^{\frac{1}{2}(\ch(z_j)+\ch(z_k)-2-m+\sum_{l=1}^{m}\ch(z_{j_m}))})
\]
uniformly in $t\in [0,1]$.

Using the fundamental theorem of calculus, we have
\begin{align} \label{FTA2}
\frac{\partial^{m+1} h_{j\overline{k}}}{\partial{Z_{j_1}}\cdots\partial{Z_{j_m}} } = \int^1_0\frac{\partial^{m+1} \widetilde{G}_{j\overline{k}}}{\partial t\, \partial Z_{j_1}\cdots \partial Z_{j_m}}dt.
\end{align}
Here, the partial differentiation commutes with the integration by the dominated convergence theorem, since the integrand is bounded above uniformly in $t\in [0,1]$. Thus, we obtain
\[
|\partial_{Z_{j_1}}\cdots\partial_{Z_{j_m}} h_{j\overline{k}}| \leq  \int^1_0\left|\frac{\partial^{m+1} \widetilde{G}_{j\overline{k}}}{\partial t\, \partial Z_{j_1}\cdots \partial Z_{j_m}}\right| dt = O(\widehat{u}^{\frac{1}{2}(\ch(z_j)+\ch(z_k)+\sum_{l=1}^{m}\ch(z_{j_m})-2-m)}).
\]
Using Proposition \ref{prop:dZ_asymp}, we obtain
\begin{align*}
&|\partial_{Z_{j_1}}\cdots\partial_{Z_{j_m}} h_{j\overline{k}}dZ_{j_1}\otimes \cdots \otimes dZ_{j_m}\otimes (dZ_j \odot d\overline{Z}_k)| \\
&=O(\widehat{u}^{\frac{1}{2}(\ch(z_j)+\ch(z_k)+\sum_{l=1}^{m}\ch(z_{j_m})-2-m)})O(\widehat{u}^{\frac{1}{2}(-\ch(z_j)-\ch(z_k)-\sum_{l=1}^{m}\ch(z_{j_m})-s(2+m))})\\
&=O(\widehat{u}^{\frac{-(2+m)(1+s)}{2}})
\end{align*}
\end{proof}

}

\subsection{Bounded Calabi-Yau geometry of $X_{CY}$}
\label{section3.2}
We move on to the problem of the existence of a bounded Calabi-Yau form. We assume the Calabi-Yau condition $d_1+\cdots+d_r=n$. We have an explicit formula for the non-vanishing holomorphic volume form in Definition \ref{volumeform}. We need to prove that there is a metric $\mathfrak{g}_s$ that makes $\nabla^k \Omega$ and $\nabla^k \Theta$ bounded for all $k$. It turns out that the size of the norm $\Omega$ and $\Theta$ is controlled by the constant $s\geq 0$ which was introduced in 
\eqref{Kahler potential} to define the K\"ahler metric $\mathfrak{g}_s$ on $\cX$. Recall that the metric $\widetilde{g}_s$ is defined by the K\"ahler potential $u^{1+s}.$


From Proposition~\ref{prop:tensor_compariton}, when $T\in \Gamma(T^{\otimes p}_\bC M \otimes (T_\bC^*)^{\otimes q} M)$ is a tensor field on a asymptotically conical manifold $(M,g)$ with a diffeomorphism $\Phi:C/K\to M/K'$ onto Riemannian cone $(C,g_C)$ as defined in Definition \ref{ACmetric}, $|T|_{g}=|\Phi^*T|_{\Phi^*g}=O(\rho_C^{-\n})$ if and only if $ |\Phi^* f|_{g_C} = O(\rho_C^{-\n})$. 

The norms of ${\Omega}$, $\Theta$, $\nabla^k{\Omega}$, and $\nabla^k\Theta$ are defined by the extension of $\mathfrak{g}_s$ to the tensor products. We first prove that $|\Theta|_{\mathfrak{g}_s}$ and $|\Omega|_{\mathfrak{g}_s}$ are bounded when $s=\frac{1}{n+r-1}$.

\begin{proposition}
The norm $|\Omega|_{\mathfrak{g}_s}$ is bounded if and only if $s\geq \frac{1}{n+r-1}$, and the norm $|\Theta|_{\mathfrak{g}_s}$ is bounded if and only if $0\leq s\leq \frac{1}{n+r-1}$.
\label{expcompvolumeform}
\end{proposition}
\begin{proof}
We let $\{V_1,\ldots,V_n\}$ be an open cover of $\pcC$ defined in Lemma \ref{prop:shrinking}. We will show that the statement is true for each open subset $V_j$. We will only prove this for $V_1$, since the proof is generally the same for other open subsets. Recall that on $V_1$, we had the holomorphic chart $(Z_2,\ldots,Z_{n+r})$ and a function $\eta_1:V_1 \to \bR$ satisfying
\[
\mu_s(e^{\eta_1(Z)}\cdot Z) = 0.
\]
Then in the proof of Theorem \ref{keylemma}, we defined
\[
\widehat{z}_j := e^{\ch(z_j)\eta_1(Z)}Z_j, \quad \widehat{u}:= u(\widehat{z}(Z)), \quad \widehat{w}_{k} = w_k(\widehat{z}(Z)).
\]
Then,
\[
\Phi^* \Omega = - dZ_2 \wedge \cdots \wedge dZ_{n+r}, \quad \Phi^*\Theta = - \partial_{Z_2}\wedge \cdots \wedge \partial_{Z_{n+r}}.
\]
From Proposition~\ref{prop:tensor_compariton}, it is enough to prove that $|\Phi^*\Omega|_{\mathfrak{g}_{s,\pcC}}$ is bounded if and only if $s\geq \frac{1}{n+r-1}$, and the norm $|\Phi^*\Theta|_{\mathfrak{g}_{s,\pcC}}$ is bounded if and only if $0\leq s\leq \frac{1}{n+r-1}$.

We let $(\mathfrak{g}_{s,\pcC}(\partial_{Z_j},\partial_{\overline{Z}_{k}}))^{2\leq j,k \leq n+r}$ be the matrix representing the Riemannian metric $\mathfrak{g}_{s,\pcC}$ on the differential forms $dZ_2,\ldots,dZ_{n+r}$. Then,
\begin{align*}
|\Phi^* \Omega|_{\mathfrak{g}_{s,\pcC}} =& |\det((\mathfrak{g}_{s,\pcC}(\partial_{Z_j},\partial_{\overline{Z}_{k}})^{2\leq j,k \leq n+r})|, \\
|\Phi^* \Theta|_{\mathfrak{g}_{s,\pcC}} =& |\det((\mathfrak{g}_{s,\pcC}(\partial_{Z_j},\partial_{\overline{Z}_{k}})_{2\leq j,k \leq n+r})|.
\end{align*}
From Hadamard's inequality for the determinant, we have
\[
|\det((\mathfrak{g}_{s,\pcC}(\partial_{Z_j},\partial_{\overline{Z}_{k}})^{2\leq j,k \leq n+r})| \leq |dZ_2|_{\mathfrak{g}_{s,\pcC}}\cdots |dZ_{n+r}|_{\mathfrak{g}_{s,\pcC}}.
\]
From Proposition \ref{prop:dZ_asymp}, we had $|dZ_j|_{\mathfrak{g}_{s,\pcC}}\asymp \widehat{u}^{-\frac{1}{2}(\ch(z_j)+s)}$, so
\[
|\det((\mathfrak{g}_{s,\pcC}(\partial_{Z_j},\partial_{\overline{Z}_{k}})^{2\leq j,k \leq n+r})| = O(\widehat{u}^{-\frac{1}{2}(\ch(z_2)+\cdots + \ch(z_{n+r})+(n+r-1)s)}).
\]
From the Calabi-Yau condition, $\ch(z_2)+\cdots + \ch(z_{n+r}) = -1$, so
\[
|\Phi^*\Omega|_{\mathfrak{g}_{s,\pcC}} = O(\widehat{u}^{-\frac{1}{2}((n+r-1)s-1)}).
\]
Hence, $|\Phi^*\Omega|_{\mathfrak{g}_{s,\pcC}}$ is bounded in and only if $(n+r-1)s-1\geq0$. 

To give a bound on $|\Phi^*\Theta|_{\mathfrak{g}_{s,\pcC}}$, we  prove that $|\partial_{Z_j}|_{\mathfrak{g}_{s,\pcC}}\asymp \widehat{u}^{\frac{1}{2}(\ch(z_j)+s)}$. For the lower bound, it is followed from Proposition \ref{prop:dZ_asymp} since
\[
1 = |dZ_j (\partial_{Z_j})|  \leq |dZ_j|_{\mathfrak{g}_{s,\pcC}}|\partial_{Z_j}|_{\mathfrak{g}_{s,\pcC}}.
\]
For the upper bound, we use \eqref{eq:compute_g1} so that
\[
|\partial_{Z_{j}}|_{\mathfrak{g}_{s,\pcC}}=\sqrt{G_{j\overline{j}}(\widehat{z})},
\]
and
\begin{align*}
|G_{j\overline{j}}| &= |\widetilde{g}_s(x_1^{\ch(z_j)}\cH\partial_{z_j},\overline{x}_1^{\ch(z_j)}\cH\partial_{\overline{z}_j})|\\
&\leq |x_1|^{2\ch(z_j)}\widetilde{g}_s(\partial_{z_j},\partial_{\overline{z}_j}) \leq (1+s)^2|x_1|^{2\ch(z_j)}u^s.
\end{align*}
By Proposition \ref{prop:shrinking}, $|\widehat{x}_1| \asymp \widehat{u}^{1/2}$ on $V_1$, so we have
\[
|\partial_{Z_j}|_{\mathfrak{g}_{s,\pcC}} = O(\widehat{u}^{\frac{1}{2}(\ch(z_j)+s)}).
\]
By Hadamard's inequality again,
\[
|\Phi^*\Theta|_{\mathfrak{g}_{s,\pcC}} \leq |\partial_{Z_2}|_{\mathfrak{g}_{s,\pcC}}\cdots|\partial_{Z_{n+r}}|_{\mathfrak{g}_{s,\pcC}} = O(\widehat{u}^{\frac{1}{2}(\ch(z_2)+\cdots + \ch(z_{n+r})+(n+r-1)s)}).
\]
By the Calabi-Yau condition, $\ch(z_2)+\cdots + \ch(z_{n+r})=-1$ and
\[
|\Phi^*\Theta|_{\mathfrak{g}_{s,\pcC}} = O(\widehat{u}^{\frac{1}{2}((n+r-1)s-1)}).
\]
Thus, $|\Phi^*\Theta|_{\mathfrak{g}_{s,\pcC}}$ is bounded if $(n+r-1)s\leq 1$.
\end{proof}

\begin{proposition}
We have $|\nabla^k \Theta|_{\mathfrak{g}_s}/|\Theta|_{\mathfrak{g}_s}=O(\rho_{\pcC}^{-k})$ and $|\nabla^k \Omega|_{\mathfrak{g}_s}/|\Omega|_{\mathfrak{g}_s}=O(\rho_{\pcC}^{-k})$.
\label{4.16}
\end{proposition}
\begin{proof}
Using Proposition~\ref{prop:tensor_compariton} and \eqref{eq:radial_coordinate}, it is enough to prove 
\[|(\nabla^{\pazocal{C}})^k \Phi^*\Theta|_{\mathfrak{g}_{s,\pazocal{C}}}/|\Phi^*\Theta|_{\mathfrak{g}_{s,\pazocal{C}}}=O\left((u|_{\mu^{-1}_s(0)})^{-k\frac{1+s}{2}}\right),\quad |(\nabla^{\pazocal{C}})^k \Phi^*\Omega|_{\mathfrak{g}_{s,\pazocal{C}}}/|\Phi^*\Omega|_{\mathfrak{g}_{s,\pazocal{C}}}=O\left((u|_{\mu^{-1}_s(0)})^{-k\frac{1+s}{2}}\right).\]

To prove that there exists $C_k$ such that $|(\nabla^{\pazocal{C}})^k \Phi^*\Theta|_{\mathfrak{g}_{s,\pazocal{C}}}< C_k(u|_{\mu^{-1}_s(0)})^{-k\frac{1+s}{2}}|\Phi^*\Theta|_{\mathfrak{g}_{s,\pazocal{C}}}$, we use induction on $k$. We first observe that in $V_1$,
$$
(\nabla^{\pazocal{C}}) \Phi^*\Theta =  \left(\sum_{2\leq j,k \leq n+r} \Gamma^j_{j,k} dZ_k\right) \otimes \Phi^*\Theta,
$$
where $\Gamma^j_{k,l}$ is the Levi-Civita connection of $\nabla^{\pazocal{C}}$.
We define $\theta:= \sum_{2\leq j,k \leq n+r} \Gamma^j_{j,k} dZ_k$. If we prove that there is $D_l>0$ such that $|(\nabla^{\pazocal{C}})^l\theta|_{\mathfrak{g}_{s,\pazocal{C}}}<D_{l}\widehat{u}^{-(l+1)\frac{1+s}{2}}$ for all $l\geq 0$, we are done since
\begin{align*}
|(\nabla^{\pazocal{C}})^k \Phi^*\Theta|_{\mathfrak{g}_{s,\pazocal{C}}}&=\left|(\nabla^{\pazocal{C}})^{k-1} (\theta \otimes \Phi^*\Theta)\right|_{\mathfrak{g}_{s,\pazocal{C}}} \\
&= \left|\sum_{l=1}^{k-1} {k\choose l} (\nabla^{\pazocal{C}})^l\theta \otimes (\nabla^{\pazocal{C}})^{k-l-1} \Phi^*\Theta\right|_{\mathfrak{g}_{s,\pazocal{C}}}\\
&\leq\sum_{l=1}^{k-1} {k\choose l} \left|(\nabla^{\pazocal{C}})^l\theta\right|_{\mathfrak{g}_{s,\pazocal{C}}} \left|(\nabla^{\pazocal{C}})^{k-l-1} \Phi^*\Theta\right|_{\mathfrak{g}_{s,\pazocal{C}}} \\
&\leq  \sum_{l=1}^{k-1} C_{k-l-1}D_l{k\choose l} (u|_{\mu^{-1}_s(0)})^{-k\frac{1+s}{2}} \left| \Phi^*\Theta\right|_{\mathfrak{g}_{s,\pazocal{C}}}.
\end{align*}
The proof of $|(\nabla^{\pazocal{C}})^k \Phi^*\Omega|_{\mathfrak{g}_{s,\pazocal{C}}}=O( (u|_{\mu^{-1}_s(0)})^{-k\frac{1+s}{2}})|\Phi^*\Omega|_{\mathfrak{g}_{s,\pazocal{C}}}$ is similar since
$$
(\nabla^{\pazocal{C}}) \Phi^*\Omega = -\left(\sum_{2\leq j,k \leq n+r} \Gamma^j_{j,k}dZ_k\right) \otimes \Phi^*\Omega=-\theta \otimes \Phi^*\Omega
$$
on $V_1$.

Now we prove that for each $k>0$, we have 
$$
|(\nabla^{\pazocal{C}})^k \theta|_{\mathfrak{g}_{s,\pazocal{C}}} =O(\rho_0^{-(k+1)\frac{1+s}{2}}(\ud{u})).
$$
We will only prove that this is true on $V_1$, since the proof is generally the same for other open subsets. $(\nabla^{\pazocal{C}})^k \theta$ is a linear combination of the form
\[
(\partial_{Z_{j_1}}\cdots \partial_{Z_{j_q}}\Gamma^l_{lm})(\nabla^\pcC)^{k_1}dZ_{j_1}\otimes \cdots (\nabla^\pcC)^{k_q}dZ_{j_q}\otimes (\nabla^\pcC)^{k_{q+1}}dZ_{m}
\]
where $q+k_1+\cdots + k_{q+1}=k$. From Lemma \ref{cor:cov_asymp}, we had
\[
|(\nabla^{\pcC})^k dZ_j|_{\mathfrak{g}_{s,\pcC}} = O(u^{-k\frac{1+s}{2}} |dZ_j|_{\mathfrak{g}_{s,\pcC}}),
\]
and by Proposition \ref{prop:dZ_asymp},
\[
|(\nabla^{\pcC})^k dZ_j|_{\mathfrak{g}_{s,\pcC}} = O(u^{-\frac{1}{2}(k(1+s)+\ch(z_j)+s)}),
\]
From Proposition \ref{prop:formula_of_Gamma}, we had
\[
\Gamma^l_{lm} \in \widehat{R}_{\ch(z_m)-1},
\]
and by Proposition \ref{differentiation1}, we have
\[
\partial_{Z_{j_1}}\cdots \partial_{Z_{j_q}}\Gamma^l_{lm} \in \widehat{R}_{\ch(z_m)-1-q+\sum_{s=1}^q\ch(z_{j_s})},
\]
and thus
\[
|\partial_{Z_{j_1}}\cdots \partial_{Z_{j_q}}\Gamma^l_{lm} | = O(\widehat{u}^{\frac{1}{2}(\ch(z_m)-1-q+\sum_{s=1}^q\ch(z_{j_s}))}).
\]
Putting all these equations together, we obtain
\begin{align*}
&|(\partial_{Z_{j_1}}\cdots \partial_{Z_{j_q}}\Gamma^l_{lm})(\nabla^\pcC)^{k_1}dZ_{j_1}\otimes \cdots (\nabla^\pcC)^{k_q}dZ_{j_q}\otimes (\nabla^\pcC)^{k_{q+1}}dZ_{m}|_{\mathfrak{g}_{s,\pcC}}\\
&=O(\widehat{u}^{-\frac{1}{2}(1+q+(k_1+\cdots+k_{q+1})(1+s)+(q+1)s)})=O(\widehat{u}^{-\frac{1}{2}(1+s)(1+k)}).
\end{align*}
\end{proof}

\begin{remark}
The Riemannian manifold $(\pazocal{C},\mathfrak{g}_{s,\pazocal{C}})$ is geodesically incomplete, and its completion $\overline{\pazocal{C}}$, $$\overline{\pazocal{C}}\cong F^{-1}(0)/U(1)$$ has a conical singularity at $0\in (\mu^{-1}_s(0)\cup\{0\})/U(1)$. The diffeomorphism $\Phi^{-1}:\cX\setminus \bC\bP^{n-1}\to \pazocal{C}$ extends to $\cX \to \overline{\pazocal{C}}$, which is a resolution of the conical singularity. 
\end{remark}

\begin{remark}
Using an argument similar to Proposition \ref{biholomorphism}, there is a $U(1)$ action on $\mu_s^{-1}(-1)$, and we have an isomorphism $\mu_s^{-1}(-1)/U(1)\cong X_{LG}$ between orbifolds which makes the following diagram commutes;
\[\begin{tikzcd}
	{\mu_s^{-1}(-1)} & {\bC^n \times (\bC^{r}\setminus \{0\})} \\
	{\mu_s^{-1}(-1)/U(1)} & {X_{LG}}
	\arrow[hook, from=1-1, to=1-2]
	\arrow[two heads, from=1-1, to=2-1]
	\arrow[two heads, from=1-2, to=2-2]
	\arrow["\cong", from=2-1, to=2-2]
\end{tikzcd}.\]
We have the following commutative diagram that relates $X_{CY}$ and $X_{LG}$:
\[\begin{tikzcd}
	&& {\bC^{n+r}\setminus\{0\}} \\
	{\bC^n \times (\bC^{r}\setminus \{0\})} & {\mu_s^{-1}(-1)} & {\mu_s^{-1}(0)\cup\{0\}} & {\mu_s^{-1}(1)} & {(\bC^n\setminus\{0\}) \times \bC^{r}} \\
	& {X_{LG}} & {(\mu_s^{-1}(0)\cup\{0\})/U(1)} & {X_{CY}}
	\arrow[hook, from=2-1, to=1-3]
	\arrow[two heads, from=2-1, to=3-2]
	\arrow[hook, from=2-2, to=2-1]
	\arrow[from=2-2, to=2-3]
	\arrow[two heads, from=2-2, to=3-2]
	\arrow[two heads, from=2-3, to=3-3]
	\arrow[from=2-4, to=2-3]
	\arrow[hook, from=2-4, to=2-5]
	\arrow[two heads, from=2-4, to=3-4]
	\arrow[hook, from=2-5, to=1-3]
	\arrow[two heads, from=2-5, to=3-4]
	\arrow[from=3-2, to=3-3]
	\arrow[from=3-4, to=3-3]
\end{tikzcd}.\]
Like $X_{LG}$ and $X_{CY}$, ${(\mu_s^{-1}(0)\cup\{0\})/U(1)}$ is also a toric variety. An interesting question is whether one can generalize the work of Li-Wen \cite{LW} to the orbifold $X_{LG}$ and develop a parallel geometric theory on $X_{LG}$ corresponding to a bounded Calabi-Yau geometry on $X_{CY}$ using the isomorphism $\mu_s^{-1}(-1)/U(1)\cong X_{LG}$.
\end{remark}

\begin{remark}
A metric on $\cX$ that makes $(\cX,J,\Omega)$ a K\"ahler manifold of bounded Calabi-Yau geometry is not unique. By a similar argument, when $\pcU$ is given a metric defined by the K\"ahler potential
$$
(|x_1|^2+\cdots+|x_n|^2)^s+|p_1|^2+\cdots+|p_r|^2,
$$
$\cX$ inherits a K\"ahler metric of bounded geometry, and it is a bounded Calabi-Yau geometry when $s=\frac{1}{n-1}$. One of the interesting problems is to analyze the moduli space of K\"ahler metrics of bounded Calabi-Yau geometry inside the K\"ahler class of $\mathfrak{g}_s$. Another interesting question is whether there is a Ricci-flat metric inside this moduli space.
\end{remark}

\subsection{Strong Ellipticity of $W$}
\label{section3.3}

Now we prove that $W$ is strongly elliptic if $d_{\max} \leq 2d_{\min}$ and $s<1$. (See Definition \ref{Definition:strell} for the definition of strong ellipticity.) The proof relies on the fact that $W_1,\ldots,W_r$ gives a smooth complete intersection in ${\bC}\bP^{n-1}$.

In Lemma \ref{se} below, we directly compute the order of growth rates of $|\nabla W|$ and $|\nabla^k W|$ with respect to $u$. (See \eqref{eq:radial_coordinate} for the relation between $\rho_\pcC$ and $u$.) According to Definition \ref{Definition:strell}, we need to prove inequality \eqref{Definition:strelleq} for $W$. For a point $[z_1,\ldots,z_{n+r}]\in \cX$, if $[z_1,\ldots,z_n]\in \bP^{n-1}$ is an element of $V(\ud{W})$, the inequality \eqref{Definition:strelleq} will be demonstrated using the fact that $W_1,\ldots,W_r$ defines a smooth complete intersection. When $[z_1,\ldots,z_n]$ is not in $V(\ud{W})$, at least one of $|W_1(\ud{x})|,\ldots,|W_r(\ud{x})|$ is nonzero, so again we have the inequality \eqref{Definition:strelleq}.

\begin{lemma}\label{se}
If $d_{\max}\leq 2d_{\min}$ and $s<1$, the holomorphic function $W$ is strongly elliptic on $(\cX,\mathfrak{g}_{s})$.
\label{Wstrell}
\end{lemma}
\begin{proof} 
We want to show that for some constant $C$, 
\begin{align*}
|\nabla W|_{\mathfrak{g}_s}^k &\geq C^k{(u|_{\mu^{-1}_s(-1)})}^{\frac{1}{2}k(d_{\min}-s)},\\
|\nabla^k W|_{\mathfrak{g}_{s}}&= O\left((u|_{\mu^{-1}_s(-1)})^{\frac{1}{2}({d_{\max}+1}-k{(1+s)})}\right).
\end{align*} 
If those asymptotics are true, since $z\to \infty$ by the definition of strong ellipticity, we may assume that $u|_{\mu^{-1}_s(-1)}>1$. Also, $k\geq 2$,
so
\begin{align*}
|\nabla W|_{\mathfrak{g}_s}^k &\geq C^k{(u|_{\mu^{-1}_s(-1)})}^{\frac{1}{2}k(d_{\min}-s)}\geq C^k{(u|_{\mu^{-1}_s(-1)})}^{d_{\min}-s}.
\end{align*} 
Also, by the definition of the big $O$, we have a constant $D$ such that
\[
|\nabla^k W|_{\mathfrak{g}_{s}}\leq D (u|_{\mu^{-1}_s(-1)})^{\frac{1}{2}({d_{\max}+1}-k{(1+s))}}\leq D (u|_{\mu^{-1}_s(-1)})^{\frac{1}{2}d_{\max}-\frac{1}{2}-s}.
\]
Now we have $d_{\min}\geq d_{\max}/2 > (d_{\max}-1)/2$, so
\[
|\nabla W|_{\mathfrak{g}_s}^k - |\nabla^k W|_{\mathfrak{g}_{s}} \geq c {(u|_{\mu^{-1}_s(-1)})}^{d_{\min}-s}
\]
for some constant $c$ whenever $u|_{\mu^{-1}_s(-1)}$ is large enough. Thus, as $z\to \infty$,
\[
|\nabla W|_{\mathfrak{g}_s}^k - |\nabla^k W|_{\mathfrak{g}_{s}}  \to \infty.
\]

We give the lower bound on $|\nabla W|_{\mathfrak{g}_s}^2$. Note that ${\nabla}^{\pcU} W$ is already horizontal, so we have $\pi^*\nabla W = {\nabla}^{\pcU}W$. Then $|\nabla W|_{\mathfrak{g}_s}^2$ is computed as follows:
\begin{align*}
&|\nabla^{\pcU} W|_{\widetilde{g}_s}^2 = \left|\sum_{j=1}^n \frac{\partial W}{\partial x_j}dx_j+\sum_{k=1}^rW_kdp_k\right|_{\widetilde{g}_s}^2\\
&=\frac{1}{1+s}{u}^{-s}\left(\sum_{j=1}^n\left|\frac{\partial W}{\partial x_j}\right|^2+\sum_{k=1}^r\left|W_k\right|^2\right)-\frac{s}{(1+s)^2}{u}^{-1-s}\left|\sum_{l=1}^{n+r}{z}_l \frac{\partial W}{\partial z_l}\right|^2.
\end{align*}
From the complex Cauchy-Schwarz inequality,
\[
\left|\sum_{l=1}^{n+r}{z}_l \frac{\partial W}{\partial z_l}\right|^2\leq {u}\;\left|\sum_{l=1}^{n+r} \frac{\partial W}{\partial z_l}\right|^2,
\]
so we have
\begin{align*}
&|\nabla^{\pcU} W|_{\widetilde{g}_s}^2 \geq\frac{1}{(1+s)^2}u^{-s}\sum_{j=1}^n\left|\frac{\partial W}{\partial x_j}\right|^2+\frac{1}{1+s}u^{-s}\sum_{k=1}^r\left|W_k\right|^2.
\end{align*}
Restricting the above inequality on ${\mu^{-1}_s(-1)}$, we obtain
\begin{align*}
&|\nabla W|_{\mathfrak{g}_s}^2 \geq\frac{1}{(1+s)^2}(u|_{\mu^{-1}_s(-1)})^{-s}\sum_{j=1}^n\left|\left.\frac{\partial W}{\partial x_j}\right|_{\mu^{-1}_s(-1)}\right|^2+\frac{1}{1+s}(u|_{\mu^{-1}_s(-1)})^{-s}\sum_{k=1}^r\left|\left.W_k\right|_{\mu^{-1}_s(-1)}\right|^2.
\end{align*}

Now we use the fact that $W_1,\ldots,W_r$ gives a smooth complete intersection. We set $S^{2n-1}_2\subset \mathbb{C}^n$ as the sphere of radius $2$;
$$
S_2^{2n-1}:=\{(x_1,\ldots,x_n)\in \bC^n:|x_1|^2+\cdots + |x_n|^2 =4\}.
$$
If $(z_1,\ldots,z_{n+r})\in \mu^{-1}_s(-1)$ and $(z_1,\ldots,z_n)\in S_{2}^{2n-1}$, $(z_{n+1},\ldots,z_{v+r})$ cannot be zero since
$$
(z_{n+1},\ldots,z_{v+r})=0 \Rightarrow |x_1|^2+\cdots+|x_n|^2 =u=w,
$$
and on $\mu^{-1}_s(-1)$, we have
\[
u^s w = \frac{1}{1+s} \Rightarrow u=\left(1+s\right)^{-\frac{1}{1+s}}<1.
\]
When $W_1=\cdots = W_r=0$, the Jacobian matrix
$$
\mathbf{J}:=\begin{pmatrix}
\frac{\partial W_1}{\partial x_1} & \cdots & \frac{\partial W_r}{\partial x_1} \\
\vdots & \ddots & \vdots \\
\frac{\partial W_1}{\partial x_n} &\cdots & \frac{\partial W_r}{\partial x_n} 
\end{pmatrix}
$$
has rank $r$ and thus defines an injective map. 
We denote by $\mathbf{v}$ the column vector 
\[\mathbf{v}:=(z_{n+1},\ldots,z_{n+r})^T.\] Thus, in the compact set $(S^{2n-1}_2\times S_1^{2n-1})\cap \{W_1=\cdots=W_r=0\}\subset \bC^{n+r}$, we have a minimum $a>0$:
\begin{align*}
a&:=\min\left\{\sum_{j=1}^n \left|\frac{\partial W}{\partial x_j}(z_1,\ldots,z_{n+r})\right|^2:(z_{n+1},\ldots,z_{n+r})\in S_{1}^{2r-1}, (z_1,\ldots,z_n)\in S^{2n-1}_2, W_1=\cdots =W_r=0\right\}, \\
&= \min\left\{\frac{| \mathbf{J}\mathbf{v}|_{\mathrm{Euc}}^2}{|\mathbf{v}|_{\mathrm{Euc}}^2}:(z_{n+1},\ldots,z_{n+r})\in \bC^r\setminus\{0\}, (z_1,\ldots,z_n)\in S^{2n-1}_2, W_1=\cdots=W_r=0\right\}.
\end{align*}
Here, $|-|_{\mathrm{Euc}}$ denotes the usual Euclidean norm. In $S^{2n-1}_2$, we define $\mathcal{W}$ as
$$
\mathcal{W}:=\left\{(z_1,\ldots,z_{n})\in S^{2n-1}_2: \sum_{j=1}^n \left|\frac{\partial W}{\partial x_j}(z_1,\ldots,z_n)\right|^2>\frac{a}{2} \textrm{ for all $(z_{n+1},\ldots,z_{n+r})\in S_1^{2r-1}$}\right\}.
$$
Then $W$ contains $\{W_1=\cdots = W_r=0\}$. Outside $\mathcal{W}$, $|W_1|^2+\cdots +|W_r|^2$ is never zero. Let $b$ be the minimum of these values in $S^{2n-1}_2\setminus \mathcal{W}$,
$$
b:=\min\left\{|W_1|^2+\cdots + |W_r|^2:(z_1,\ldots,z_n)\in S^{2n-1}_2\setminus\mathcal{W}\right\}.
$$
Since $S^{2n-1}_2\setminus \mathcal{W}$ is compact, the minimum $b$ is positive.

Now let $t$ be any real number with 
$
t\geq 1,
$ and consider an arbitrary point $(z_1,\ldots,z_{n+r})\in \mu^{-1}_s(-1)$ with $|z_1|^2+\cdots+|z_n|^2 = 4t^2$. We let $(z_1',\ldots,z_n') \in S^{2n-1}_2$ be such that $t z'_j = z_j$, for $j=1,\ldots,n.$ When $(z_1',\ldots,z_n')\in \mathcal{W}$, together with \eqref{eq:radial_coordinate},
\begin{align}
\sum_{j=1}^n \left|\frac{\partial W}{\partial x_j}(z_1,\ldots,z_{n+r})\right|^2 =& \sum_{j=1}^n \left|\sum_{k=1}^r z_{n+k} \frac{\partial W_k}{\partial x_j}(z_1,\ldots,z_n)\right|^2 \nonumber\\
=& \sum_{j=1}^n \left|\sum_{k=1}^r t^{d_k-1}z_{n+k} \frac{\partial W_k}{\partial x_j}(z_1',\ldots,z_n')\right|^2\nonumber\\
\geq& \frac{a}{2}\sum_{k=1}^r \left|t^{d_k-1}z_{n+k}\right|^2\geq \frac{a}{2}t^{2d_{\min}-2}\sum_{k=1}^r\left|z_{n+k}\right|^2. \label{ineq1}
\end{align}
For $(z_1,\ldots,z_{n+r})\in \mu^{-1}_s(-1)$, we have
\[
|z_1|^2+\cdots + |z_n|^2 = \frac{(u|_{\mu^{-1}_s(-1)})^{-s}}{1+s}+d_1|z_{n+1}|^2+\cdots + d_r|z_{n+r}|^2,
\]
which implies that
\begin{align*}
u|_{\mu^{-1}_s(-1)} &= \frac{(u|_{\mu^{-1}_s(-1)})^{-s}}{1+s}+(d_1+1)|z_{n+1}|^2+\cdots + (d_r+1)|z_{n+r}|^2\\
&\leq \frac{(u|_{\mu^{-1}_s(-1)})^{-s}}{1+s}+(d_{\max}+1)(|z_{n+1}|^2+\cdots + |z_{n+r}|^2).
\end{align*}
We may assume that $u|_{\mu^{-1}_s(-1)}$ is large enough so that
\[
(u|_{\mu^{-1}_s(-1)})^{-s} \leq u|_{\mu^{-1}_s(-1)},
\]
and we get
\begin{align}\label{ineq2}
\frac{s}{(d_{\max}+1)(1+s)}u|_{\mu^{-1}_s(-1)} \leq |z_{n+1}|^2+\cdots + |z_{n+r}|^2.
\end{align}
Putting \eqref{ineq1} and \eqref{ineq2} together, we obtain
\[
\sum_{j=1}^n \left|\frac{\partial W}{\partial x_j}(z_1,\ldots,z_{n+r})\right|^2\geq \frac{a}{2}t^{2d_{\min}-2}\sum_{k=1}^r\left|z_{n+k}\right|^2\geq\frac{as}{2(1+s)(d_{\max}+1)} t^{2d_{\min}-2}u|_{\mu^{-1}_s(-1)}
\]
When $(z_1',\ldots,z'_n)$ is outside $\mathcal{W}$,
\begin{align*}
\sum_{k=1}^r |W_k(z_1,\ldots,z_n)|^2 =& \sum_{k=1}^r t^{2d_k}|W_k(z_1',\ldots,z_n')|^2\\
\geq& \sum_{k=1}^r t^{2d_{\min}}|W_k(z_1',\ldots,z_n')|^2\geq t^{2d_{\min}}b
\end{align*}
The value of $t$ and that of $u|_{\mu^{-1}_s(-1)}$ is related by the following inequality: 
\begin{align*}
4t^2 &= |z_1|^2 +\cdots + |z_n|^2 \\
&= \frac{1}{2}\left(|z_1|^2 +\cdots + |z_n|^2+d_1|z_{n+1}|^2 +\cdots + d_r|z_{n+r}|^2+\frac{1}{1+s}(u|_{\mu^{-1}_s(-1)})^{-s}\right)\\
&\geq \frac{1}{2}(|z_1|^2 +\cdots + |z_{n+r}|^2)\\
&=\frac{1}{2}u|_{\mu^{-1}_s(-1)}.
\end{align*}
Thus, 
$$
|\nabla W|^2_{\mathfrak{g}_s} \geq \min\left(\frac{as}{2\cdot 8^{d_{\min}-1}(1+s)^3(d_{\max}+1)},\frac{b}{8^{d_{\min}}(1+s)}\right)(u|_{\mu^{-1}_s(-1)})^{d_{\min}-s}.
$$

Next, we prove 
\[
|\nabla^k W|_{\mathfrak{g}_{s}}= O\left((u|_{\mu^{-1}_s(-1)})^{\frac{1}{2}({d_{\max}+1}-k{(1+s)})}\right).
\]
By Propositions~\ref{prop:tensor_compariton} and \ref{prop:shrinking}, it is enough to prove that the following is true on each open subset $V_j$
\[
|\nabla^\pazocal{C} (\Phi_* W)|_{\mathfrak{g}_{s,\pazocal{C}}}=O\left((u|_{V_j})^{\frac{d_{\max}+1}{2}-k\frac{1+s}{2}}\right).
\]
We only prove that above is true on $V_1$:
\[
|\nabla^\pazocal{C} (\Phi_* W)|_{\mathfrak{g}_{s,\pazocal{C}}}=O\left(\widehat{u}^{\frac{d_{\max}+1}{2}-k\frac{1+s}{2}}\right).
\]
We have
\begin{align*}
&|(\nabla^\pazocal{C})^k (\Phi_* W)|_{{g}_\pazocal{C}} \\
&\leq \sum_{j=1}^r |(\nabla^\pazocal{C})^k (\Phi_* P_jW_j)|_{\mathfrak{g}_{s,\pazocal{C}}} \\
&\leq \sum_{\substack{1\leq j\leq r\\l+k_1+\cdots+k_k=k\\
2\leq q_1,\cdots,q_{l}\leq n+r}}{k \choose l, k_1,\ldots,k_{m}}\left|\frac{\partial^{l}(P_jW_j)}{\partial Z_{q_1}\cdots\partial Z_{q_{l}}}\right|\left|(\nabla^\pazocal{C})^{k_1} dZ_{q_1}\right|_{\mathfrak{g}_{s,\pazocal{C}}}\cdots\left|(\nabla^\pazocal{C})^{k_l}dZ_{q_l}\right|_{\mathfrak{g}_{s,\pazocal{C}}}
\end{align*}
From Propositions~\ref{prop:dZ_asymp} and \ref{cor:cov_asymp}, we obtain
\[
\left|(\nabla^\pazocal{C})^{k}dZ_{j}\right|_{\mathfrak{g}_{s,\pazocal{C}}} = O(\widehat{u}^{-\frac{1}{2}(k(1+s)+\ch(z_j)+s)})
\]
for each $j=1,\ldots,n+r$ and $k=0,1,\ldots$. Obviously, $P_jW_j \in \widehat{R}_{-\ch(z_{n+j})+1}$, where $\widehat{R}$ is defined in Definition~\ref{def:degree1}, and by Proposition~\ref{differentiation1},
\[
\left|\frac{\partial^{l}(P_jW_j)}{\partial Z_{q_1}\cdots\partial Z_{q_{l}}}\right| = O(\widehat{u}^{\frac{1}{2}(\ch(z_{q_1})+\cdots+\ch(z_{q_l})-\ch(z_{n+j})+1-l)}).
\]
Putting these all together, we get
\begin{align*}
\left|\frac{\partial^{l}(P_jW_j)}{\partial Z_{q_1}\cdots\partial Z_{q_{l}}}\right|\left|(\nabla^\pazocal{C})^{k_1} dZ_{q_1}\right|_{\mathfrak{g}_{s,\pazocal{C}}}\cdots\left|(\nabla^\pazocal{C})^{k_l}dZ_{q_l}\right|_{\mathfrak{g}_{s,\pazocal{C}}}=O(\widehat{u}^{\frac{1}{2}(-\ch(z_{n+j})+1-k(1+s))})
\end{align*}
Since $\ch(z_{n+j})=-d_j$, we have
\[
|(\nabla^\pazocal{C})^k (\Phi_* W)|_{\mathfrak{g}_{s,\pazocal{C}}} = \sum_{j=1}^{n+r} O(\widehat{u}^{\frac{1}{2}(d_j+1-k(1+s))})= O(\widehat{u}^{\frac{1}{2}(d_{\max}+1-k(1+s))}).
\]
\end{proof}

\begin{remark}\label{remark3.40}
When $2d_{\min}< d_{\max}$, there are some counterexamples $W_1,\ldots,W_r$ such that $W$ is not strongly elliptic. These counterexamples occur when hypersurfaces defined by each homogeneous polynomials $W_1,\ldots, W_r$ contain some singularities that are not included in $V(W_1,\ldots,W_r)$. 
\end{remark}
{
\begin{remark}\label{remark3.41}
The condition $d_{\max} \leq 2d_{\min}$ required for the strong ellipticity of the holomorphic function $W:\cX \to \bC$ is not a fundamental geometric obstruction inherent to the Calabi-Yau complete intersection $V(\ud{W})$, but rather an artifact of the specific $U(1)$-invariant K\"ahler potential $u^{1+s}$ chosen for the ambient metric $\mathfrak{g}_s$. This particular potential treats all coordinate directions symmetrically, which leads to unbalanced radial growth rates in the covariant derivatives of $W$ when the degrees $d_k, \ k=1, \cdots, r$ of the defining homogeneous polynomials vary significantly. Consequently, if the maximum degree exceeds twice the minimum degree, the gradient decay of the lower-degree components fails to compensate for the higher-order derivative bounds of the higher-degree components, violating the strong ellipticity criteria. It is plausible that this restriction could be relaxed or eliminated by employing a more general, anisotropic family of $U$(1)-invariant K\"ahler potentials $\widetilde{u}(\ud{z})$. By carefully tuning $\widetilde{u}(z)$, one might equalize the asymptotic behavior of $|\nabla W|_{\mathfrak{g}_s}^k$ and $|\nabla^k W|_{\mathfrak{g}_s}$. 
Such a modification would necessitate a corresponding adjustment in the K\"ahler reduction process and the subsequent asymptotic estimates, offering a promising direction for future research to extend the applicability of the $L^2$-Hodge theory to all smooth projective Calabi-Yau complete intersections regardless of their multidegree ratios.
\end{remark}
}
\section{Computation of the cohomology}
\label{section4}
In this section, we prove Theorem \ref{Main2}; it is a generalization of \cite{Paper1} which deals with the case $r=1$. The basic strategies for the current proof and \cite{Paper1} are the same, but they contain much heavier technical computations due to the fact that $r>1$.
{ The computation of the cohomology is organized into three distinct phases. First, in Section \ref{section4.1}, we evaluate the second page of the algebraic spectral sequence, denoted by $\tilde{E}_2^{p,q}$, by utilizing the Koszul-type complex computations of Adolphson and Sperber \cite{AS}. Second, in Section \ref{section4.2}, we establish an isomorphism between the analytic spectral sequence $E_2^{p,q}$ and the algebraic spectral sequence $\tilde{E}_2^{p,q}$ at the second page, applying the GAGA principle of Serre \cite{SerreGAGA} to bridge the analytic and algebraic domains. Finally, in Section \ref{section4.3}, we analyze the higher differentials $d_k$ for $k \geq 2$, demonstrating that the sequences degenerate after the $r$-th page to yield the final isomorphic cohomology groups.}

The primitive cohomology of $V(\ud{W})$ was calculated by Griffiths (r = 1) \cite{PG}, Terasoma ($d_1=\cdots=d_r$) \cite{Terasoma} and Konno (general cases) \cite{Konno}; moreover, see \cite{Dimca} for a nice summary of results. The primitive cohomology is described by the charge 0 component of the Jacobian ring $R(W)$. The Jacobian ring of $W$ is defined as follows:
$$
R(W) := \frac{\bC[\ud{x},\ud{p}]}{\left(\frac{\partial W}{\partial x_1},\ldots, \frac{\partial W}{\partial x_n}, \frac{\partial W}{\partial p_1},\ldots, \frac{\partial W}{\partial p_r}\right)}.
$$
We introduce the \emph{charge} and \emph{weight}, which are additive gradings on $\bC[\ud{x},\ud{p}]$,
\begin{equation}
\begin{aligned}
\ch(x_1)&=\cdots=\ch(x_n)=1, & \qquad \ch(p_j)&=-d_j,\quad j=1,\ldots,r,\\
\wt(x_1)&=\cdots=\wt(x_n)=0, & \wt(p_1)&=\cdots=\wt(p_r)=1.
\end{aligned}
\label{definition:chandwt}
\end{equation}
We write $\bC[\ud{x},\ud{p}]_j$ for a charge $j$ component of $\bC[\ud{x},\ud{p}]$, and $\bC[\ud{x},\ud{p}]_{(k)}$ for the weight $k$ component. We denote $\bC[\ud{x},\ud{p}]_{j,(k)}$ by $\bC[\ud{x},\ud{p}]_j\cap \bC[\ud{x},\ud{p}]_{(k)}$.
From the theory of toric varieties, if $\cM$ is a quasi-coherent sheaf of $\cO_{\cX}$-module, there is a corresponding graded $\bC[\ud{x},\ud{p}]$-module $M$ whose grading is compatible with the charge grading on $\bC[\ud{x},\ud{p}]$ \cite[Section 5.3]{Cox}. We call such grading the charge grading on $M$, and write $M_k$ for the charge $k$ component in $M$. Then, $M_0$ corresponds to the global sections of $\cM$. $M$ is also given the weight grading corresponding to the weight grading of $\bC[\ud{x},\ud{p}]$ arising from the filtration
$$
M\supseteq (p_1,\ldots,p_r)M \supseteq (p_1,\ldots,p_r)^2M \supseteq \cdots
$$
We write $M_{(k)}$ for the weight $k$ component of $M$. We denote $M_{j,(k)}$ by $M_{j}\cap M_{(k)}$. 

The charge and weight gradings on $\bC[\ud{x},\ud{p}]$ naturally induce gradings of the charge and weight on $R(W)$. The charge $0$ component is
$$
R(W)_{0} = \frac{\bC[\ud{x},\ud{p}]_{0}}{\left(\frac{\partial W}{\partial x_1},\ldots, \frac{\partial W}{\partial x_n}, \frac{\partial W}{\partial p_1},\ldots, \frac{\partial W}{\partial p_r}\right)\cap \bC[\ud{x},\ud{p}]_{0}}.
$$
The non-primitive cohomology can be computed using the weak Lefschetz theorem. In summary, we have the following result:
$$
H^p(V(\ud{W});\bC) \cong \begin{cases}
\bC & \textrm{if $0\leq p\leq 2(n-r-1)$, $p\neq n-r-1$, $p$ is even,} \\
R(W)_{0} & \textrm{if $p=n-r-1$, $n-r-1$ is odd,}\\
R(W)_{0} \oplus \bC & \textrm{if $p=n-r-1$, $n-r-1$ is even,}\\
0 & \textrm{otherwise.}
\end{cases}
$$

In the previous sections, we considered $\cX$ as a complex manifold. In this section, we will consider $\cX$ as an algebraic variety over $\bC$, endowed with the Zariski topology and the structure sheaf $\cO_{\cX}$ whose value on an open subset $U\subset \cX$ is an algebraic (in other words, polynomial) function on $U$. The complex manifold we have denoted by $\cX$ in Section \ref{section3} will be denoted by $\cX^{an}$, which is the analytification of $\cX$; see \cite[p. 7]{SerreGAGA}. Then $\cX^{an}$ is also a ringed space whose topology is the Euclidean topology, and the stalk $(\cO_{\cX}^{an})_{[\ud{x},\ud{p}]}$ at $[\ud{x},\ud{p}]\in \cX^{an}$ consists of germs of holomorphic functions at $[\ud{x},\ud{p}]$. Consequently, we will distinguish the projective variety $\bP^{n-1}_{\bC}$ and the compact complex manifold $(\bP^{n-1}_{\bC})^{an}=\bC\bP^{n-1}$. 
We are given a morphism $h:(\cX^{an},\cO_{\cX}^{an})\to (\cX,\cO_{\cX})$ of ringed spaces, which maps a point in $\cX^{an}$ to its corresponding maximal ideal in $\cX$. There is an analytification functor $(-)^{an}:QCoh(\cO_{\cX})\to QCoh(\cO_{\cX}^{an})$ between the category of quasi-coherent sheaves on $(\cX,\cO_{\cX})$ and $(\cX^{an},\cO_{\cX}^{an})$ defined as
$$
\cF^{an} = h^{-1} \cF \otimes_{h^{-1}\cO_{\cX}}\cO_{\cX}^{an}.
$$
See \cite[p. 16]{SerreGAGA} for more details.

Using the spectral sequence of the double complex where the vertical differential is $\overline{\partial}$ and the horizontal differential is $dW\wedge$, we have a spectral sequence whose 0th, 1st, and 2nd pages are
\begin{align*}
E_0^{p,q} &= \cA^{p,q}(\cX^{an}),\\
E_1^{p,q} &= H^q(\cA^{p,\bullet}(\cX^{an}),\overline{\partial}),\\
E_2^{p,q} &= H^p(H^q(\cA^{\bullet,\bullet}(\cX^{an}),\overline{\partial}),dW),
\end{align*}
where we recall the definition of $\cA^{p,q}(\cX^{an})$ in \eqref{spaceofsmoothforms}. We let $\Omega_{\cX}:=\Omega_{\cX/\bC}$ be a sheaf of K\"ahler differentials on $\cX$ which is a coherent sheaf of $\cO_{\cX}$-modules. Its analytification $(\Omega_{\cX})^{an}$ is a sheaf of holomorphic sections on $\cX^{an}$. Accordingly, $\Omega_{\cX}^p$ ($(\Omega_{\cX}^p)^{an}$, respectively) is a sheaf of algebraic (holomorphic, respectively) differential forms defined by the exterior powers of $\Omega_{\cX}$ (($\Omega_{\cX})^{an}$, respectively):
$$
\Omega_{\cX}^p := \bigwedge^p \Omega_{\cX}.
$$
By Dolbeault's theorem, we have
$$
H^q(\cA^{p,\bullet}(\cX^{an}),\overline{\partial}) \cong H^q(\cX^{an},(\Omega_{\cX}^{p})^{an}).
$$
Hence we have
\begin{align*}
E_1^{p,q} &\cong H^q(\cX^{an},(\Omega_{\cX}^{p})^{an}),\\
E_2^{p,q} &\cong H^p(H^q(\cX^{an},(\Omega_{\cX}^{\bullet})^{an}),dW).
\end{align*}

The cohomology is easier to manipulate in the algebraic setting; we have another spectral sequence of the double complex where the vertical differential is the \v{C}ech differential, and the horizontal differential is defined as $dW\wedge$. $\cU:=\{U_1,\ldots,U_n\}$ is an open cover of $\cX$ defined in Definition \ref{volumeform}. Then the 0th page of the spectral sequence is defined by the \v{C}ech complex:
\begin{align*}
\widetilde{E}_0^{p,q} &= \check{C}^{q}(\cU,\Omega_{\cX}^p),\\
\widetilde{E}_1^{p,q} &= H^q(\cX,\Omega^{p}_{\cX}),\\
\widetilde{E}_2^{p,q} &= H^p(H^q(\cX,\Omega^{\bullet}_{\cX}),dW).
\end{align*}
In Section \ref{section4.1}, we compute the second page of this spectral sequence. We use the work of Adolphson and Sperber \cite{AS} to derive the following lemma.
\begin{lemma} \label{ETWO}
$$
H^p(H^q(\cX,\Omega_{\cX}^{\bullet}),dW)\cong\begin{cases}
\bC & \textrm{if $1\leq p=q\leq n-1$ or $q=0,p=2k-1,2\leq k\leq r$,}\\
R(W)_{0} & \textrm{if $q=0, p=n+r-1$}, \\
0 &\textrm{otherwise.}
\end{cases}
$$
\end{lemma}
The cohomologies of the sheaf of holomorphic differential forms $(\Omega^{p}_{\cX})^{an}$ and the sheaf of algebraic differential forms $\Omega^{p}_{\cX}$ are not the same in general, because $\cX$ is not compact. We use Lemma \ref{ETWO} (the finite-dimensionality of $\widetilde{E}_2^{p,q}$) to prove the following lemma, which connects the algebraic cohomology with the analytic cohomology, i.e. $E_2^{p,q}\cong \widetilde{E}_2^{p,q}$:
\begin{lemma}
\label{lemma5.2}
For each $p,q$, we have an isomorphism
$$H^p(H^q(\cX^{an},(\Omega_{\cX}^{\bullet})^{an}),dW) \cong H^p(H^q(\cX,\Omega_{\cX}^{\bullet}),dW).$$
\end{lemma}

In Section \ref{section4.3}, we compute the higher differentials $d_k^{p,q}:E_k^{p,q} \to E_k^{p+k,q-k+1} $ for $k>2$. We show that for each $p=1,2,\ldots,r-1$, $d_{p+1}^{p,p}$ is an isomorphism, and all other differentials $d_k^{p,q}$ are trivial for $k\geq 2$ (Lemma \ref{higherdiffan}). Thus, if we take the cohomology along the higher differentials, we get the following isomorphisms:
$$
{E}_{\infty}^{p,q}  \cong \begin{cases}
\bC & \textrm{if $r\leq p=q \leq n-1$,}\\
R(W)_{0} & \textrm{if $q=0$, $p=n+r-1$,} \\
0 & \textrm{otherwise.}
\end{cases}
$$
Since all the differentials are trivial after the $r$th page, we obtain the following theorem, which proves Theorem \ref{Main2}.
\begin{theorem}
\label{cohomologycmp}
$$
H^{p}(PV^\bullet(\cX),\overline{\partial}_W) \cong \begin{cases}
\bC & \textrm{if $p=-n+r+1+2k\neq 0$, $0\leq k \leq n-r-1$}\\
R(W)_{0} & \textrm{if $p=0$ and $n-r-1$ is odd,}\\
R(W)_{0}\oplus \bC & \textrm{if $p=0$ and $n-r-1$ is even,}\\
0 & \textrm{otherwise.}
\end{cases}
$$
\end{theorem}

\subsection{Computation of the second page of the spectral sequence}
\label{section4.1}

We prove Lemma 
\ref{ETWO}. Using the Leray spectral sequence for the natural projection map
\begin{eqnarray*}
    \pi:=\pi_{CY}:\cX \to \bP_{\C}^{n-1},
\end{eqnarray*}
we have the following lemma.
\begin{proposition}
\label{leray}
For any sheaf $\cF$ of $\cO_{\cX}$-modules, we have an isomorphism
\begin{align*}
H^q(\cX,\cF) \simeq H^q(\bP^{n-1}_{\bC},\pi_* \cF).
\end{align*}
\end{proposition}
\begin{proof}
It follows directly from the fact that $\pi$ is an affine morphism and \cite[Lemma 01F4]{Stacks}.
\end{proof}

We define $\cO_{\cX}(d) := \pi^*\cO_{\bP^{n-1}_{\bC}}(d)$ for an integer $d$.
\begin{proposition}[Euler's sequence] \label{lone}
We have a short exact sequence of $\cO_{\cX}$-modules,
\begin{align}
0\to \Omega_{\cX} \to \{\cO_{\cX}(-1)\}^{\oplus n}\oplus \bigoplus_{k=1}^r\cO_{\cX}(d_k)\to \cO_{\cX} \to 0, \label{eulersequence1} \\
0\to \Omega_{\cX}^{p} \to \bigwedge^p \left(\{\cO_{\cX}(-1)\}^{\oplus n}\oplus \bigoplus_{k=1}^r\cO_{\cX}(d_k)\right)\to \Omega_{\cX}^{p-1} \to 0. \label{eulersequence2}
\end{align}
\end{proposition}
\begin{proof}
The short exact sequence \eqref{eulersequence2} is induced by taking the exterior powers of the sequence \eqref{eulersequence1}, together with the fact that $\cO_{\cX}$ is a locally free sheaf. We only need to prove \eqref{eulersequence1}. This is done by defining the exact sequence locally: on $U_j$ defined by $x_j\neq 0$, define an affine coordinate $X_k=x_k/x_j$ for $k=1,\ldots,n$ and $P_l = x_j^{d_l} p_l$ for $l=1,\ldots,r$. Then $\Omega_{\cX}(U_j)$ is a free $\cO(U_j)$-module of rank $n+r-1$ whose basis is $\{dX_1,\ldots,dX_{j-1},dX_{j+1},\ldots,dX_n,dP_1,\ldots,dP_r\}$. We denote an element $(f_1,\cdots,f_n,g_1,\cdots,g_r)\in \left(\{\cO_{\cX}(-1)\}^{\oplus n}\oplus \bigoplus_{k=1}^r\cO_{\cX}(d_k)\right)(U_j)$ by
$$
f_1 dx_1 + \cdots + f_n dx_n + g_1 dp_1 + \cdots + g_r dp_r.$$
The map
$$
\Omega_{\cX}(U_j) \to \left(\{\cO_{\cX}(-1)\}^{\oplus n}\oplus \bigoplus_{k=1}^r\cO_{\cX}(d_k)\right)(U_j)
$$
is then defined as
$$
dX_k \mapsto \frac{1}{x_j} dx_k - \frac{x_k}{x_j^2} dx_j, \quad dP_l \mapsto x_j^{d_l}dp_l + d_l x_j^{d_l-1}p_l dx_j,
$$
for $k=1,\ldots,n$ and $l=1,\ldots,r$. The map
$$
\left(\{\cO_{\cX}(-1)\}^{\oplus n}\oplus \bigoplus_{k=1}^r\cO_{\cX}(d_k)\right)(U_j)\to \cO_{\cX}(U_j)
$$
is defined as
\begin{align}
\label{eulersequencesurj}
dx_k \mapsto x_k, \quad dp_l \mapsto -d_l\cdot p_l,
\end{align}
for $k=1,\ldots,n$ and $l=1,\ldots,r$. We leave it as an exercise to the reader to verify that these maps induce morphisms between sheaves of $\cO_{\cX}$-modules and that the sequence \eqref{eulersequence1} is exact.
\end{proof}
Since $\pi$ is an affine morphism, we have the same exact sequences after taking the direct image functor $\pi_*$ to \eqref{eulersequence1} and \eqref{eulersequence2}, i.e. we have the exact sequences of $\cO_{\bP^{n-1}_\bC}$-modules,
\begin{align}
0\to \pi_*\Omega_{\cX} \to \pi_*\{\cO_{\cX}(-1)\}^{\oplus n}\oplus \bigoplus_{k=1}^r\cO_{\cX}(d_k)\to \pi_*\cO_{\cX} \to 0, \notag\\
0\to \pi_*\Omega_{\cX}^{p} \to \bigwedge^p \left(\{\pi_*\cO_{\cX}(-1)\}^{\oplus n}\oplus \bigoplus_{k=1}^r\pi_*\cO_{\cX}(d_k)\right)\to \pi_*\Omega_{\cX}^{p-1} \to 0. \label{eulersequence3}
\end{align}
The sheaf $\pi_* \cO_{\cX}(d)$ over $\bP^{n-1}_{\bC}$ can be computed as follows. 
\begin{proposition}
For any $d \in \Z$, we have
$$
\pi_* \cO_{\cX}(d) \cong \bigoplus_{k_1,\ldots,k_r\geq 0} \cO_{\bP_{\bC}^{n-1}}\left(\sum_{\alpha=1}^rd_\alpha k_\alpha +d \right).
$$
\label{lemma5.6}
\end{proposition}
\begin{proof}
We compute $\pi_* \cO_{\cX}(d)$ on the affine open subset $V_j$ of $\bP_{\bC}^{n-1}$ defined by $x_j\neq 0$:
\begin{align*}
\pi_* \cO_{\cX}(d) (V_j) &=  \cO_{\cX}(d)(\pi^{-1}(V_j)) \\
&= \bC\left[x_1,\ldots,x_n,p_1,\ldots,p_r,\frac{1}{x_j}\right]_{d} \\
&\cong \bigoplus_{k_1,\ldots,k_r\geq 0} \bC\left[x_1,\ldots,x_n,\frac{1}{x_j}\right]_{\sum_{\alpha=1}^r d_\alpha k_\alpha+d} p_1^{k_1}\cdots p_r^{k_r} \\
&\cong \bigoplus_{k_1,\ldots,k_r\geq 0} \cO_{\bP_{\bC}^{n-1}}\left(\sum_{\alpha=1}^rd_\alpha k_\alpha +d \right)(V_j).
\end{align*}
\end{proof}

\begin{proposition}\cite[{Lemma 01XT}]{Stacks}\label{lthree}
Let $n$ be a positive integer. We have
$$
H^q (\bP^{n-1}_{\bC}, \cO_{\bP_{\bC}^{n-1}}(d)) \cong \begin{cases}
\bC[\ud x]_{d} & \textrm{if } q=0, d\geq 0, \\
\left(\frac{1}{x_1\cdots x_n}\bC\left[\frac{1}{x_1},\ldots,\frac{1}{x_n}\right]\right)_{d} & \textrm{if }  q=n-1, d\leq  -n, \\
0 & \textrm{otherwise}. 
\end{cases}
$$
\label{lemma5.7}
\end{proposition}

\begin{proof}[Proof of Lemma \ref{ETWO}]
By Propositions \ref{leray}, \ref{lemma5.6}, and \ref{lthree}, the $p$th sheaf cohomology of the middle term of \eqref{eulersequence3} is nonzero only if either $q=0$ or $q=n-1, p=n$:
\begin{align*}
H^q\left(\cX,\bigwedge^p \left(\{\cO_{\cX}(-1)\}^{\oplus n}\oplus \bigoplus_{k=1}^r\cO_{\cX}(d_k)\right)\right) \cong 
\begin{cases}
\left(\Omega_{\bC[\ud{x},\ud{p}]/\bC}^p\right)_{0} & \textrm{if $q=0$,} \\
\bC & \textrm{if $q=n-1$, $p=n$,} \\
0 & \textrm{otherwise.}
\end{cases}
\end{align*}
Recall that the charge is defined on the graded $\bC[\ud{x},\ud{p}]$-module which corresponds to a sheaf of $\cO_{\cX}$-module \cite[Section 5.2]{Cox}. The module $\Omega_{\bC[\ud{x},\ud{p}]/\bC}^p$ is a graded $\bC[\ud{x},\ud{p}]$-module which corresponds to the sheaf $$\bigwedge^p \left(\{\cO_{\cX}(-1)\}^{\oplus n}\oplus \bigoplus_{k=1}^r\cO_{\cX}(d_k)\right),$$
and the charge grading on $\Omega_{\bC[\ud{x},\ud{p}]/\bC}^p$ is described by the following equality:
$$
\ch(x_j)= \ch(dx_j)=1, \qquad \ch(q_k)=\ch(dq_k)=-d_j,
$$
for $j=1,\ldots,n$, $k=1,\ldots,r$. 

Applying the sheaf cohomology to the Euler sequences (Proposition \ref{lone}) for $q=0$, we have
$$
0\to H^0(\cX,\Omega_{\cX}^{p})\to(\Omega_{\bC[\ud{x},\ud{p}]/\bC}^p)_{0}\to H^0(\cX,\Omega_{\cX}^{p-1})\to H^1(\cX,\Omega_{\cX}^{p}) \to 0.
$$
When $p=n$, we have
\begin{align*}
0\to H^{n-2}(\cX,\Omega_{\cX}^{n-1})\to H^{n-1}(\cX,\Omega_{\cX}^{n}) \to \bC \to H^{n-1}(\cX,\Omega_{\cX}^{n-1}) \to H^{n}(\cX,\Omega_{\cX}^{n}) \to 0.
\end{align*}
If $q\neq 0$, $(q,p)\neq (n-1,n-1)$, and $(q,p)\neq (n-2,n-1)$, then we have
\begin{align}
H^q(\cX,\Omega_{\cX}^{p}) \cong H^{q+1}(\cX,\Omega_{\cX}^{p+1}).
\label{diagiso}
\end{align}
If $p\geq n+r$, $\Omega_{\cX}^{p}=0$, so for $p$ large enough, $H^q(\cX,\Omega^{p}_{\cX})=0$. Thus, we have the following computation:
\begin{align}
\tilde E_1^{p,q}=H^q(\cX,\Omega_{\cX}^{p}) \cong \begin{cases}
\bC & \textrm{if $p=q$, $1\leq p \leq n-1,$} \\
H^0(\cX,\Omega_{\cX}^{p}) & \textrm{if $q=0$, $0\leq p \leq n+r-1$,} \\
0 & \textrm{otherwise}.
\end{cases}
\label{algcohpage1}
\end{align}
To take the cohomology along $dW$, we again consider the Euler sequence (Proposition \ref{lone}).
\begin{equation}
\label{eulersequencedW}
\begin{tikzcd}
0 \arrow[r] & {H^0(\cX,\mathcal{O}_{\cX})} \arrow[r, "dW"] \arrow[d, "\cong"] & {H^0(\cX,\Omega_{\cX}^{1})} \arrow[r, "dW"] \arrow[d, hook]                   & {H^0(\cX,\Omega_{\cX}^{2})} \arrow[r, "dW"] \arrow[d, hook]                              & \cdots \\
0 \arrow[r] & {\mathbb{C}[\underline{x},\underline{p}]_0} \arrow[r, "dW"]                       & {(\Omega_{\mathbb{C}[\underline{x},\underline{p}]/\mathbb{C}}^1)_0} \arrow[r, "dW"] \arrow[d] & {(\Omega_{\mathbb{C}[\underline{x},\underline{p}]/\mathbb{C}}^2)_0} \arrow[r, "dW"] \arrow[d, two heads] & \cdots \\
            &                                                                                   & {H^0(\cX,\mathcal{O}_{\cX})} \arrow[d, two heads] \arrow[r, "dW"]           & {H^0(\cX,\Omega_{\cX}^{1})} \arrow[r, "dW"]                                              & \cdots \\
            &                                                                                   & {H^1(\cX,\Omega_{\cX}^{1})\cong\mathbb{C}}                                    &                                                                                                          &       
\end{tikzcd}
\end{equation}
where $(\Omega^1_{\bC[\ud{x},\ud{p}]/\bC})_0\to {H^0(\cX,\mathcal{O}_{\cX})}$ is defined by the morphism in the Euler sequence \eqref{eulersequence1}; see \eqref{eulersequencesurj}. Thus, the cokernel of this map is represented by the constant function in $H^0(\cX,\cO_{\cX})\cong \bC[\ud{x},\ud{p}]_0$. By reducing the term $H^0(\cX,\cO_{\cX})$ in the third row into 
$$
\tilde{H}^0(\cX,\cO_{\cX}):=\frac{H^0(\cX,\cO_{\cX})}{\bC},
$$
we obtain the exact sequence of cochain complexes:
\[
\begin{tikzcd}
0 \arrow[r] & {H^0(\cX,\mathcal{O}_{\cX})} \arrow[r, "dW"] \arrow[d, "\cong"] & {H^0(\cX,\Omega_{\cX}^{1})} \arrow[r, "dW"] \arrow[d, hook]                              & {H^0(\cX,\Omega_{\cX}^{2})} \arrow[r, "dW"] \arrow[d, hook]                              & \cdots \\
0 \arrow[r] & {\mathbb{C}[\underline{x},\underline{p}]_0} \arrow[r, "dW"]                       & {(\Omega_{\mathbb{C}[\underline{x},\underline{p}]/\mathbb{C}}^1)_0} \arrow[r, "dW"] \arrow[d, two heads] & {(\Omega_{\mathbb{C}[\underline{x},\underline{p}]/\mathbb{C}}^2)_0} \arrow[r, "dW"] \arrow[d, two heads] & \cdots \\
            &                                                                                   & {\tilde{H}^0(\cX,\mathcal{O}_{\cX})} \arrow[r, "dW"]                                   & {H^0(\cX,\Omega_{\cX}^{1})} \arrow[r, "dW"]                                              & \cdots
\end{tikzcd}\]
The cohomology $H^p(\Omega^\bullet_{\bC[\ud{x},\ud{p}]/\bC},dW)$ is completely computed in \cite[Theorem 1.6]{AS}. The only nonzero term occurs when $p=2r$, $ n+r-1$, $n+r$. In particular, $H^{2r}(\Omega^\bullet_{\bC[\ud{x},\ud{p}]/\bC},dW)$ is generated by $dW_1\wedge dp_1 \wedge \cdots \wedge dW_r \wedge dp_r$. We apply the cohomology along $dW$ and get an induced long exact sequence,
\begin{align*}
0\to \widetilde{E}^{1,0}_2\to H^1(\Omega^{\bullet}_{\bC[\ud{x},\ud{p}]/\bC},dW) \to H^0(\tilde{H}^0(\cX,\Omega^{\bullet}_{\cX}),dW) \to \widetilde{E}^{2,0}_2 \to \cdots.
\end{align*}

 
We divide our proof into three cases: $r=1$,  $n+r-1>0$, or $n+r-1=0$. When $r=1$, by \cite[Theorem 1.6]{AS}, $H^2(\Omega_{\bC[\ud{x},\ud{p}]/\bC}^\bullet,dW)$ is generated by $dW_1\wedge dp_1=\sum_{j=1}^n \frac{\partial W_1}{\partial x_j} dx_j\wedge dp_1$. By the morphism in the Euler sequence $(\Omega^2_{\bC[\ud{x},\ud{p}]})_0\twoheadrightarrow  H^0(\cX,\Omega_{\cX}^1)$ described in 
\eqref{eulersequencesurj}, this element is mapped to
\begin{align*}
\sum_{j=1}^n \frac{\partial W_1}{\partial x_j} dx_j\wedge dp_1 &\mapsto \sum_{j=1}^n \frac{\partial W_1}{\partial x_j} x_j dp_1 +\sum_{j=1}^n \frac{\partial W_1}{\partial x_j} d_1 p_1 dx_j \\
&=\sum_{j=1}^n d_1 W_1 dp_1 +\sum_{j=1}^n \frac{\partial W_1}{\partial x_j} d_1 p_1 dx_j = d_1 \cdot d(p_1W_1) = d_1 \cdot dW.
\end{align*}
Thus, we conclude that $H^k(H^0(\cX,\Omega^\bullet_{\cX}),dW)=0$ for $k=0,1,\ldots,n-1$. When $k=n$, the proof is the same as the case of $r>1$.

We now deal with the case that $n-r-1 > 0$. If we take the cohomology along $dW$, we get the following isomorphisms:
$$
\bC\cong \frac{\Ker\left(dW|_{H^0(\cX,\Omega_{\cX}^{1})}\right)}{dW\wedge \tilde{H}^0(\cX,\cO_{\cX})}\cong {H^3(H^0(\cX,\Omega_{\cX}^{\bullet}),dW)}.
$$
Note that $dW$ represents a generator of $$\frac{\Ker\left(dW|_{H^0(\cX,\Omega_{\cX}^{1})}\right)}{dW\wedge \tilde{H}^0(\cX,\cO_{\cX})}.$$ 
We let $\vartheta$ be the preimage of $dW$ along the surjection $(\Omega_{\mathbb{C}[\underline{x},\underline{p}]/\mathbb{C}}^2)_0 \twoheadrightarrow  {H^0(\cX,\Omega_{\cX}^{1})},$ described in 
\eqref{eulersequencesurj}, or explicitly,
\begin{align}
\vartheta := \sum_{j=1}^r \frac{1}{d_j}dW_j \wedge dp_j.
\label{generatorofcoh}
\end{align}
By the definition of the connecting homomorphism, the generator of $H^3(H^0(\cX,\Omega_{\cX}^{\bullet}),dW)$ is represented by $dW\wedge\vartheta\in {H^0(\cX,\Omega_{\cX}^{3})}$. Since $H^p(\Omega^{\bullet}_{\bC[\ud{x},\ud{p}]/\bC},dW)$ vanishes for $p<2r$, we have
\begin{align}
\label{oddiso}
\widetilde{E}^{3,0}_2\cong \widetilde{E}^{5,0}_2 \cong \cdots \cong \widetilde{E}^{2r-1,0}_2\cong \bC.
\end{align}
Each of them is generated by $dW\wedge\vartheta^{\wedge k}$, $k=1,\ldots, r-1$. This sequence of isomorphisms stops at $2r-1$ since $H^{2r}(\Omega^{\bullet}_{\bC[\ud{x},\ud{p}]/\bC},dW)$ is not zero, and the morphism in the Euler sequence induces an isomorphism
$$
H^{2r}(\Omega^{\bullet}_{\bC[\ud{x},\ud{p}]/\bC},dW) \to {H^0(\cX,\Omega_{\cX}^{2r-1})}
$$
which maps $dW_1\wedge dp_1 \wedge \cdots \wedge dW_r \wedge dp_r$ to $c\cdot dW\wedge\vartheta^{\wedge {(r-1)}}$ for some nonzero constant $c$. When $p=n+r-1$ or $p=n+r$, we have isomorphisms
$$
H^{n+r-1}(\Omega^{\bullet}_{\bC[\ud{x},\ud{p}]/\bC},dW)\stackrel{\cong}{\leftarrow}{H^{n+r-1}(H^0(\cX,\Omega_{\cX}^{\bullet}),dW)} \stackrel{\cong}{\rightarrow} H^{n+r}(\Omega^{\bullet}_{\bC[\ud{x},\ud{p}]/\bC},dW).
$$
They all come from the Euler sequence,
\begin{align*}
0 \to H^0(\cX,\Omega_{\cX}^{n+r-1}) \to \Omega^{n+r-1}_{\bC[\ud{x},\ud{p}]/\bC} \to H^0(\cX,\Omega_{\cX}^{n+r-2}) \to 0, \\
0=H^0(\cX,\Omega_{\cX}^{n+r}) \to \Omega^{n+r}_{\bC[\ud{x},\ud{p}]/\bC} \to H^0(\cX,\Omega_{\cX}^{n+r-1}) \to 0.
\end{align*}
This reconstructs the map $\theta$ defined in \cite[Theorem 4.4]{AS}.

Now we assume that $n-r-1=0$ so that $2r-1=n+r-2$. Then we have an isomorphism
$$
{H^{n+r-1}(H^0(\cX,\Omega_{\cX}^{\bullet}),dW)} \stackrel{\cong}{\rightarrow} H^{n+r}(\Omega^{\bullet}_{\bC[\ud{x},\ud{p}]/\bC},dW)\cong R(W)_{0}.
$$
and a short exact sequence
\[
\begin{tikzcd}
0 \arrow[r] & {\widetilde{E}^{n+r-1,0}_2} \arrow[r] & {H^{n+r-1}(\Omega^{\bullet}_{\bC[\ud{x},\ud{p}]/\bC},dW)} \arrow[r] \arrow[d, "\cong"] & {\widetilde{E}^{n+r-2,0}_2} \arrow[r] \arrow[d, "\cong"] & 0 \\
            &                           & R(W)_{0}\oplus \bC                                                                  & \bC                                         &  
\end{tikzcd}.\]
By comparing the dimensions, we still have the same isomorphisms as \eqref{oddiso}. Thus we finish the proof of Lemma \ref{ETWO}:
$$
H^p(H^q(\cX,\Omega_{\cX}^{\bullet}),dW)\cong\begin{cases}
\bC, & \textrm{if $1\leq p=q\leq n-1$ or $q=0$, $p=2k-1$, $2\leq k\leq r$,}\\
R(W)_{0}, & \textrm{if $q=0$, $p=n+r-1$}, \\
0, &\textrm{otherwise.}
\end{cases}
$$
\end{proof}

\begin{remark}\label{remark4.8}
For the first page in analytic setting $E_1^{p,q}=H^q(\cX^{an},\Omega^{an,p}_{\cX})$, we have a formula similar to \eqref{algcohpage1}:
\begin{align}
H^q(\cX^{an},(\Omega_{\cX}^{p})^{an}) \cong \begin{cases}
\bC & \textrm{if $p=q$, $1\leq p \leq n-1,$} \\
H^0(\cX^{an},(\Omega_{\cX}^{p})^{an}) & \textrm{if $q=0$, $0\leq p \leq n+r-1$,} \\
0 & \textrm{otherwise}.
\end{cases}
\end{align}
We do not have an isomorphism between $H^q(\cX,\Omega_{\cX}^{p})$ and $H^q(\cX^{an},(\Omega_{\cX}^{p})^{an})$ since $H^0(\cX^{an},(\Omega_{\cX}^{p})^{an})$ is generally larger than $H^0(\cX,\Omega_{\cX}^{p})$. Induced by the pullback of differential forms $(\Omega_{\cX}^p)^{an}$ on $\cX^{an}$ along $h:(\cX^{an},\cO_{\cX}^{an})\to (\cX,\cO_{\cX})$.  we have a map
\begin{align}
f_1^{p,q}:H^q(\cX,\Omega_{\cX}^{p}) \to H^q(\cX^{an},(\Omega_{\cX}^{p})^{an})
\label{spectralsequencemap}
\end{align}
for each $p,q$.
\end{remark}

\subsection{Algebraic setting versus analytic setting}
\label{section4.2}
Here we prove Lemma \ref{lemma5.2}. We start with the definition of a grading, which we call the weight, on $\pi_*\Omega^{p}_{\cX}$ whose graded pieces are coherent sheaves of $\cO_{\bP_{\bC}^{n-1}}$ modules, and $dW$ is homogeneous of weight $1$. From Propositions \ref{lone} and \ref{lemma5.6}, $\pi_*\Omega^{p}_{\cX}$ is a submodule of
\begin{align*}
\bigwedge^p \left(\{\pi_*\cO_{\cX}(-1)\}^{\oplus n}\oplus \bigoplus_{k=1}^r\pi_*\cO_{\cX}(d_k)\right) \cong \bigoplus_{\substack{0\leq s\leq n \\ 
 0 \leq t \leq r\\s+t=p\\ k_1,\ldots,k_r\geq 0\\1\leq j_1<\cdots<j_t\leq r}} \left(\cO_{\bP^{n-1}_{\bC}}\left(\sum_{\alpha=1}^r d_\alpha k_\alpha -s+\sum_{\alpha=1}^t d_{j_\alpha }\right)\right)^{\oplus {n \choose s}}.
\end{align*}
We define the \emph{weight} for each component as $k_1+\cdots+k_r+t$: 
\begin{align*}
&\left(\bigwedge^p \left(\{\pi_*\cO_{\cX}(-1)\}^{\oplus n}\oplus \bigoplus_{k=1}^r\pi_*\cO_{\cX}(d_k)\right)\right)_{(k)} \\&:= \bigoplus_{\substack{0\leq s\leq n \\ 
 0 \leq t \leq r\\s+t=p\\ 1\leq j_1 < \cdots < j_t \leq r\\k_1,\ldots,k_r\geq 0\\k_1+\cdots +k_r+t=k}} \left(\cO_{\bP^{n-1}_{\bC}}\left(\sum_{\alpha=1}^r d_\alpha k_\alpha -s+\sum_{\alpha=1}^t d_{j_\alpha }\right)\right)^{\oplus {n \choose s}}.
\end{align*}
Then the graded piece $(\pi_*\Omega^{q}_{\cX})_{(k)}$ of $\pi_*\Omega^{q}_{\cX}$ is defined as the intersection,
$$
\left(\pi_*\Omega^{p}_{\cX}\right)_{(k)} := \pi_*\Omega^{p}_{\cX} \cap \left(\bigwedge^p \left(\{\pi_*\cO_{\cX}(-1)\}^{\oplus n}\oplus \bigoplus_{k=1}^r\pi_*\cO_{\cX}(d_k)\right) \right)_{(k)}.
$$
This is the extension of the weight on $\bC[\ud{x},\ud{p}]$, defined in \eqref{definition:chandwt}, to differential forms on $\cX$:
\begin{align*}
\wt(x_1)=\cdots=\wt(x_n)=\wt(dx_1)=\cdots=\wt(dx_n)=0, \\
\wt(p_1)=\cdots=\wt(p_r)=\wt(dp_1)=\cdots=\wt(dp_r)=1.
\end{align*}
For a fixed $l$ and $p\geq0$, there are only finitely many choices for $k_1,\ldots,k_r\geq 0$, $0\leq t\leq r$ satisfying $k_1+\cdots+k_r+t=l$. Thus, $(\pi_*\Omega^{p}_{\cX})_{(k)}$ is a subsheaf of a locally free sheaf of finite rank. This implies that $(\pi_*\Omega^{p}_{\cX})_{(k)}$ is coherent. Also, $dW$ is a homogeneous operator since it increases the weight by $1$. Since the homogeneous components are all coherent, we can use Serre's GAGA \cite[Th\'eor\`eme 1]{SerreGAGA}. 

Since the cohomology is finite dimensional from Lemma \ref{ETWO}, it must be concentrated in the component of weight$\leq l$ for $l$ large enough, which gives the following proposition.
\begin{proposition}
There is a positive integer $l$ such that
$$
H^p(H^q(\bP^{n-1}_\bC,\pi_*\Omega_{\cX}^{\bullet}),dW) \cong H^p(H^q(\bP^{n-1}_\bC,(\pi_*\Omega_{\cX}^{\bullet})_{(\leq l)}),dW).
$$
In particular, for all $k>l$, we have
$$
H^p(H^q(\bP^{n-1}_\bC,(\pi_*\Omega_{\cX}^{\bullet})_{(k)}),dW)=0.
$$
\label{deg<l}
\end{proposition}

\begin{remark}
\label{remark4.10}
In the following commutative diagram
\[\begin{tikzcd}
	{\cX^{an}} & {\cX} \\
	{\bC\bP^{n-1}} & {\bP_{\bC}^{n-1}}
	\arrow[from=1-1, to=1-2]
	\arrow["{\pi^{an}}", from=1-1, to=2-1]
	\arrow["\pi", from=1-2, to=2-2]
	\arrow["{\iota^{an}}", shift left=3, from=2-1, to=1-1]
	\arrow[from=2-1, to=2-2]
	\arrow["\iota", shift left=3, from=2-2, to=1-2]
\end{tikzcd}\]
where $\iota$ and $\iota^{an}$ are the embeddings along the zero sections, if $\cF$ is a quasi-coherent sheaf of $\cO_{\bP_{\bC}^{n-1}}$-module, $(\iota_* \cF)^{an}$ and $\iota^{an}_* \cF^{an}$ are not the same in general. The sheaf $(\iota_* \cF)^{an}$ consists of functions that are analytic on $\cX^{an}$, but $\iota^{an}_* \cF^{an}$ consists of functions that are analytic on $\bC\bP^{n-1}$ and algebraic in the direction of the fiber. Thus, we have an injective morphism $\iota^{an}_* \cF^{an}\hookrightarrow (\iota_* \cF)^{an}$ of sheaves of $\cO_{\cX}^{an}$-modules in general. Similarly, if $\cG$ is a quasi-coherent sheaf of $\cO_{\cX}$-module, $(\pi^*\cG)^{an}$ and $\pi^{an}_* \cG^{an}$ are not the same in general. There are two different functors $(-)^{an}$: one maps $QCoh(\cO_{\cX})$ to $QCoh(\cO_{\cX}^{an})$, and the other maps $QCoh(\cO_{\bP^{n-1}_{\bC}})$ to $QCoh(\cO_{\bP^{n-1}_{\bC}}^{an})$. We will abusively denote two functors by the same notation $(-)^{an}$. In addition, we will abusively denote $\iota^{an}$ and $\pi^{an}$ by $\iota$ and $\pi$.
\end{remark}

\begin{proposition}
For all $k$,
$$
H^p(H^q(\bP^{n-1}_\bC,(\pi_*\Omega_{\cX}^{\bullet})_{(k)}),dW) \cong H^p(H^q(\bC\bP^{n-1},((\pi_*\Omega_{\cX}^{\bullet})_{(k)})^{an}),dW).
$$
In particular,
$$
H^p(H^q(\bP^{n-1}_\bC,(\pi_*\Omega_{\cX}^{\bullet}),dW) \cong H^p(H^q(\bC\bP^{n-1},(\pi_*\Omega_{\cX}^{\bullet})^{an}),dW).
$$
\label{lemma5.10}
\end{proposition}
\begin{proof}
The first isomorphism on each homogeneous component is derived from Serre's GAGA theorem \cite[Th\'eor\`eme 1]{SerreGAGA} using the fact that $(\pi_*\Omega_{\cX}^{\bullet})_{(k)}$ is a coherent sheaf on $\bP^{n-1}_{\bC}$. The second isomorphism is obvious, since $\pi_*\Omega_{\cX}^{\bullet}$ is a direct sum of homogeneous components, $(-)^{an}$ is an additive functor, and the cohomology commutes with the direct sum.
\end{proof}

\begin{proof}[Proof of Lemma \ref{lemma5.2}]
Recall that $\iota:\bC\bP^{n-1}\to \cX^{an}$ is the closed immersion along the zero section. Then $\iota$ is an affine morphism, and again by the Leray spectral sequence \cite[Lemma 01F4]{Stacks}, for each $p,q$, we have
\begin{align}
\label{cohcomparison}
H^p(H^q(\bC\bP^{n-1},(\pi_*\Omega_{\cX}^{\bullet})^{an}),dW) \cong H^p(H^q(\cX^{an},\iota_*(\pi_*\Omega_{\cX}^{\bullet})^{an}),dW).
\end{align}

Now we are given a series of morphisms between cochain complexes:
\[\begin{tikzcd}
	{(H^q(\cX^{an},\iota_*(\pi_*\Omega_{\cX}^\bullet)^{an}),dW)} \\
	{(H^q(\bC\bP^{n-1},(\pi_*\Omega_{\cX}^{\bullet}) ^{an}),dW)} \\
	{(H^q(\bP^{n-1}_{\bC},\pi_*\Omega_{\cX}^\bullet),dW)} \\
	{(H^q(\cX,\Omega_{\cX}^{\bullet}),dW)} &&& {(H^q(\cX^{an},(\Omega_{\cX}^{\bullet})^{an}),dW)}.
	\arrow["{\textrm{\eqref{cohcomparison}}}"', from=1-1, to=2-1]
	\arrow["\cong", from=1-1, to=2-1]
	\arrow["\textrm{\eqref{qisom}}",bend left=30,,from=1-1, to=4-4]
	\arrow["{\textrm{Proposition \ref{lemma5.10}}}"', from=2-1, to=3-1]
	\arrow["\cong", from=2-1, to=3-1]
	\arrow["{\textrm{Proposition \ref{leray}}}"', from=3-1, to=4-1]
	\arrow["\cong", from=3-1, to=4-1]
	\arrow["{f_1^{\bullet,q}}", from=4-1, to=4-4]
\end{tikzcd}\]
The morphism $f_1^{p,q}$ is induced by $h:(\cX^{an},\cO_{\cX}^{an})\to(\cX,\cO_{\cX})$ in Remark \ref{remark4.8}. We want to show that $f_1^{\bullet,q}$ is a quasi-isomorphism. To show this, we will verify that \begin{align}\label{qisom}
(\check{H}^q(\cU,\iota_*(\pi_*\Omega_{\cX}^{\bullet})^{an}), dW) \to (\check{H}^q(\cU,(\Omega_{\cX}^{\bullet})^{an}), dW),
\end{align}
is a quasi-isomorphism. We first consider the sheaf cohomologies as \v{C}ech cohomologies.

$\cU=\{U_1,\ldots,U_n\}$ is the open cover of $\cX$ defined in Definition \ref{volumeform}. In other words, the open set $U_j$ is biholomorphic to $\bC^{n+r-1}$ with coordinate functions being $x_1/x_j,\ldots,x_n/x_j,$ $ x_j^{d_1}p_1,\ldots,x_j^{d_r}p_r$. Their intersections $U_{i_0,\ldots,i_q}=U_{i_0}\cap \cdots\cap U_{i_q}$ are all Stein. We let $\delta$ be the \v{C}ech differential corresponding to the open cover $\cU$. Thus, by Leray's theorem and Cartan's theorem B \cite[Th\'eor\`eme B]{SerreCartan}, we have the canonical isomorphisms
$$
H^q(\cX^{an},\iota_*(\pi_*\Omega_{\cX}^{\bullet})^{an}) \cong \check{H}^q(\cU,\iota_*(\pi_*\Omega_{\cX}^{\bullet})^{an}),
$$
$$
H^q(\cX^{an},(\Omega_{\cX}^{\bullet})^{an}) \cong \check{H}^q(\cU,(\Omega_{\cX}^{\bullet})^{an}).
$$
Computing $\iota_*(\pi_*\Omega_{\cX}^{p})^{an}$ at $U_{i_0,\ldots,i_q}$, we obtain
\begin{align*}
(\iota_*(\pi_*\Omega_{\cX}^{p})^{an})(U_{i_0,\ldots,i_q}) &= (\pi_*\Omega_{\cX}^{p})^{an}(\pi(U_{i_0,\ldots,i_q})) \\
&=\pi_*\Omega_{\cX}^{p} (h\circ\pi(U_{i_0,\ldots,i_q}))\otimes_{\cO_{\bP_{\bC}^{n-1}}(h\circ\pi(U_{i_0,\ldots,i_q}))}\cO_{\bP_{\bC}^{n-1}}^{an}(\pi(U_{i_0,\ldots,i_q}))\\
&=\Omega_{\cX}^{p} (h(U_{i_0,\ldots,i_q}))\otimes_{\cO_{\bP_{\bC}^{n-1}}(h\circ\pi(U_{i_0,\ldots,i_q}))}\cO_{\bP_{\bC}^{n-1}}^{an}(\pi(U_{i_0,\ldots,i_q})).
\end{align*}
Thus, sections of $\iota_*((\pi_*\Omega_{\cX}^{p})^{an})$ are algebraic $p$-forms on $\cX$ tensored with holomorphic functions on $\bC\bP^{n-1}$. In other words, they are holomorphic along the base space $\bC\bP^{n-1}$, and algebraic along the fiber direction. We can interpret them as holomorphic forms on $\cX^{an}$ via
$$
\Omega_{\cX}^{p} (h(U_{i_0,\ldots,i_q}))\otimes_{\cO_{\bP_{\bC}^{n-1}}(h\circ\pi(U_{i_0,\ldots,i_q}))}\cO_{\bP_{\bC}^{n-1}}^{an}(\pi(U_{i_0,\ldots,i_q})) \to (\Omega_{\cX}^{p})^{an} (U_{i_0,\ldots,i_q}).
$$ 
{
Thus, by the map of cochain complexes in \eqref{qisom}, an element $[\nu]\in\check{H}^q(\cU,(\Omega_{\cX}^{p})^{an})$ is represented by a collection of sections
$$
\nu_{i_0,\ldots,i_q}\in (\Omega_{\cX}^{p})^{an}(U_{i_0,\ldots,i_q}),
$$
$$
\nu:=\prod_{i_0,\ldots,i_q}\nu_{i_0,\ldots,i_q}\in \prod_{i_0,\ldots,i_q}(\Omega_{\cX}^{p})^{an}(U_{i_0,\ldots,i_q}),
$$
such that the \v{C}ech differential $\delta\nu$ is zero.} 

We can represent each element $\nu_{i_0,\ldots,i_q}$ as a convergent power series along the fiber direction,
$$
\nu_{i_0,\ldots,i_q} = \sum_{\substack{s+t=p \\1\leq j_1<\cdots<j_r\leq r \\ k_1,\ldots,k_r\geq 0}} \eta_{s,t,\ud{j},\ud{k}}\wedge (P_1^{k_1}\cdots P_r^{k_r}  dP_{j_1}\wedge \cdots \wedge dP_{j_t}),
$$
where $\eta_{s,t,\ud{j},\ud{k}}$ is a holomorphic $s$-form on $\bC^{n-1-s}\times (\bC^*)^s$, and $P_j=x_{i_0}^{d_j}p_j$ is the affine coordinate on $U_{i_0,\ldots,i_q}$. Rearranging the power series with respect to the weight, we have
$$
\nu_{i_0,\ldots,i_q} = \sum_{k=0}^\infty \sum_{\substack{s+t=p \\1\leq j_1<\cdots<j_r\leq r\\ k_1+\cdots+k_r+t=k}} \eta_{s,t,\ud{j},\ud{k}}\wedge (P_1^{k_1}\cdots P_r^{k_r}  dP_{j_1}\wedge \cdots \wedge dP_{j_t}).
$$
We define $(\nu_{i_0,\ldots,i_q})_{(k)}$ as the weight $k$ component,
$$
(\nu_{i_0,\ldots,i_q})_{(k)} := \sum_{\substack{s+t=p \\1\leq j_1<\cdots<j_r\leq r\\k_1,\ldots,k_r\geq 0\\k_1+\cdots+k_r+t=k}} \eta_{s,t,\ud{j},\ud{k}}\wedge (P_1^{k_1}\cdots P_r^{k_r}  dP_{j_1}\wedge \cdots \wedge dP_{j_t}),
$$
{We let 
$$\nu_{(k)}\in \check{C}^q(\cU,\iota_*(\pi_*\Omega^\bullet_{\cX})_{(k)}^{an})=\prod_{i_0,\ldots,i_q}(\Omega_{\cX}^p)_{(k)}^{an}(U_{i_0,\ldots,i_q})$$
be the weight $k$ component of the $\nu$, i.e.
$$
\nu_{(k)} = \prod_{i_0,\ldots,i_q}(\nu_{i_0,\ldots,i_q})_{(k)}.
$$}
If $\nu$ is in the kernel of the \v{C}ech differential $\delta$, $\nu_{(k)}$ is also in the kernel of the \v{C}ech differential and so represents an element in the cohomology $H^q(\cX^{an},\iota_*(\pi_*\Omega^\bullet_{\cX})_{(k)}^{an}).$ 

If $dW\wedge \nu_{i_0,\ldots,i_q}=0$ for all $i_0,\ldots,i_q$, since $dW$ is also homogeneous with respect to the weight, we have $dW\wedge (\nu_{i_0,\ldots,i_q})_{(k)}=0$ for all $i_0,\ldots, i_q$. Thus, $\nu_{(k)}$ represents an element of the cohomology
$$
[\nu_{(k)}]\in H^p(H^q(\cX^{an},\iota_*(\pi_*\Omega_{\cX}^{\bullet})_{(k)}^{an}),dW).
$$
However, by Proposition \ref{deg<l}, there is a positive integer $l$ such that for $k>l$, 
$$H^p(H^q(\cX^{an},\iota_*(\pi_*\Omega_{\cX}^{\bullet})_{(k)}^{an},dW)=0.$$ 
Thus, for $k$ large enough, $\nu_{(k)}= dW\wedge \mu_{(k-1)}$ with $\mu_{(k-1)} \in \check{H}^{q}(\cU,\iota_*(\pi_*\Omega_{\cX}^{p-1})_{(k-1)}^{an})$. Hence, if $dW\wedge \nu=0$, it must be equivalent to the finite sum modulo $dW$,
\begin{align}
\nu_{i_0,\ldots,i_q} \equiv \sum_{k=0}^l \sum_{\substack{s+t=p \\1\leq j_1<\cdots<j_r\leq r\\k_1,\ldots,k_r\geq 0\\ k_1+\cdots+k_r+t=k}} \eta_{s,t,\ud{j},\ud{k}}\wedge (P_1^{k_1}\cdots P_r^{k_r}  dP_{j_1}\wedge \cdots \wedge dP_{j_t}) \pmod{dW}.
\label{qisomissurj}
\end{align}
The element on the right-hand side sits inside the image of the map \eqref{qisom}. Thus, the map \eqref{qisom} induces a surjection on cohomology. 

To show that \eqref{qisom} is injective on the cohomologies, we again observe that $dW$ is a homogeneous map. Suppose that the element described in \eqref{qisomissurj} is zero in the cohomology $H^p(H^q(\cX^{an},(\Omega^{\bullet}_{\cX})^{an}),dW)$. In other words, $\nu$ is in the image of $dW$. However, by homogeneity, each graded component $\nu_{(k)}$ must be in the image of $dW$, which implies that $\nu$ represents zero in $H^p(H^q(\cX,\Omega^{\bullet}_{\cX}),dW)$. Thus, the kernel of \eqref{qisom} on the cohomology is zero.
\end{proof}

\subsection{Higher terms in the spectral sequence}
\label{section4.3}

Our next task is to compute the higher differentials $d_k$, $k\geq 2$ in the spectral sequence $E_k^{p,q}$.

\begin{figure}[htbp]
\centering
\resizebox{\linewidth}{!}{%
\begin{tikzcd}[
    ampersand replacement=\&,
    row sep=1.5em,
    column sep=1.2em,
    execute at end picture={
        \path (\tikzcdmatrixname-8-1.east) -- coordinate (qx) (\tikzcdmatrixname-8-2.west);
        \path (\tikzcdmatrixname-8-2.south) -- coordinate (py) (\tikzcdmatrixname-9-2.north);
        \coordinate (O) at (qx |- py);
        \draw[->, thick] (O) -- ([yshift=1.5em]O |- \tikzcdmatrixname-1-1.north) node[above] {$q$};
        \draw[->, thick] (O) -- ([xshift=1.5em]O -| \tikzcdmatrixname-8-16.east) node[right] {$p$};
    }
]
n-1    \&   \&   \&   \&     \&   \&     \&        \&     \&     \&        \&         \&        \& \bC \&        \& \\
\vdots \&   \&   \&   \&     \&   \&     \&        \&     \&     \&        \& \iddots \&        \&     \&        \& \\
r      \&   \&   \&   \&     \&   \&     \&        \&     \& \bC \&        \&         \&        \&     \&        \& \\
r-1    \&   \&   \&   \&     \&   \&     \&        \& \bC \arrow[to=8-12, "d_r" description, blue, thick, shorten >=2pt, shorten <=2pt] \& \& \& \& \& \& \& \\
\vdots \&   \&   \&   \&     \& \iddots \& \&      \&     \&     \&        \&         \&        \&     \&        \& \\
2      \&   \&   \& \bC \arrow[to=8-7, "d_3" description, blue, thick, shorten >=2pt, shorten <=2pt] \& \& \& \& \& \& \& \& \& \& \& \& \\
1      \&   \& \bC \arrow[to=8-5, "d_2" description, blue, thick, shorten >=2pt, shorten <=2pt] \& \& \& \& \& \& \& \& \& \& \& \& \& \\
0      \& 0 \& 0 \& 0 \& \bC \& 0 \& \bC \& \cdots \& 0   \& 0   \& \cdots \& \bC     \& \cdots \& 0   \& \cdots \& R(W)_0 \\
       \& 0 \& 1 \& 2 \& 3   \& 4 \& 5   \& \cdots \& r-1 \& r   \& \cdots \& 2r-1    \& \cdots \& n-1 \& \cdots \& n+r-1
\end{tikzcd}%
}
\caption{The $E_2$ page of the spectral sequence. The blue arrows denote the trajectories of the non-trivial higher differentials $d_{p+1}^{p,p} : E_{p+1}^{p,p} \to E_{p+1}^{2p+1, 0}$.}
\label{fig:e2_page}
\end{figure}

\begin{lemma}
\label{higherdiffan}
The only nonzero higher differentials in the spectral sequence is
$$
{d}_{p+1}^{p,p}:{E}_{p+1}^{p,p}\cong H^p(H^p(\cX^{an},(\Omega_{\cX}^{\bullet})^{an}),dW)\to {E}^{2p+1,0}_{p+1} \cong H^{2p+1}(H^0(\cX^{an},\Omega^{an}_{\cX}),dW),
$$
for $p=1,\ldots,r-1$. Moreover, ${d}_{p+1}^{p,p}$ are all isomorphisms between vector spaces for $p=1,\ldots,r-1$.
\end{lemma}
\begin{proof}
From the description of the $E_2$-page, the only higher differentials that are possibly nonzero are $$d_{p+1}^{p,p}:E^{p,p}_{p+1}\to E^{2p+1,0}_{p+1}, \ p=1, \ldots, n-1,$$ and $$d_{\frac{n+r}{2}}^{\frac{n+r}{2}-1,\frac{n+r}{2}-1}:E_{\frac{n+r}{2}}^{\frac{n+r}{2}-1,\frac{n+r}{2}-1}\to E_\frac{n+r}{2}^{n+r-1,0}, \text{ where } n+r \text{ is even}.$$ We will show that $d_{p+1}^{p,p}$ is an isomorphism for each $p=1,\cdots,r-1$ and the zero map for each $p=r,\cdots, n-1$, and $d_{\frac{n+r}{2}}^{\frac{n+r}{2}-1,\frac{n+r}{2}-1}$ is the zero map when $n+r$ is even.

Since $d_2^{p,p},\ldots,d_p^{p,p}$ are all trivial, we have
\begin{align*}
E^{p,p}_{p+1} \cong E_1^{p,p} = H^p(\cA^{p,\bullet}(\cX^{an}),\overline{\partial}) \cong \bC,\\
E^{2p+1,0}_{p+1} \cong E_2^{2p+1,0} = H^{2p+1}(H^0(\cA^{\bullet,\bullet}(\cX^{an}),\overline{\partial}),dW) \cong \bC.
\end{align*}
In the proof of Lemma \ref{ETWO}, we showed that $H^{2p+1}(H^0(\cA^{\bullet,\bullet}(\cX^{an}),\overline{\partial}),dW)$ is generated by the $(2p+1,0)$-form $dW\wedge \vartheta^p$. We will first explicitly find a generator of $H^p(H^p(\cA^{\bullet,\bullet}(\cX^{an}),\overline{\partial}),dW)$.

Let $\omega$ be the Fubini-Study $2$-form which is the generator of $H^{1,1}(\bC\bP^{n-1};\bC)$. We let $\cV=\{V_1,\ldots,V_{n}\}$ be the open cover of $\bC\bP^{n-1}$ defined by $V_j=\{[x_1,\ldots,x_n]\in \bC\bP^{n-1}:x_j\neq 0\}$. Under the identification of the Cech cohomology and the Dolbeault cohomology on $\bC\bP^{n-1}$, we have a smooth differential form $\theta_j$ on each open subset $V_j\in \cV$  such that
\begin{align}
\theta_j|_{V_j\cap V_k} - \theta_{k}|_{V_j \cap V_k} &= \frac{dx_j}{x_j} - \frac{dx_k}{x_k},
\label{dolbeaultcech}\\
\overline{\partial}\theta_j &= c\cdot\omega,
\label{dolbeaultcech2}
\end{align}
for some nonzero constant $c$. Using the pullback of the differential along $\pi$, we have an element $[\pi^*(\omega)]\in H^1(\cA^{1,\bullet}(\cX^{an}))$. To be specific, we define $\theta_j$ as follows:
$$
\theta_j := \sum_{\substack{1\leq k\leq n\\k\neq j}}\frac{\overline{X}_{k}dX_{k}}{1+\sum_{\substack{1\leq l\leq n\\l\neq j}}|X_{l}|^2},
$$
where $X_k=x_k/x_j$ is the affine coordinate of $V_j$. Since $\pi^*(\omega)$ has the same local description as $\omega$, we will use abusively the notation $\omega$ for $\pi^*(\omega)$. 

We now show that $[\omega^{\wedge p}]$ is a generator of $H^p(\cA^{p,\bullet}(\cX^{an}),\overline{\partial})$. Using the closed immersion $\iota:\bC\bP^{n-1}\hookrightarrow \cX$ defined in Remark \ref{remark4.10}, we show that $\pi^*$ yields an injective map $H^p(\cA^{p,\bullet}(\bC\bP^{n-1}),\overline{\partial})\to H^p(\cA^{p,\bullet}(\cX^{an}),\overline{\partial})$: since $\pi\circ \iota = \mathrm{id}_{\bC\bP^{n-1}}$, $\iota^* \circ \pi^*$ yields the identity map on the cohomology $H^p(\cA^{p,\bullet}(\bC\bP^{n-1}),\overline{\partial})$ and so $\pi^*:H^p(\cA^{p,\bullet}(\bC\bP^{n-1}),\overline{\partial})\to H^p(\cA^{p,\bullet}(\cX^{an}),\overline{\partial})$ is injective. We will show that $[\omega^{\wedge p}]$ is mapped to $c_p [dW\wedge \vartheta^{\wedge p}]$ with nonzero constant $c_p$ by $d_{p+1}^{p,p}:E_{p+1}^{p,p}\to E_{p+1}^{2p+1,0}$ for $p=1,\ldots,r-1$.

Dolbeault's isomorphism $H^1(\cA^{1,\bullet}(\cX^{an}),\overline{\partial})\cong H^{1}(\cX^{an},\Omega^{an,1}_{\cX})$ maps $c\cdot\omega$ to
$$
\left[\left(\frac{dx_j}{x_j} - \frac{dx_k}{x_k}\right)_{(j,k)}\right]\in H^{1}(\cX^{an},(\Omega_{\cX}^1)^{an})\cong H^{1}(\cX,\Omega_{\cX}^{1}).
$$
On $U_j$ (which is defined in Definition \ref{volumeform}), we consider the $(2,0)$-form $\alpha$ whose restriction on $U_j$ is
$$
\alpha|_{U_j}=\sum_{l=1}^r \frac{1}{d_l} d\left(\frac{W_l}{x_j^{d_l}}\right)\wedge d(x_j^{d_l}p_l) - dW \wedge \theta_j.
$$
A crucial equality for the proof is the following: for $j,k\in \{1,\ldots,n\}$,
\begin{align}
\sum_{l=1}^r\frac{1}{d_l}\left(d\left(\frac{W_l}{x_j^{d_l}}\right)\wedge d(x_j^{d_l}p_l)-d\left(\frac{W_l}{x_k^{d_l}}\right)\wedge d(x_k^{d_l}p_l)\right) = dW\wedge \left(\frac{dx_j}{x_j}-\frac{dx_k}{x_k}\right),
\label{higherdiffeq}
\end{align}
which is derived from the following elementary equation: for each $j=1,\ldots,n$,
\begin{align}
\sum_{l=1}^r \frac{1}{d_l} d\left(\frac{W_l}{x_j^{d_l}}\right)\wedge d(x_j^{d_l}p_l) = \sum_{l=1}^r \frac{1}{d_l} dW_l \wedge dp_l + dW \wedge \frac{dx_j}{x_j}.
\label{higherdiffeq2}
\end{align}
An element $(\frac{dx_j}{x_j}-\frac{dx_k}{x_k})_{j,k}\in \check{C}^1(\cU,\Omega_{\cX}^{1})$ represents a class in $H^1(\cX,\Omega_{\cX}^{1})\cong H^1(\cX^{an},(\Omega_{\cX}^{1})^{an})$. We have an element $\left(\sum_{l=1}^r\frac{1}{d_l}{d\left(\frac{W_l}{x_j^{d_l}}\right)\wedge d(x_j^{d_l}p_l)}\right)_{j}$ of $\check{C}^0(\cU,\Omega_{\cX}^{2})$, which is not in the kernel of $\delta$. However, after taking $dW$ on \eqref{higherdiffeq2}, we get an element of $E_2^{3,0}=H^3(H^0(\cX^{an},(\Omega^{\bullet}_{\cX})^{an}),dW)$:
\begin{align}
\label{higherdiffeq3}
dW\wedge \left(\sum_{l=1}^r\frac{1}{d_l}{d\left(\frac{W_l}{x_j^{d_l}}\right)\wedge d(x_j^{d_l}p_l)}\right) = dW\wedge \vartheta,
\end{align}
where $\vartheta$ is defined in \eqref{generatorofcoh}. 

$\alpha_j$ defines a global $(2,0)$-form since by \eqref{higherdiffeq} and \eqref{dolbeaultcech},
$$
\alpha|_{U_j \cap U_k} - \alpha|_{U_k \cap U_j} = 0.
$$
Moreover, by \eqref{higherdiffeq3} and \eqref{dolbeaultcech2}, we have
\begin{align*}
dW\wedge \alpha = dW\wedge \vartheta, \\
\overline{\partial}\alpha = c\cdot dW\wedge \omega.
\end{align*}
We have shown that $c$ is a nonzero constant in \eqref{dolbeaultcech2}. Thus, along the staircase
\[\begin{tikzcd}
	{\cA^{1,1}(\cX^{an})} & {\cA^{2,1}(\cX^{an})} \\
	& {\cA^{2,0}(\cX^{an})} & {\cA^{3,0}(\cX^{an}),}
	\arrow["dW", from=1-1, to=1-2]
	\arrow["{\overline{\partial}}"', from=2-2, to=1-2]
	\arrow["dW", from=2-2, to=2-3]
\end{tikzcd}\]
the Fubini-Study $2$-form $\omega$ is transferred as follows:
\begin{equation}
\label{cdrelation}
\begin{tikzcd}
	{c\cdot\omega} & {c\cdot dW\wedge \omega} \\
	& \alpha & {dW\wedge \vartheta.}
	\arrow["dW", maps to, from=1-1, to=1-2]
	\arrow["{\overline{\partial}}"', maps to, from=2-2, to=1-2]
	\arrow["dW", maps to, from=2-2, to=2-3]
\end{tikzcd}
\end{equation}

When $2\leq p \leq r-1$, we apply the relation \eqref{cdrelation} consecutively to obtain $d_{p+1}^{p,p}([\omega^{\wedge p}])$. For convenience, we write $\omega_0 = c\omega$. Then, \eqref{cdrelation} is applied to $\omega_0^{\wedge p}$ as follows:
\[\begin{tikzcd}
	{\omega_0^{\wedge p}} & {dW\wedge \omega_0^{\wedge p}} \\
	& {\alpha\wedge \omega_0^{\wedge (p-1)}} & {\alpha \wedge dW\wedge \omega_0^{\wedge (p-1)}} \\
	&& {\frac{1}{2!}\alpha^{\wedge 2}\wedge \omega_0^{\wedge (p-1)}} & \ddots \\
	&&& {\frac{1}{p!}\alpha^{\wedge p}} & {\frac{1}{p!}dW\wedge \vartheta^{\wedge p}}
	\arrow["dW", maps to, from=1-1, to=1-2]
	\arrow["{\overline{\partial}}"', maps to, from=2-2, to=1-2]
	\arrow["dW", maps to, from=2-2, to=2-3]
	\arrow["{\overline{\partial}}"', maps to, from=3-3, to=2-3]
	\arrow["dW", maps to, from=3-3, to=3-4]
	\arrow["{\overline{\partial}}"', maps to, from=4-4, to=3-4]
	\arrow["dW", maps to, from=4-4, to=4-5]
\end{tikzcd}.\]
Thus, $d_{p+1}^{p,p}$ maps the generator $[\omega]$ of $E_{p+1}^{p,p}$ to a nonzero element in $E_{p+1}^{2p+1,0}$.

We automatically have ${d}_{p+1}^{p,p}=0$ for $p\geq r$ since $dW\wedge \vartheta^{\wedge p}=0$ for $p\geq r$. Thus, when $n+r$ is even, $(n+r)/2>r$ and $d_{\frac{n+r}{2}}^{\frac{n+r}{2}-1,\frac{n+r}{2}-1}$ is a zero map.
\end{proof}

\begin{remark}
\begin{enumerate}
\item
According to the proof of Lemma \ref{higherdiffan}, it is possible to find the explicit representative of $H^{2k}(\cA^\bullet(\cX^{an}),\overline{\partial}+dW)$ for each $r\leq k \leq n-1$:
$$
\sum_{j=0}^r \frac{(-1)^j}{j!}\alpha^{\wedge j}\wedge \omega_0^{\wedge (k-j)}
$$
is a representative for a generator of $H^{2k}(\cA^\bullet(\cX^{an}),\overline{\partial}+dW)$. One can show that it is in the kernel of $\overline{\partial}+dW$ via the relations
\begin{align*}
dW\wedge \alpha^{\wedge (j-1)} \wedge \omega_0^{\wedge
(k-j+1)} &= \frac{1}{j}\overline{\partial}(\alpha^{\wedge j}\wedge \omega_0^{\wedge(k-j)}), \qquad (j=1,\ldots,r-1),\\
\overline{\partial}\omega_0^{\wedge k}&=0,\\
dW\wedge \alpha^r \wedge \omega_0^{\wedge(k-r)}&=dW\wedge \vartheta^r\wedge \omega_0^{\wedge(k-r)}=0.
\end{align*}
\item
The equations \eqref{higherdiffeq}, \eqref{higherdiffeq2}, and \eqref{higherdiffeq3} can also be used to prove the same fact for $\widetilde{E}_k^{p,q}$, that is, $\widetilde{d}_{p+1}^{p,p}$ are all isomorphisms between vector spaces for $p=1,\ldots,r-1$, and all other $\widetilde{d}_k^{p,q}$'s are zero. This implies that there is a map $f_k^{p,q}:\widetilde{E}_k^{p,q}\to E_k^{p,q}$ where $f_1$ is defined in \eqref{spectralsequencemap}, $f_2$ is defined in Lemma \ref{lemma5.2}, and $f_k$ is induced from $f_{k-1}$ for each $k=3,4,\ldots$. Since $f_2$ is an isomorphism, by the mapping lemma \cite[Lemma 5.2.4]{Weibel}, $f_3,f_4,\ldots$ are all isomorphisms. 
Since both $E_k^{p,q}$ and $\widetilde{E}_{k}^{p,q}$ are convergent spectral sequences, both spectral sequences $E_k^{p,q}$ and $\widetilde{E}_{k}^{p,q}$ compute $H^{p+q}(PV(\cX),\overline{\partial}_W)$.
{
\item
From the commutative diagram 
\eqref{eulersequencedW}, we obtain the following isomorphism:
\begin{align}
\widetilde{E}_1^{1,1} \dashrightarrow H^{0}(\cX,\cO_{\cX}) \to \frac{H^{0}(\cX,\Omega_{\cX}^{1})}{dW\wedge \tilde{H}^{0}(\cX,\cO_{\cX})} \to \widetilde{E}_2^{3,0}.
\label{d2new}
\end{align}
In the proof of Lemma \ref{higherdiffan}, we observe that this is the same as $\widetilde{d}_2^{1,1}$. In the proof of Lemma \ref{ETWO}, two isomorphisms \eqref{diagiso} and \eqref{oddiso} on the spectral sequence $\widetilde{E}_2^{p,q}$ are described. Taking the composition of \eqref{diagiso}, \eqref{d2new}, and \eqref{oddiso}, for $1\leq p \leq r-1$, we have
$$
\widetilde{E}_1^{p,p}\cong \cdots \cong \widetilde{E}_1^{2,2} \cong \widetilde{E}_1^{1,1} \stackrel{d_2^{1,1}}{\to} \widetilde{E}_2^{3,0} \cong \widetilde{E}_2^{5,0} \cong \cdots \cong \widetilde{E}^{2p+1,0}_2.
$$
Since $\widetilde{d}_k^{p,p}$ and $\widetilde{d}_k^{2p+1,0}$ are zero for $k=2,\ldots,p$, we have 
\begin{align}
\label{dpnew}
\widetilde{E}_{p+1}^{p,p} \cong \widetilde{E}_2^{p,p} \cong \widetilde{E}^{2p+1,0}_2 \cong \widetilde{E}^{2p+1,0}_{p+1}.
\end{align}
In the proof of Lemma \ref{higherdiffan}, we observed that the higher differentials $\widetilde{d}_{p+1}^{p,p}$ are actually the composition \eqref{dpnew}.
}
\end{enumerate}
\end{remark}

\begin{comment}
\begin{proof}[Proof of Theorem \ref{cohomologycmp}]
Recall we have $f_1^{p,q}:\widetilde{E}_1^{p,q}\to E_1^{p,q}$ and $f_2^{p,q}:\widetilde{E}_2^{p,q}\to E_2^{p,q}$ in \eqref{spectralsequencemap} and Lemma \ref{lemma5.2}. By Lemmata \ref{higherdiffalg} and \ref{higherdiffan}, we have a map $f:\widetilde{E}\to E$ between spectral sequences: for $k=1,2,\cdots$ and each $p,q$
\begin{align*}
f_k^{p,q}: \widetilde{E}_{k}^{p,q}\to {E}_{rk}^{p,q},\\
d_k^{p,q} \circ f_k^{p,q} = f_{k}^{p+k,q-k+1}\circ d_{k}^{p,q}.
\end{align*}
(See \cite[p. 123]{Weibel}.) Using the comparison theorem of the spectral sequences \cite[Theorem 5.2.12]{Weibel}, we conclude that $E_k^{p,q}$ and $\widetilde{E}_k^{p,q}$ converge to the same cohomology $H^{p+q}(PV(\cX),\overline{\partial}_W).$
\end{proof}
\end{comment}

We show that $R(W)_{0}$ can be embedded not only as a subspace, but also as a subring inside $H(PV^\bullet(\cX),\overline{\partial}_W)$.

\begin{proposition}
There is an injective ring homomorphism
$$
R(W)_{0} \hookrightarrow H^{0}(PV^\bullet(\cX),\overline{\partial}_W).
$$
\end{proposition}
\begin{proof}
Since $\bC[\ud{x},\ud{p}]_{0}$ is contained in the kernel of $\overline{\partial}_W$ in $PV^{0,0}(\cX)$, an element in the image must represent a cohomology class in $H^{0}(PV^\bullet(\cX),\overline{\partial}_W)$. Suppose $v\in PV^{0,1}(\cX)$ is a holomorphic vector field, i.e. $\overline{\partial}(v)=0$. Then, $\{W,v\}\in \bC[\ud{x},\ud{p}]_{0}$ is contained in the Jacobian ideal of $W$. Thus, the map $\bC[\ud{x},\ud{p}]_{0}\to H^0(PV^{\bullet}(\cX),\overline{\partial}_W)$ factors through the quotient by the Jacobian ideal
$$
\bC[\ud{x},\ud{p}]_{0}\twoheadrightarrow R(W)_{0} \hookrightarrow H^0(PV^{\bullet}(\cX),\overline{\partial}_W).
$$
Since $PV^{0,0}(\cX)$ is closed under the multiplication, $H^0(PV^\bullet(\cX),\overline{\partial}_W)$ is also closed under the multiplication. Since $\bC[\ud{x},\ud{p}]\to PV^{0,0}(\cX)$ is a ring homomorphism, we get the conclusion.
\end{proof}

\begin{remark}\label{compare}
There is another Frobenius manifold structure on $H(V(\ud{W});\bC)$ due to Barannikov-Kontsevich \cite{BK}. When $r=1$, we showed \cite{Paper1} that the isomorphism in Theorem \ref{Main2} is a ring isomorphism; in particular, $R(W)_{0}\cong H_{\mathrm{prim}}(V(\ud{W});\bC)$ is a ring isomorphism, where the ring structure on $H_{\mathrm{prim}}(V(\ud{W});\bC)$ is given by Barannikov-Kontsevich. When $r>1$, we do not currently know whether \eqref{main:iso} is a ring isomorphism; the method of the proof in \cite{Paper1} relies on a computation of the residue map for a hypersurface due to Carlson-Griffiths \cite{CG} but an analogous computation is not available for a higher codimensional case ($r>1$ case).  
\end{remark}

\section{Appendix}\label{section5}

\subsection{Positivity of the injective radius of an asymptotically conical manifold}
\label{appendix2}

In this section, we give a technical detail skipped in the proof of Proposition \ref{ACisbg}.

\begin{proposition}
Let $(M,g)$ be an asymptotically conical manifold together with the Riemannian cone $(C,g_C)$ and the diffeomorphism $\Phi:C\setminus K'\to M\setminus K$ such that $K'$ and $K$ are compact, and for some $\a>0$
\[
|(\nabla^C)^k (\Phi^* g - g_C)|_{g_C} = O(\rho_C^{-\a-k})
\]
for each $k=0,1,\ldots.$ Then the injectivity radius of $(M,g)$ is positive.
\end{proposition}
\begin{proof}
From the standard fact on Riemannian cones, there is a compact Riemannian manifold $(L,g_L)$ such that there is an isometry $(C,g_C)\cong(\bR_{>0}\times L,\rho^2g_L+d\rho^2)$, where $\rho$ is a standard coordinate on $\bR_{>0}$. Hence we identify $(C,g_C)$ to $(\bR_{>0}\times L,\rho^2g_L+d\rho^2)$. Via this identification, $\rho$ is identified to the radial coordinate $\rho_C$ on $C$.

We take $\rho_0>0$ large enough so that $(0,\rho_0)\times L\subset K'$, where $K'$ is a compact subset in $C$ given by the fact that $(M,g)$ is asymptotically conical. For each $\rho>\rho_0$, we consider a one parameter family of metrics $g_\rho$ on \[P:=[1/2,2]\times L\] given by
\[
g_\rho := (\rho^{-2}\delta_\rho^*\Phi^* g)|_P
\]
We have $|(\nabla^{C})^k (\Phi^* g-g_C)|_{g_C}\to 0$ as $\rho \to \infty$ for each $k$, so $g_\rho\to g_C$ on $P$ in $C^\infty$-topology. To avoid confusion, we write $t$ for the coordinate function on $[\tfrac12,2]$.

However, $P$ itself is not a closed manifold, so we cannot directly apply {\cite[Chapter~6, Theorem~1]{Ehr74}} to $P$. We use the ``doubling trick" to construct a closed manifold $\widehat{P}$ together with new metrics $\widehat{g}_C$ and $\widehat{g}_{\rho}$. Let $\overline{P}$ be the copy of $P$ with flipped orientation. As a topological space, $\widehat{P}=P\cup_{\partial P} \overline{P}$ is given by gluing $P$ and $\overline{P}$ along the boundary, and it is homeomorphic to $S^1\times P$. From the standard fact of differential topology, $\widehat{P}$ can be given a smooth structure so that it is diffeomorphic in the compact subset $P'\subset P$ such that $P'\cap\partial P=\emptyset$; say $P':=[\tfrac34,\tfrac43]\times L$. Now we move onto to the problem of defining the Riemannian metric $\widehat{g}_C$ on $\widehat{P}$. On $[\tfrac12,2]$, we consider three smooth cutoff functions $\chi_1,\chi_2,\chi_3$ such that $\chi_1+\chi_2+\chi_3=1$ and
\begin{align*}
\chi_1&\equiv 1\quad\textrm{on $[\tfrac12,\tfrac45]$},& \quad \chi_1&\equiv 0\quad\textrm{on $[\tfrac56,2]$}, \\
\chi_2&\equiv 1\quad\textrm{on $[\tfrac45,\tfrac54]$},& \quad \chi_2&\equiv 0\quad\textrm{on $[\tfrac12,\tfrac56]\cup[\tfrac65,2]$}, \\
\chi_3&\equiv 1\quad\textrm{on $[\tfrac65,2]$},& \quad \chi_3&\equiv 0\quad\textrm{on $[\tfrac12,\tfrac54]$}.
\end{align*}
Note that $\chi_2$ is $1$ on the interval larger than $[\tfrac34,\tfrac43]$, because we need to control the injectivity radius (which is a value determined by the metric on the neighborhood of a point), and we do not want $P'$ to touch the boundary of the space where the metric $\widehat{g}_C$ differs from $g_C$. Now we define $\widehat{g}_C$ on $P$ as follows:
\[
\widehat{g}_C := dt^2 + (\chi_1(t)(\tfrac56)^2+\chi_2(t)t^2+\chi_3(t)(\tfrac65)^2)g_L.
\]
We write 
\[P'' := [\tfrac45,\tfrac54]\times L, \qquad \partial P'' = \{\tfrac45,\tfrac54\}\times L.\] Hence, $\widehat{g}_C=g_C$ on $P'$, and on $P\setminus P''$, $\widehat{g}_C$ is just a cylindrical metric. Thus, by defining $\widehat{g}_C$ on $\overline{P}$, we can smoothly attach $\widehat{g}_C$ onto the whole space $\widehat{P}$. Finally, we define $\widehat{g}_{\rho}$ on $P$ as
\[
\widehat{g}_\rho = \chi_2(t) g_{\rho} +  (\chi_1(t)+\chi_3(t))dt^2+(\chi_1(t)(\tfrac56)^2+\chi_3(t)(\tfrac65)^2)g_L.
\]
One can check that this is positive definite. Then $\widehat{g}_{\rho}$ is cylindrical near the boundary, so it can be extended to $\widehat{P}$. Finally, we can easily see that in $\widehat{P}$, $\widehat{g}_{\rho}\to \widehat{g}_C$ in the $C^\infty$ topology since $\widehat{g}_{\rho}-\widehat{g}_C= \chi_2(\rho)(g_\rho-g_C)$. 

We use {\cite[Chapter~6, Theorem~1]{Ehr74}} on $\widehat{P}$: for each $(t,x)\in P'$,
\[
\inj((t,x),\widehat{g}_C) \leq \liminf_{\rho\to \infty}\inj((t,x),\widehat{g}_\rho).
\]
Note that inside $P'$, $\widehat{g}_C=g_C$ and $\widehat{g}_\rho=g_\rho$. However, we must consider that the injectivity radius (even the pointwise version) contains global information; hence, we do not want our geodesic ball in $P'$ to touch $t=\tfrac45$ or $t=\tfrac54$. 
 Set
\[
d_0:=\min_{(s,y)\in P'}\dist_{\widehat g_C}((s,y),\partial P'')>0.
\]
Hence, inside $P'$, we define $m$ as
\[
m = \min_{(s,y)\in P'} (\inj((s,y),\widehat{g}_C),d_0)>0.
\]
Then,
\[
m\leq \liminf_{\rho\to \infty}\inj((t,x),\widehat{g}_\rho),
\]
and we have $c_1>0$ such that when $\rho>c_1$,
\[
\inj((t,x),\widehat{g}_\rho) \geq \frac{m}{2}.
\]

Since $\widehat g_\rho\to \widehat g_C$ in $C^0(\widehat P)$ and $\widehat P$ is compact, for $\rho\gg1$
we have $(1-\varepsilon_\rho)\widehat g_C(X,X)\le \widehat g_\rho(X,X)\le (1+\varepsilon_\rho)\widehat g_C(X,X)$ with
$\varepsilon_\rho=\sup_{\widehat{P}} \|\widehat{g}_C^{-1}(\widehat{g}_\rho-\widehat{g})\|_{\widehat{g}_C}\stackrel{\rho\to\infty}{\longrightarrow}0$. Hence,
\[
\min_{(s,y)\in P'}\dist_{\widehat g_\rho}((s,y),\partial P'')\ \longrightarrow\ d_0.
\]
In particular, there exists $c_2$ such that for $\rho\ge c_2$,
\[
\min_{(s,y)\in P'}\dist_{\widehat g_\rho}((s,y),\partial P'')\ \ge\ \frac{d_0}{2}\ \ge\ \frac{m}{2},
\]
since $m\le d_0$ by definition. Take $c=\max(c_1,c_2)$. Then, for all $\rho\ge c$, we have
\[
\inj((t,x),g_\rho)\geq \min\left(\inj((t,x),\widehat{g}_\rho),\min_{(s,y)\in P'}\dist_{\widehat{g}_{\rho}}((s,y),\partial P'')\right)\geq \frac{m}{2}.
\]
Thus, we put $g_\rho=\rho^{-2}\delta_\rho^*\Phi^*g$ and
\[
\inj(\delta_\rho(t,x),\Phi^* g)=\rho\;\inj((t,x),\rho^{-2}\delta_\rho^*\Phi^* g)\geq \frac{m}{2}\rho.
\]
\end{proof}

\subsection{Simultaneous weighted blow-up of $X_{CY}$ and $X_{LG}$} 
\label{section5.2}
In this section, we return to the case where $X_{CY}$ can be singular, i.e. $w_1,\ldots,w_n$ are no longer equal to $1$'s. Furthermore, we assume that $\gcd(d_1,\ldots,d_r)=\gcd(w_1,\ldots,w_n)=1$.
{
Extending Li-Wen's theory to orbifolds using the simultaneous weighted blow-up in this appendix is an anticipated future direction for general weighted projective CY complete intersections.
}

There is a stack theoretic definition of the weighted blow-up. However, we will deal with the simplicial toric varieties for simplicity of presentation: we will work on the coarse moduli space of stacks, i.e. classical algebraic geometry. Let us give a classical definition of the weighted blow-up for a simplicial toric variety.

\begin{definition}
Let $\sigma$ be a cone and $\rho$ be a ray that is not a face of $\sigma$. A \emph{weighted star subdivision} $\Sigma_{\sigma,\rho}$ of a cone $\sigma$ along the ray $\rho$ is defined as a fan of cones generated by $\rho$ and a proper face of $\sigma$ not containing $\rho$,
$$
\Sigma_{\sigma,\rho} = \bigcup_{\substack{\tau\lneq  \sigma\\ \rho^\circ\nsubseteq\tau^\circ}}(\textrm{Faces of }\tau + \rho)
$$
where $\tau^\circ$ is the interior of $\tau$.

A fan $\widetilde{\Sigma}$ is called the \emph{weighted star subdivision} of $\Sigma$ along the ray $\rho\in \widetilde{\Sigma}$ if 
$$
\widetilde{\Sigma} = \{\sigma\in\Sigma:\rho \nsubseteq\sigma\} \cup\left( \bigcup_{\substack{\rho \subseteq \sigma\\\rho^{\circ}\nsubseteq\sigma^\circ}} \Sigma_{\sigma,\rho}\right).
$$
\end{definition}

\begin{definition}[Weighted blow-up] Let $X$ be a simplicial toric variety defined by a fan $\Sigma$. A birational morphism $\phi:\widetilde{X}\to X$ is called the \emph{weighted blow-up} if the fan $\widetilde{\Sigma}$ of $\widetilde{X}$ is the weighted star subdivision of $\Sigma$ along the ray $\rho$ contained in the cone $\sigma$ of $\Sigma$, and $\phi$ is represented by the natural surjection $\widetilde{\Sigma}\to\Sigma$.
\end{definition}

See \cite[Section 2.2]{Kerr} for the equivalence between the stack theoretic definition and the classical definition on toric varieties.

\begin{definition}
Let $\widetilde{X}$ be a toric variety defined by the group quotient of $\bC\times (\bC^n\setminus\{0\})\times (\bC^r \setminus \{0\})$ by the group action of $(\bC^*)^2$, where $(\lambda,\mu)\in (\bC^*)^2$ acts on $(z,x_1,\ldots,x_n,p_1,\ldots,p_r)$ as
\begin{align}
(\lambda,\mu)\cdot (z,x_1,\ldots,x_n,p_1,\ldots,p_r) = ((\lambda\mu)^{-1}z,\lambda^{w_1}x_1,\ldots,\lambda^{w_n}x_n,\mu^{d_1}p_1,\ldots,\mu^{d_r}p_r),
\label{groupaction}
\end{align}
\begin{align*}
\widetilde{X} = \frac{\bC\times (\bC^n\setminus\{0\})\times (\bC^r \setminus \{0\})}{(z,x_1,\ldots,x_n,p_1,\ldots,p_r)\sim (\lambda,\mu)\cdot(z,x_1,\ldots,x_n,p_1,\ldots,p_r)}.
\end{align*}
\end{definition}

We define morphisms
$\phi_{CY}:\widetilde{X}\to X_{CY}$ and 
$\phi_{LG}:\widetilde{X}\to X_{LG}$ by
\begin{align}
\begin{array}{c}
\phi_{CY}[z,x_1,\ldots,x_n,p_1,\ldots,p_r] = [x_1,\ldots,x_n,z^{d_1}p_1,\ldots,z^{d_r}p_r],  \\
\phi_{LG}[z,x_1,\ldots,x_n,p_1,\ldots,p_r] = [z^{w_1}x_1,\ldots,z^{w_n}x_n,p_1,\ldots,p_r].
\end{array}\label{phiCY}
\end{align}
To show that this is well-defined, we need to show that the image of the group orbit of $(\bC^*)^2$ in $\bC\times (\bC^n\setminus\{0\})\times (\bC^r\setminus\{0\})$ is the orbit of the charge action:
\[\begin{tikzcd}
	{[z,x_1,\ldots,x_n,p_1,\ldots,p_r]} & {[x_1,\ldots,x_n,z^{d_1}p_1,\ldots,z^{d_r}p_r]} \\
	{[(\lambda\mu)^{-1}z,\lambda^{w_1}x_1,\ldots,\lambda^{w_n}x_n,\mu^{d_1}p_1,\ldots,\mu^{d_r}p_r]} & {[\lambda^{w_1}x_1,\ldots,\lambda^{w_n}x_n,\lambda^{-d_1}z^{d_1}p_1,\ldots,\lambda^{-d_r}z^{d_r}p_r]}
	\arrow["{\phi_{CY}}", mapsto, from=1-1, to=1-2]
	\arrow["{(\lambda,\mu)\cdot}"', no head, from=1-1, to=2-1]
	\arrow[shift left, no head, from=1-1, to=2-1]
	\arrow["{\lambda\cdot}", no head, from=1-2, to=2-2]
	\arrow[shift right, no head, from=1-2, to=2-2]
	\arrow["{\phi_{CY}}",mapsto, from=2-1, to=2-2]
\end{tikzcd},\]
\[\begin{tikzcd}
	{[z,x_1,\ldots,x_n,p_1,\ldots,p_r]} & {[z^{w_1}x_1,\ldots,z^{w_n}x_n,p_1,\ldots,p_r]} \\
	{[(\lambda\mu)^{-1} z,\lambda^{w_1}x_1,\ldots,\lambda^{w_n}x_n,\mu^{d_1}p_1,\ldots,\mu^{d_r}p_r]} & {[\mu^{-w_1}z^{w_1}x_1,\ldots,\mu^{-w_n}z^{w_n}x_n,\mu^{d_1}p_1,\ldots,\mu^{d_r}p_r]}
	\arrow["{\phi_{LG}}", mapsto, from=1-1, to=1-2]
	\arrow["{(\lambda,\mu)\cdot}"', no head, from=1-1, to=2-1]
	\arrow[shift left, no head, from=1-1, to=2-1]
	\arrow["{\mu^{-1}\cdot}", no head, from=1-2, to=2-2]
	\arrow[shift right, no head, from=1-2, to=2-2]
	\arrow["{\phi_{LG}}", mapsto, from=2-1, to=2-2]
\end{tikzcd}.\]

The algebraic variety $\widetilde{X}$ is a line bundle $\cO_{\bP(\ud{w})\times \bP(\ud{d})}(-1,-1)$ over $\bP(\ud{w})\times \bP(\ud{d})$. We have an embedding $\bP(\ud{w})\times \bP(\ud{d}) \hookrightarrow \widetilde{X}$ by the zero section. When $z\neq 0$, $\phi_{CY}|_{\{z\neq0\}}$ is injective onto $X_{CY}\setminus \bP(\ud{w})$ since
\begin{align*}
\phi_{CY}[1,x_1,\ldots,x_n,p_1,\ldots,p_r]=[x_1,\ldots,x_n,p_1,\ldots,p_r].
\end{align*}
Similarly, when $z\neq 0$, $\phi_{LG}|_{\{z\neq0\}}$ is injective onto $X_{LG}\setminus \bP(\ud{d})$. Thus, $\phi_{CY}$ and $\phi_{LG}$ are birational morphisms.

\begin{theorem}
\label{blowup}
The morphisms $\phi_{CY}$ and $\phi_{LG}$ satisfy the following properties.
\begin{enumerate}
\item 
$\phi_{CY}:\widetilde{X}\to X_{CY}$ and 
$\phi_{LG}:\widetilde{X}\to X_{LG}$ are weighted blow-ups.
\item 
The exceptional divisor of $\phi_{CY}$ (respectively, $\phi_{LG}$) is given by the inverse image of zero section $\bP(\ud{w})\hookrightarrow X_{CY}$ (respectively, $\bP(\ud{d})\hookrightarrow X_{LG}$):
\label{exdiv}
\begin{align*}
\phi_{CY}^{-1}(\bP(\ud{w})) = \bP(\ud{w})\times \bP(\ud{d})=\phi_{LG}^{-1}(\bP(\ud{d})) .
\end{align*}
\end{enumerate}
\end{theorem}
\begin{remark}
When $r=1$ and $d=1$, $\bP(\ud{d})$ is a singleton, so $\phi^{-1}_{CY}(\bP(\ud{w})) = \bP(\ud{w})$, and $\phi_{CY}$ is an isomorphism. Thus, a morphism
\begin{align}X_{CY}\cong \widetilde{X}\to X_{LG}\label{r=1blowup}\end{align}
is the blow-up morphism. 
\end{remark}

In summary, we have the following commutative diagram:
$$
\begin{tikzcd}
& \bP(\ud{w})\times \bP(\ud{d}) \ar[rd,"p_2"] \ar[ld,"p_1"'] \ar[dd,"\textrm{ex. div.}",hook]\\ 
\bP(\ud{w}) \ar[dd,"\textrm{zero}"',hook] & & \bP(\ud{d}) \ar[dd,"\textrm{zero}",hook] \\
& \widetilde{X} \ar[rd,->>,"\phi_{LG}"] \ar[ld,->>,"\phi_{CY}"'] \\
X_{CY} & & X_{LG}\\
& \frac{(\bC^n\setminus\{0\})\times(\bC^r\setminus\{0\})}{\bC^*} \ar[ur,hook] \ar[uu,hook] \ar[lu,hook]
\end{tikzcd}.
$$

For the proof of Theorem \ref{blowup}, we need to write down the fans explicitly.
We fix the fans of $\bP(\ud{w})$ and $\bP(\ud{d})$. Let $N_x$ be the lattice of one parameter subgroup of $\bP(\ud{w})$, and $M_x$ be the lattice of characters of $\bP(\ud{w})$. Let $\rho_{x_1},\ldots,\rho_{x_n}\subset (N_x)_{\bR}$ be the rays of the fan of $\bP(\ud{w})$. The open affine toric variety in $\bP(\ud{w})$ is defined by the cone $\sigma_{x_j},$
$$
\sigma_{x_j} = \sum_{ j'\neq j} \rho_{x_{j'}}\subset (N_x)_{\bR}.
$$
Let $u_{x_j}\in N_x$ be the minimal ray generator of $\rho_{x_j}$. We want our class group $\Cl(\bP(\ud{w}))$ to be torsion-free, so the cokernel of the following linear map
$$
M_x\to \bigoplus_{j=1}^n\bZ\rho_{x_j}, \quad m\mapsto \sum_{j=1}^n \langle m,u_{x_j}\rangle \rho_{x_j},
$$
is a torsion-free abelian group of rank 1. Specifically, we have an element $\rho_{x_0}$ in $\bigoplus_{j=1}^n\bZ\rho_{x_j}$:
$$
\rho_{x_0} = \sum_{j=1}^n \alpha_j \rho_{x_j},
$$
which generates $\Cl(\bP(\ud{w}))$ and satisfies $[\rho_{x_j}]= w_j [\rho_{x_0}]$ modulo $M_x$. We further assume that $\gcd(w_1,\ldots,w_n)=1$ so that we have integers $a_1,\ldots,a_n$ such that
$$
a_1 w_1 + \cdots + a_n w_n = 1.
$$

Similarly, for $\bP(\ud{d})$, we define the lattices $N_p$ and $M_p$, rays $\rho_{p_1},\ldots,\rho_{p_n}\subset (N_p)_{\bR}$, minimal ray generators $u_{p_j}\in N_p$, cones $\sigma_{p_j}\subset (N_p)_{\bR}$, the class group generator $\rho_{p_{0}}\in \bZ^{\Sigma_{\bP(\ud{d})}(1)}$ generating $\Cl(\bP(\ud{d}))$ with
$$
\rho_{p_0} = \sum_{k=1}^r \beta_k \rho_{p_k},
$$
and integers $b_1,\ldots,b_r$ such that
$$
b_1 d_1 + \cdots + b_r d_r = 1.
$$

Now we let $N=N_x\times N_p \times \bZ$. We will define a fan $\widetilde{\Sigma}$ in $N_{\bR}$. We first define ray generators.
\begin{align*}
u_0 &= ({0}_{N_x},{0}_{N_p},1), \\
u_{x_j}' &= (u_{x_j},{0}_{N_p},a_j), \qquad \textrm{for $j=1,\ldots,n$},\\
u_{p_k}' &= ({0}_{N_x},u_{p_k},b_k), \qquad \textrm{for $j=1,\ldots,n$}.
\end{align*}
For each $j=1,\ldots,n$ and $k=1,\ldots,r$, we define $\rho_0$, $\rho_{x_j}'$, $\rho_{p_k}'$ as rays generated by $u_0$, $u_{x_j}'$, $u_{p_k}'$, respectively. They constitute $\widetilde{\Sigma}(1)$. The maximal $(n+r-1)$-dimensional cone $\sigma_{jk}\subset N_\mathbf{R}$ is defined as follows:
$$
\sigma_{jk} = \bR_{\geq 0} u_0 + \left(\sum_{j' \neq j} \bR_{\geq 0} u_{x_j'}\right) + \left(\sum_{k' \neq k} \bR_{\geq 0} u_{p_k'} \right).
$$
Thus, there are $n+r+1$ rays. We define the fan $\widetilde{\Sigma}$ as
$$
\widetilde{\Sigma} = \bigcup_{\substack{1\leq j \leq n \\ 1\leq k \leq r}}\{\textrm{faces of $\sigma_{jk}$}\}.
$$
\begin{lemma}
The toric variety $X_{\widetilde{\Sigma}}$ defined by the fan $\widetilde{\Sigma}$ is isomorphic to $\widetilde{X}$:
$$X_{\widetilde{\Sigma}}=\widetilde{X}.$$
\label{fan}
\end{lemma}
\begin{proof}
We have the following identification,
$$
X_{\widetilde{\Sigma}} = \frac{\bC^{n+r+1}\setminus E_{\widetilde{\Sigma}}}{G_{\widetilde{\Sigma}}},
$$
where $E_{\widetilde{\Sigma}}$ is the exceptional set and $G_{\widetilde{\Sigma}} = \Hom_{\mathrm{Grp}}(\Cl(X_{\widetilde{\Sigma}}),\bC^*)$. The exceptional set is the intersection of zero sets of irrelevant ideals $(x^{\hat{\sigma}})$ for each maximal cone $\sigma$ generated by the monomial $x^{\hat{\sigma}}$,
$$
x^{\hat{\sigma}} = \prod_{\substack{\rho\in {\widetilde{\Sigma}}(1)\\ \rho\nprec\sigma} } x_\rho.
$$
The irrelevant ideal of $\sigma_{jk}$ is $(x_j p_k)$. Thus, the exceptional set is
$$
E_{\widetilde{\Sigma}} = (\bC\times \{0\}\times \bC^r) \cup (\bC\times \bC^n\times \{0\}),
$$
and we have
$$
X_{\widetilde{\Sigma}} = \frac{\bC\times (\bC^n\setminus\{0\})\times (\bC^r\setminus\{0\})}{G_{\widetilde{\Sigma}}}.
$$

We know that the class group is of rank $2$ and $G_{\widetilde{\Sigma}}\cong (\bC^*)^2$, since the torus is $(n+r-1)$-dimensional, and there are $n+r+1$ rays. We pick $\rho_{x_0}'$ and $\rho_{p_0}'$ for the generators of the class group,
\begin{align*}
\rho_{x_0}' = \sum_{j=1}^n \alpha_j \rho_{x_j}', \\
\rho_{p_0}' = \sum_{k=1}^r \beta_k \rho_{p_k}'.
\end{align*}
Then $M=M_x\times M_p\times \{0\}$, we have
\begin{align*}
[\rho_{x_j}'] = w_j [\rho_{x_0}'],\\
[\rho_{p_k}'] = d_k [\rho_{p_0}'],
\end{align*}
for each $j=1,\ldots,n$ and $k=1,\ldots,r$.

Finally, if we send $(0_{M_x},0_{M_p},1)$ along the map $M\to \bZ^{{\widetilde{\Sigma}}(1)}$, we have
$$
\rho_0 + \sum_{j=1}^n a_j \rho_{x_j}' + \sum_{k=1}^n b_k \rho_{p_j}'.
$$
However, modulo the module $M$, we have
\begin{align*}
&\rho_0 + \sum_{j=1}^n a_j \rho_{x_j}' + \sum_{k=1}^n b_k \rho_{p_j}' \mod M\\
&= \rho_0 + \sum_{j=1}^n a_j w_j \rho_{x_0}' + \sum_{k=1}^n b_k d_j \rho_{p_0}' \mod M\\
&= \rho_0 + \rho_{x_0}' + \rho_{p_0}'.
\end{align*}
Hence we find that the class group of $X_{\widetilde{\Sigma}}$ is a free abelian group of rank $2$, and the degree of $z$ is indeed $(-1,-1)$. The $G_{\widetilde{\Sigma}}$-action on $\bC\times (\bC^n\setminus\{0\})\times (\bC^r\setminus\{0\})$ is exactly the one defined in \eqref{groupaction}.
Thus we get the desired conclusion.
\end{proof}

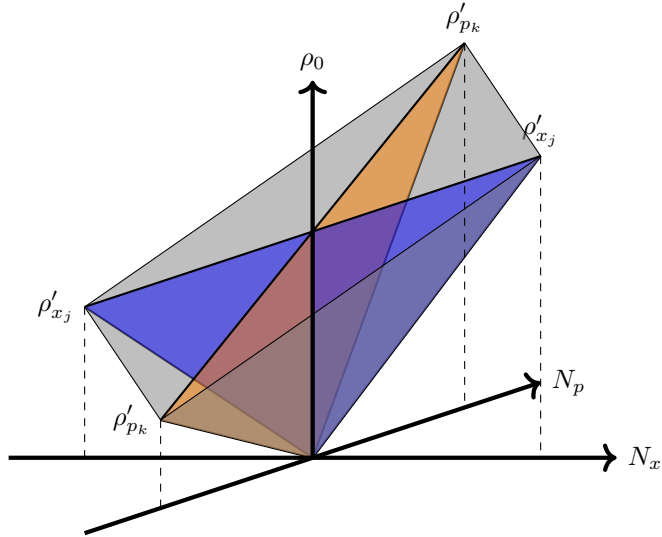
\begin{figure}[!ht]
\label{fan1}
\centering
\begin{tikzpicture}

\draw[thick] (0,0) -- (-3,2) node[anchor=east]{$\rho_{x_j}'$};
\draw[thick] (0,0) -- (-2,0.5) node[anchor=east]{$\rho_{p_k}'$};
\draw[thick] (0,0) -- (3,4) node[anchor=south]{$\rho_{x_j}'$};
\draw[thick] (0,0) -- (2,5.5) node[anchor=south]{$\rho_{p_k}'$};

\fill[gray,semitransparent] (0,0)--(-3,2)--(-2,0.5)--cycle;
\fill[gray,semitransparent] (0,0)--(3,4)--(2,5.5)--cycle;
\fill[gray,semitransparent] (0,0)--(2,5.5)--(-3,2)--cycle;

\fill[orange,semitransparent] (0,0)--(2,5.5)--(0,3)--cycle;
\fill[blue,semitransparent] (0,0)--(-3,2)--(3,4)--cycle;
\fill[orange,semitransparent] (0,0)--(0,3)--(-2,0.5)--cycle;

\fill[gray,semitransparent] (0,0)--(-2,0.5)--(3,4)--cycle;
\draw[thick] (-3,2)--(3,4) ;
\draw[thick] (2,5.5)--(-2,0.5) ;

\draw[ultra thick,->] (-4,0) -- (4,0) node[anchor=west]{$N_x$};
\draw[ultra thick,->] (-3,-1) -- (3,1) node[anchor=west]{$N_p$};
\draw[ultra thick,->] (0,0) -- (0,5) node[anchor=south]{$\rho_0$};

\draw[] (-3,2)--(-2,0.5) ;
\draw[] (-2,0.5)--(3,4) ;
\draw[] (3,4)--(2,5.5) ;
\draw[] (2,5.5)--(-3,2) ;
\draw[dashed](-3,2) -- (-3,0);
\draw[dashed](3,4) -- (3,0);
\draw[dashed](2,5.5) -- (2,0.666);
\draw[dashed](-2,0.5) -- (-2,-0.666);
\end{tikzpicture}
\caption{The Fan of $\widetilde{X}$}
\end{figure}

In Figure 1, the rays are tilted so that the variable $z$ corresponding to the ray $\rho_0$ has degree $(-1,-1)$. If we neglect the last $\bZ$ component in $N=N_x\times N_p \times \bZ$, the fan $\Sigma_{\bP(\ud{w})}$ along the $N_x$ direction and $\Sigma_{\bP(\ud{d})}$ along the $N_p$ direction are simply multiplied so that the base space becomes the product $\bP(\ud{w})\times \bP(\ud{d})$.

\begin{figure}[!ht]
\centering
\begin{tikzpicture}
\draw[thick] (0,0) -- (-3,2) node[anchor=east]{$\rho_{x_j}'$};
\draw[thick] (0,0) -- (-2,0.5) node[anchor=east]{$\rho_{p_k}'$};
\draw[thick] (0,0) -- (3,4) node[anchor=south]{$\rho_{x_j}'$};
\draw[thick] (0,0) -- (2,5.5) node[anchor=south]{$\rho_{p_k}'$};
\fill[gray,semitransparent] (0,0)--(-3,2)--(-2,0.5)--cycle;
\fill[gray,semitransparent] (0,0)--(3,4)--(2,5.5)--cycle;
\fill[gray,semitransparent] (0,0)--(2,5.5)--(-3,2)--cycle;
\fill[blue,semitransparent] (0,0)--(-3,2)--(3,4)--cycle;
\fill[gray,semitransparent] (0,0)--(-2,0.5)--(3,4)--cycle;
\draw[thick] (-3,2)--(3,4) ;
\draw[ultra thick,->] (-4,0) -- (4,0) node[anchor=west]{$N_x$};
\draw[ultra thick,->] (-3,-1) -- (3,1) node[anchor=west]{$N_p$};
\draw[] (-3,2)--(-2,0.5) ;
\draw[] (-2,0.5)--(3,4) ;
\draw[] (3,4)--(2,5.5) ;
\draw[] (2,5.5)--(-3,2) ;
\draw[dashed](-3,2) -- (-3,0);
\draw[dashed](3,4) -- (3,0);
\draw[dashed](2,5.5) -- (2,0.666);
\draw[dashed](-2,0.5) -- (-2,-0.666);
\end{tikzpicture}
\label{fan2}
\caption{After Collapsing the Orange Walls}
\end{figure}
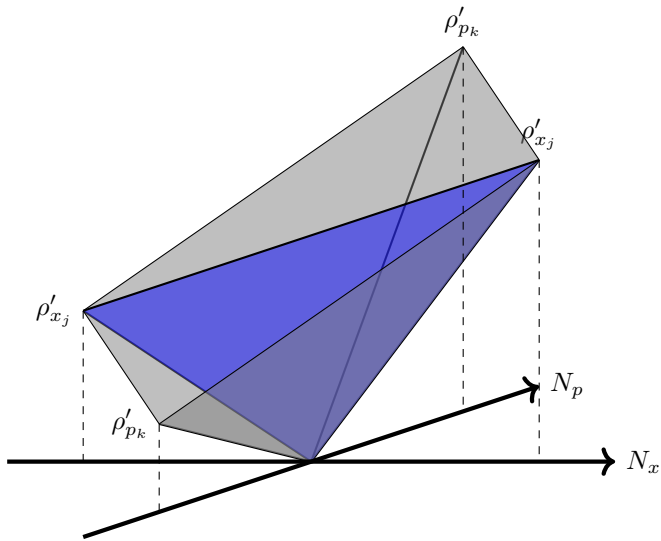

Now imagine collapsing the orange walls. The ray $\rho_0$ is now removed and some cones that were separated by an orange wall are merged. we can consider collapsing all the coordinates of $N_x$ and $\bZ$ onto the $N_p$-direction, which gives the surjection onto $\bP(\ud{d})$. In other words, this fan is giving a bundle over $\bP(\ud{d})$. We call this fan $\Sigma_{LG}$. The maximal cones of this fan are
$$
\sigma_k = \left(\sum_{j=1}^n \rho'_{x_j}\right) + \left(\sum_{\substack{1\leq k' \leq r\\k'\neq k}} \rho'_{p_k}\right).
$$

\begin{lemma}
$X_{LG} \cong X_{\Sigma_{LG}}$ and $X_{CY} \cong X_{\Sigma_{CY}}$. Thus, $\widetilde{X}$, $X_{CY}$, and $X_{LG}$ are all simplicial toric varieties.
\end{lemma}
\begin{proof}
The proof of $X_{CY} \cong X_{\Sigma_{CY}}$ is identical to that of $X_{LG} \cong X_{\Sigma_{LG}}$. We will only prove $X_{LG} \cong X_{\Sigma_{LG}}$. The variables $x_j$ are no longer contained in the generators of the irrelevant ideal of $\Sigma_{LG}$, since they are all contained in the maximal cone. Thus, the exceptional set is $p_1=\cdots=p_r=0$. On the other hand, since one ray is removed, the rank of the class group is also decreased by $1$. This time, if we send $(0_{M,x},0_{M,p},1)$ along the map $M\to \bZ^{\Sigma_{LG}(1)}$, we obtain
$$
(0_{M,x},0_{M,p},1)\mapsto \sum_{j=1}^n a_j \rho_{x_j}' + \sum_{k=1}^n b_k \rho_{p_j}',
$$
which is equivalent to
$$
\rho'_{x_0}+\rho_{p_0}'
$$
modulo $M$.
Thus, by fixing $\rho'_{x_0}$ as the generator of the class group, we get the charge grading $M\to \bZ$ which is defined in \eqref{definition:chandwt},
$$
\ch x_j = w_j,\, \quad \ch p_k = -d_k.
$$
Hence, the toric variety $X_{\Sigma_{LG}}$ is described as the following quotient,
$$
\frac{\bC^n\times (\bC^r\setminus\{0\})}{(x_1,\ldots,x_n,p_1,\ldots,p_r)\sim(\lambda^{w_1}x_1,\ldots,\lambda^{w_1}x_n,\lambda^{-d_1}p_1,\ldots,\lambda^{-d_r}p_r)}.
$$
The zero section $x_1=\cdots=x_n$ is exactly given by the blue wall. Thus, we get the desired conclusion.
\end{proof}

Similarly, we can define the fan $\Sigma_{CY}$ by collapsing the blue walls in Figure 1. By the same argument, we have the following proposition.
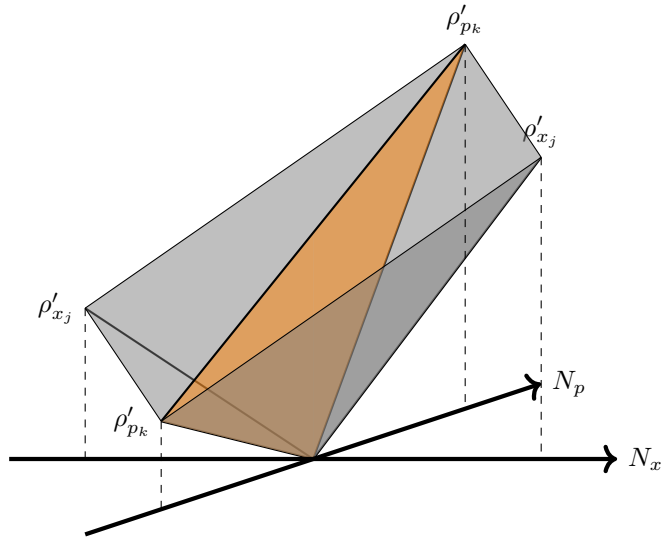
\begin{figure}[!ht]
\centering
\begin{tikzpicture}

\draw[thick] (0,0) -- (-3,2) node[anchor=east]{$\rho_{x_j}'$};
\draw[thick] (0,0) -- (-2,0.5) node[anchor=east]{$\rho_{p_k}'$};
\draw[thick] (0,0) -- (3,4) node[anchor=south]{$\rho_{x_j}'$};
\draw[thick] (0,0) -- (2,5.5) node[anchor=south]{$\rho_{p_k}'$};

\fill[gray,semitransparent] (0,0)--(-3,2)--(-2,0.5)--cycle;
\fill[gray,semitransparent] (0,0)--(3,4)--(2,5.5)--cycle;
\fill[gray,semitransparent] (0,0)--(2,5.5)--(-3,2)--cycle;

\fill[orange,semitransparent] (0,0)--(2,5.5)--(0,3)--cycle;
\fill[orange,semitransparent] (0,0)--(0,3)--(-2,0.5)--cycle;

\fill[gray,semitransparent] (0,0)--(-2,0.5)--(3,4)--cycle;
\draw[thick] (2,5.5)--(-2,0.5) ;

\draw[ultra thick,->] (-4,0) -- (4,0) node[anchor=west]{$N_x$};
\draw[ultra thick,->] (-3,-1) -- (3,1) node[anchor=west]{$N_p$};

\draw[] (-3,2)--(-2,0.5) ;
\draw[] (-2,0.5)--(3,4) ;
\draw[] (3,4)--(2,5.5) ;
\draw[] (2,5.5)--(-3,2) ;
\draw[dashed](-3,2) -- (-3,0);
\draw[dashed](3,4) -- (3,0);
\draw[dashed](2,5.5) -- (2,0.666);
\draw[dashed](-2,0.5) -- (-2,-0.666);
\end{tikzpicture}
\label{fan3}
\caption{After Collapsing the Blue Walls}
\end{figure}

\begin{proof}[Proof of Theorem \ref{blowup}]
We already have maps $\Sigma_{LG}\to \Sigma$ and $\Sigma_{CY}\to \Sigma$, and they correspond to $\phi_{LG}$ and $\phi_{CY}$ since they are identity maps on the torus $\bC^{n+r-1}$. We only need to prove that the star subdivisions of $\Sigma_{LG}$ and $\Sigma_{CY}$ are $\Sigma$.

The maximal cones of $\Sigma_{LG}$ is
$$
\sigma_k = \left( \sum_{j=1}^n \rho_{x_j}'\right) + \left( \sum_{k\neq k'}\rho'_{p_k}\right).
$$
The generator $u_0$ of the ray $\rho_0$ is contained in the interior of the cone
$$
\sum_{k=1}^r\rho'_{p_k}.
$$
Since for each $k'=1,\ldots,r$, $\left( \sum_{k\neq k'}\rho'_{p_k}\right)$ is a face of $\sum_{k=1}^r\rho'_{p_k}$, it never contains $u_0$. On the other hand, $\left( \sum_{j=1}^n \rho_{x_j}'\right)$ contains $u_0$. Thus, $\widetilde{\Sigma}$ is the weighted star subdivision of $\Sigma_{LG}$ along the ray $\rho_0$ which sits inside the interior of the cone $\sum_{k=1}^r\rho'_{p_k}\in \Sigma_{LG}$.

Similarly, $\widetilde{\Sigma}$ is the weighted star subdivision of $\Sigma_{CY}$ along the ray $\rho_0$ which sits inside the interior of the cone $\sum_{j=1}^n\rho'_{x_j}\in \Sigma_{CY}$.
\end{proof}


\begin{thebibliography}{99}


\bibitem{AS}
A. Adolphson and S.~I. Sperber, {\it On the Jacobian ring of a complete intersection,} J. Algebra {\bf 304} (2006), no.~2, 1193--1227.

\bibitem{Bah} E. Bahuaud et al., Well-posedness of nonlinear flows on manifolds of bounded geometry, Ann. Global Anal. Geom. {\bf 65} (2024), no.~4, Paper No. 25, 39 pp.; MR4742191


\bibitem{BK}
Barannikov, S.; Kontsevich, M.:
\newblock{\it Frobenius manifolds and formality of Lie algebras of polyvector fields},
\newblock Internat. Math. Res. Notices {\bf 4} (1998) 201--215.
\bibitem{Bara} Barannikov, Serguei: {\it Non-commutative periods and mirror symmetry in higher dimensions}, Comm. Math. Phys. 228 (2002), no. 2, 281--325.

\bibitem{CG}
Carlson, J.; Griffiths, P.:{\it Infinitesimal variations of Hodge structure and the global Torelli problem}, in Journees de geometrie algebrique, Angers, juillet 1979, Sijthoff and Noordhoff, Alphen aan den Rijn, 1980, 51--76.

\bibitem{CN} A. Chiodo; J. Nagel: {\it The hybrid Landau-Ginzburg models of Calabi-Yau complete intersections}, in {\it Topological recursion and its influence in analysis, geometry, and topology}, 103--117, Proc. Sympos. Pure Math., 100, Amer. Math. Soc., Providence, RI. 

\bibitem{Coevering} C. van~Coevering, Examples of asymptotically conical Ricci-flat K\"ahler manifolds, Math. Z. {\bf 267} (2011), no.~1-2, 465--496; MR2772262

\bibitem{ConHein2} R.~J. Conlon and H.-J. Hein, Asymptotically conical Calabi-Yau manifolds, I, Duke Math. J. {\bf 162} (2013), no.~15, 2855--2902; MR3161306

\bibitem{ConHein} R.~J. Conlon and H.-J. Hein, Classification of asymptotically conical Calabi-Yau manifolds, Duke Math. J. {\bf 173} (2024), no.~5, 947--1015; MR4740213


\bibitem{Cox} D.~A. Cox, J.~B. Little and H. Schenck, {\it Toric varieties}, Graduate Studies in Mathematics, 124, Amer. Math. Soc., Providence, RI, 2011

\bibitem{Dimca} A. Dimca, {\it Residues and cohomology of complete intersections}, Duke Math. J. {\bf 78} (1995), no.~1, 89--100; MR1328753

\bibitem{Dubrovin}
Dubrovin,  B. A.:{\it Geometry of $2$D topological field theories}, In  Integrable systems and quantum groups (Montecatini Terme, 1993),  
Lecture Notes in Math.\ vol.\ 1620, 120--348. Springer-Verlag, Berlin (1996).
arXiv:hep-th/9407018.


\bibitem{Ehr74}
P.~E.~Ehrlich,
\emph{Continuity properties of the injectivity radius function},
Compositio Mathematica \textbf{29} (1974), no.~2, 151--178.

\bibitem{JEi} J. Eichhorn, {\it Global analysis on open manifolds}, Nova Sci. Publ., New York, 2007.

\bibitem{JEl} J. Eldering, {\it Normally hyperbolic invariant manifolds}, Atlantis Studies in Dynamical Systems, 2, Atlantis Press, Paris, 2013.

\bibitem{RG} S. Gallot, D. Hulin and J. Lafontaine, {\it Riemannian geometry}, third edition, 
Universitext, Springer, Berlin, 2004; MR2088027

\bibitem{Her02} Hertling, C.: {\it Frobenius manifolds and moduli spaces for singularities.} Cambridge Tracts in Mathematics, 151. Cambridge University Press, Cambridge, 2002. x+270 pp. 

\bibitem{HKLR}
N.~J.~Hitchin, A.~Karlhede, U.~Lindstr{\"o}m, and M.~Ro{\v c}ek,
\emph{Hyperk\"ahler metrics and supersymmetry},
Commun. Math. Phys. \textbf{108} (1987), no.~4, 535--589.


\bibitem{HR} H. Hopf and W. Rinow, Ueber den Begriff der vollst\"andigen differentialgeometrischen Fl\"ache, Comment. Math. Helv. {\bf 3} (1931), no.~1, 209--225; MR1509435



\bibitem{Joyce}
D. Joyce: {\it Asymptotically locally Euclidean metrics with holonomy $SU(m)$}, Ann. Glob. Anal. Geom. 19 (1) (2001) 55--73.

\bibitem{Kerr} G. Kerr, {\it Weighted blowups and mirror symmetry for toric surfaces}, Adv. Math. {\bf 219} (2008), no.~1, 199--250; MR2435423

\bibitem{Kirwan} F.~C. Kirwan, {\it Cohomology of quotients in symplectic and algebraic geometry}, Mathematical Notes, 31, Princeton Univ. Press, Princeton, NJ, 1984; MR0766741


\bibitem{Konno} K. Konno, {\it On the variational Torelli problem for complete intersections}, Compositio Math. {\bf 78} (1991), no.~3, 271--296.


\bibitem{LW} Li, Si; Wen, Hao: \emph{On the $L^2$-Hodge theory of Landau--Ginzburg models}, Advances in Mathematics 396 (2022) 108165.


\bibitem{OV} Ogus, A, Vologodsky, V.:{\it Nonabelian Hodge theory in characteristic $p$}, Publ. Math. 106(1) (2007) 1--138.

\bibitem{MWM} J.~E. Marsden and A.~D. Weinstein, Reduction of symplectic manifolds with symmetry, Rep. Mathematical Phys. {\bf 5} (1974), no.~1, 121--130; MR0402819


\bibitem{Paper1}
Jeehoon Park; Jaewon Yoo, {\it $L^2$-Hodge Theoretic Construction of Frobenius Manifolds for Calabi-Yau Smooth Projective Hypersurfaces}, preprint.


\bibitem {PG}
Griffiths, Phillip A.: {\it On the periods of certain rational integrals. I, II,}
Ann. of Math. (2) 90 (1969), 460-495; ibid. (2) 90 (1969), 496-541.
\bibitem{Sa}
Sabbah, C.: {\it On a twisted de Rham complex}, Tohoku Math. J. 51 (1) (1998) 125--140.

\bibitem{Sab}
Sabbah, C.: {\it D\'eformations isomonodromiques et vari\'et\'es de Frobenius, une introduction.} Centre de Mathematiques, Ecole Polytechnique, U.M.R. 7640 du C.N.R.S., no. 2000-05, 251 pages.

\bibitem{Stacks} The stacks project authors: \emph{Stacks project}, \url{https://stacks.math.columbia.edu} (2018).

%




\bibitem{Saito} Saito, K.: \newblock{\it Primitive forms for an universal unfolding of a functions with isolated critical point}, \newblock Journ.\ Fac.\ Sci.\ Univ.\ Tokyo, Sect. IA Math.  {\bf 28} no.3 (1981) 777--792. 



\bibitem{SerreCartan} J.-P. Serre, {\it Cohomologie et fonctions de variables complexes,} in {\it S\'eminaire Bourbaki, Vol.\ 2}, Exp.\ No.\ 71, 213--218, Soc. Math. France, Paris.

\bibitem{SerreGAGA} J.-P. Serre, {\it G\'eom\'etrie alg\'ebrique et g\'eom\'etrie analytique,} Ann. Inst. Fourier (Grenoble) {\bf 6} (1955/56), 1--42.

\bibitem{Terasoma} T. Terasoma, Infinitesimal variation of Hodge structures and the weak global Torelli theorem for complete intersections, Ann. of Math. (2) {\bf 132} (1990), no.~2, 213--235.


\bibitem{TY90}
G.~Tian and S.-T.~Yau, Complete K\"ahler manifolds with zero Ricci curvature. I, \emph{J. Amer. Math. Soc.} \textbf{3} (1990), no.~3, 579--609.

\bibitem{TY91}
G.~Tian and S.-T.~Yau, Complete K\"ahler manifolds with zero Ricci curvature. II, \emph{Invent. Math.} \textbf{106} (1991), no.~1, 27--60.

\bibitem{Velez} A.-Q. V\'elez, {\it Mckay correspondence for Landau-Ginzburg models}, Communications in Number theory and Physics, Volume 3, Number 1, 173-208, 2009.

\bibitem{Voisin} C. Voisin, {\it Hodge theory and complex algebraic geometry. I}, translated from the French original by Leila Schneps, 
Cambridge Studies in Advanced Mathematics, 76, Cambridge Univ. Press, Cambridge, 2002.

\bibitem{Weibel} Weibel, Charles A.: {\it An introduction to homological algebra.} Cambridge Studies in Advanced Mathematics, 38. Cambridge University Press, Cambridge, 1994. xiv+450 pp. ISBN: 0-521-43500-5; 0-521-55987-1 

\bibitem{Yau} S.-T. Yau, Calabi's conjecture and some new results in algebraic geometry, Proc. Nat. Acad. Sci. U.S.A. {\bf 74} (1977), no.~5, 1798--1799.

\bibitem{Yau2} S.-T. Yau, On the Ricci curvature of a compact K\"ahler manifold and the complex Monge-Amp\`ere equation. I, Comm. Pure Appl. Math. {\bf 31} (1978), no.~3, 339--411; MR0480350














%






%

%
%







%


%




%




%
%
%


%
%




%


%




%
%


%

%

%
%


%


\end{thebibliography}
\end{document}